\pdfoutput=1
\documentclass[twoside,11pt]{article}
\usepackage{amssymb,amsfonts}
\usepackage{amsmath, mathtools, mathrsfs}
\usepackage{smile}
\usepackage[linesnumbered,ruled,vlined]{algorithm2e}
\usepackage{tikz}
\usetikzlibrary{patterns}
\usepackage{bbm}
\usetikzlibrary{shapes}
\usetikzlibrary{plotmarks}
\usepackage{xcolor}
\usepackage{capt-of}
\usepackage{fullpage}
\usepackage[colorlinks=true,
            linkcolor=blue,
            urlcolor=blue,
            citecolor=blue]{hyperref}
\usepackage[ruled]{algorithm2e}

\usepackage[colorinlistoftodos]{todonotes}
\usepackage{multirow}
\usepackage{caption}
\usepackage{perpage}

\newcommand{\npar}{\vspace{0.5em}\par }
\newcommand{\opar}{\vspace{0.3em}\par }

\newtheorem*{theorem*}{Theorem}
\newtheorem*{corollary*}{Corollary}
\newtheorem*{lemma*}{Lemma}
\newtheorem*{proposition*}{Proposition}
\newtheorem*{remark*}{Remark}

\theoremstyle{definition}
\newtheorem*{definition*}{Definition}

\SetKwInput{KwInput}{Input}
\SetKwInput{KwOutput}{Output}
\SetKwRepeat{KwRepeat}{repeat}{until}
\SetKw{Break}{break}
\SetKw{Continue}{continue}
\SetKw{Import}{import}

\MakePerPage{footnote}

\begin{document}

\begin{center}
    {\LARGE Robust Signal Detection with Quadratically Convex Orthosymmetric Constraints}

    {\large
        \begin{center}
            Yikun Li and Matey Neykov
        \end{center}}

    {Department of Statistics and Data Science\\ Northwestern University\\Evanston, IL 60208\\\texttt{yikunli2028@u.northwestern.edu}~~~~~~~~\texttt{mneykov@northwestern.edu}}
\end{center}

\begin{abstract}
This paper studies the problem of robust signal detection in Gaussian noise under quadratically convex orthosymmetric (QCO) constraints. We consider a minimax testing framework where the signal belongs to a QCO set and is separated from zero in Euclidean norm, while an adversary is allowed to arbitrarily corrupt a fraction $\epsilon $ of the samples. We establish the minimax separation radius between the null and alternative purely in terms of the constraint geometry, sample size, corruption rate, and noise scale. Our analysis argues that the Kolmogorov widths of the constraint set play a central role in determining the detection limits, paralleling to classic results in estimation problem. The derived lower bounds exhibit phase transitions with respect to the corruption rate and confirm that robust testing is statistically easier than robust estimation. While the information-theoretic upper bound is achieved by a computationally intractable test, we develop a polynomial-time algorithm that achieves the minimax lower bound up to logarithmic factors. Unlike prior work, our algorithm handles signals of arbitrary Euclidean length while respecting the QCO constraints. Finally, we extend these results to the robust $\ell _{p}$ norm testing for $1 \le p < 2$.
\end{abstract}

\tableofcontents

\newpage

\section{Introduction}

\subsection*{Data Contamination and Robust Statistics}

One of the most basic tasks in statistics is to perform inference about the unknown parameters used in data generation. However, a variety of reasons can lead the observed data to deviate from the original samples of the true mechanism. These reasons could range from measurement errors in some percentage of the observations, to even malicious tampering with the data. For example, in certain biological datasets, inconsistencies in the data can be natural but beyond the level of white noise \citep[see][e.g.]{marti2019recentrifuge, fibbi2023coclustering}, while in machine learning security applications, contamination can be artificial or malicious, which is also called \textit{data poisoning attacks} \citep[see][e.g.]{ramirez2022poisoningattacksdefensesartificial}; see also \cite{diakonikolas2019recentadvancesalgorithmichighdimensional} for other examples.

In this work, we take a stylized and classical statistical inference problem and study it from both constrained and robust perspectives. Specifically, we will study the Gaussian sequence model, where each coordinate $X _{ij}$ of a single sample $X _{i}$ is drawn independently from $\mathcal{N}\left(\mu _{j},\sigma ^{2}\right)$, and $X _{i}, i=1,2,\dots$ are i.i.d. Moreover, a constraint $K$ is imposed for the mean vector $\mu :=\left(\mu _{1},\mu _{2},\dots,\mu _{d}\right)$ s.t. $\mu \in K$, where $K$ is some appropriate set that belongs to $\mathbb{R}^{d}$ (or $K \subseteq \ell _{2}$ in the infinite-dimensional case). This extra information about the parameter encodes the prior knowledge of the statistician, and facilitates the inference meanwhile. Finally, and importantly, we also assume that there is a potential source of contamination, which is called the adversary. It is powerful in the sense that it has the knowledge of the true model and the access to the original data. For some positive constant $0 \le \epsilon <c _{0}<\frac{1}{2}$ known to both the statistician and the adversary, the adversary is allowed to inspect the original data and modify up to $\epsilon$ fraction of the samples arbitrarily. The statistician, however, does not know which samples have been modified. The goal of the statistician is to perform inference about $\mu$ based on the contaminated data. A more formal explanation of the the problem we study can be found in Section \ref{section: problem formulation and our contributions}.

Performing inference under constraints and even contamination could be a formidable task. Most of the related works focus on the estimation and testing problem under several simple constraints with no contamination (\cite{10.1214/aos/1176348890, 10.1145/3406325.3451084, 10.1145/3583680}). One crucial reason for this is that most of the methods and techniques heavily depend on the assumption of the true model behind the data (for example, likelihood ratio test and tests depending on the asymptotic distributions). Consider the following special case where $K=\mathbb{R}^{d}$ and $\epsilon =0$. Given $N$ samples $X _{1},X _{2},\dots,X _{N}\in \mathbb{R}^{d}$ from $\mathcal{N}(\mu ,\sigma ^{2}\mathbf{I}_d)$ with unknown $\mu $ and known $\sigma \in \mathbb{R}^{+}$. In the testing problem, we are asked to distinguish between two hypotheses:
\begin{equation}\label{eq: original testing problem}
    \begin{aligned}
        & H _{0}: \mu =0,\\ 
        & H _{1}: \left\lVert \mu \right\rVert_{2}^{}\ge \rho .
    \end{aligned}
\end{equation}
In the estimation problem, we are asked to give an estimate $\hat{\mu }$ of $\mu $ to achieve a small error under a certain loss function (a natural one here is the $\ell_2$ loss). $\bar{X}=\frac{1}{N}\sum\limits_{i=1}^{N}X _{i}$ is a sufficient statistic for the model. The testing and estimation based on $\bar{X}$ is widely known to be optimal for both problems. Specifically, to control errors within the level of (small) constant, one needs $\rho \gtrsim \frac{d ^{1/4}}{\sqrt{N}}\sigma $ for the testing problem (in the minimax sense) (\cite{ingster2003nonparametric, diakonikolas2022gaussian}) while $\inf\limits_{\hat \mu}\sup\limits_{\mu \in \mathbb{R}^d}\mathbb{E}\left\lVert\hat \mu - \mu \right\rVert_{2}^{2}\gtrsim \frac{d}{N}\sigma^2$ if one wants to optimally estimate $\mu$ in $\ell_{2}$ norm (\cite{10.1214/aos/1193342380}, also see Example 15.4 in \cite{wainwright2019high}). However, as mentioned, the property of $\bar{X}$ is very sensitive to the authenticity of the samples. The method quickly breaks down with even modification on a single sample (the adversary may modify any single datum, thereby forcing the mean to be zero and causing considerable type-II error.). This is precisely the major reason why testing and estimation problems are very challenging when studied from a robust perspective. 

Below, we provide an informal introduction to several key concepts that play a central role in our study --- namely, the strong contamination model, the Kolmogorov $k$-width, and quadratically convex orthosymmetric (QCO) sets. Throughout this paper, we assume that the constraint set $K$ is a QCO set, and we operate under the strong contamination model. As will be shown later, the Kolmogorov widths of $K$ determine the critical radius of the testing problem, i.e., the minimal $\rho $ under which the problem is testable in the minimax sense.

\subsection*{Strong Contamination Model}

Throughout the years, different types of adversary models have been proposed in the field of robust statistics. One of the earliest ones, the \textit{Huber contamination model}, was introduced by \cite{10.1214/aoms/1177703732}, where instead of directly sampling from the authentic distribution $\mathcal{P}$, it is assumed that there exists some error distribution $\mathcal{Q}$ selected by the adversary, and the data are sampled from the mixture $(1-\epsilon )\mathcal{P}+\epsilon \mathcal{Q}$. Huber's contamination model can be regarded as a case of the larger \textit{oblivious adversary}, where while the adversary has the adaptivity to some extent, the contamination process is guaranteed to be independent with the original samples. In contrast, the adversary in the \textit{strong contamination model} (\cite{diakonikolas2019recentadvancesalgorithmichighdimensional}) is granted with strong adaptivity, specifically, with knowledge of the original data, the true model and even the algorithms. The adversary is allowed to inspect and modify the original data based on its knowledge. In this model, the fraction of the contaminated samples must be controlled to obtain useful information. Note that this quantity is not bounded in the Huber contamination model (although the expected number of corruptions is $\epsilon N$). Other types of contamination are also studied, including the perturbation of the sampling distribution with bounded TV distance (which also includes Huber's contamination model). See \cite{diakonikolas2019recentadvancesalgorithmichighdimensional} for an overview.

\subsection*{Kolmogorov $k$-Width}

Introduced by \cite{kolmogoroff1936uber} in the study of function spaces and approximation theory, the Kolmogorov $k$-width (also known as the Kolmogorov $N$-width) characterizes how well an element can be reached by some linear approximations of intrinsic dimension $k$, given that the element belongs to a linear space and some restrictions are applied. For a given $k$, it selects the best $k$-dimensional projection operator under the measure of minimax error. Readers can refer to Definition \ref{def_Kolmogorov_width} for a formal presentation. Kolmogorov $k$-width is a classic concept and has seen usage in various fields. See \cite{Floater_2021} for a recent work in the problem of low-rank approximation of matrices using the Kolmogorov $k$-width. \cite{EVANS20091726} uses this concept to study the $k$-method in isogeometric analysis. Kolmogorov $k$-width is also popular in the problems related to partial differential equations (PDEs), where the decay rate of the Kolmogorov $k$-width is intensely studied (\cite{10.1093/imanum/dru066, GREIF2019216}). See also \cite{PAPAPICCO2022114687} for an application in neural networks and PDEs, and a very recent work, \cite{Arbes2025}, for an application in optimal transportation and PDEs. Finally, see \cite{GREIF2019216} and \cite{510617} for its power in the wave problem and wavelet representations. One can find more details and properties of the Kolmogorov $k$-width in \cite{pinkus2012n}.

\subsection*{QCO Sets}

Quadratically convex orthosymmetric (QCO) sets were first introduced in \cite{donoho1990minimax} (see also \cite[Chapter 4.8]{johnstone2019function}), where the authors give many examples of such sets, such as hyperrectangles, ellipsoids (both aligned with the coordinate axes) and sets of the form $\{\theta \left\lvert\right. \sum\limits_{i \in I}^{} a_i \psi(\theta_i^2) \le 1\}$, where $\psi$ is a given convex function and $a_i \in \mathbb{R}^{+},i=1,2,\dots$. The last example implies that all unit balls in (weighted) $\ell_p$ norms, for $p \ge 2$ are QCO sets. The set of QCO constraints is therefore much larger than simply the set of all ellipsoids. See Definition \ref{def_QCO_set} for the formal definition of a QCO set. In contrast to the present work, which focuses on testing, \cite{donoho1990minimax} characterized the minimax estimation rates in the corresponding Gaussian sequence model under an uncorrupted setting. It turns out that the Kolmogorov widths of the set $K$ (see Definition \ref{def_Kolmogorov_width}) determine the minimax rate of estimation. Later, by taking $\epsilon = 0$, we will observe that the Kolmogorov widths also determine the minimax rate of testing with QCO constraints, though the resulting rate is strictly smaller. Furthermore, we shall compare the minimax testing rate under QCO constrains with corruption to the corresponding minimax estimation rate (recently derived by \cite{prasadan2025informationtheoreticlimitsrobust}) and establish the same fundamental conclusion---testing is intrinsically easier than estimation.

\subsection{Problem Formulation and Our Contributions}\label{section: problem formulation and our contributions}

\begin{figure}[htbp]
    \centering
    \includegraphics[width=\textwidth]{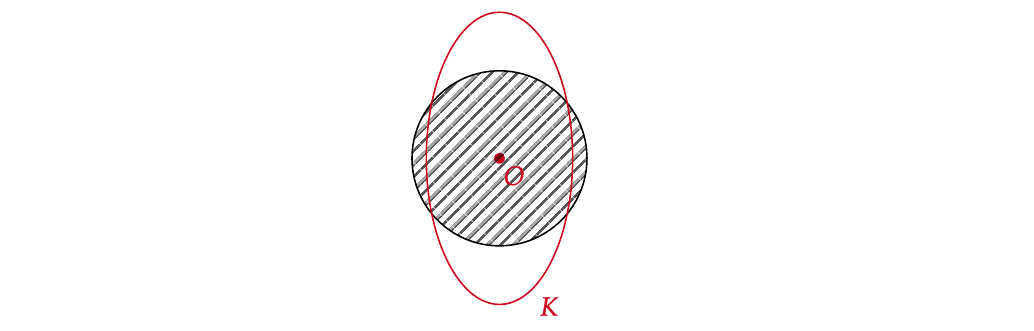}
    \caption{Testing $H_0: \mu = 0$ vs $H_1: \|\mu\|_2 \ge \rho, \mu \in K$ for the axis aligned ellipsoid (which is a special case of the QCO set) $K$ based on a single observation $X$. Upon applying a rotation to the data, one can use the test we develop in Section \ref{subsection: upper bounds} to solve the problems optimally for rotated (non-axis aligned) ellipsoids or more generally QCO sets (which are still centered around the origin). See our first remark in the \hyperref[section: discussion]{discussion} section for more details.}\label{figure: explanation of the studied problem}
\end{figure}

We now present the formal definition of the problem of interest. Specifically, we study the problem of minimax optimal signal detection in Gaussian noise (i.e., Gaussian mean testing) under the QCO constraints and strong $\epsilon $-contamination model in both the high-dimensional and the classic scenarios. Let $\tilde{\mathbf{X}}:=\left\{\tilde{X}_{1},\dots,\tilde{X}_{N}\right\}$ be the set of original i.i.d. samples, where $\tilde{X}_{i}:=\left(\tilde{X} _{i1},\tilde{X}_{i2},\dots,\tilde{X}_{id}\right)^{\top } $, $\tilde{X}_{ij}=\mu _{j}+\sigma \cdot \xi _{ij}$ and $\mu :=\left(\mu _{1},\mu _{2},\dots,\mu _{d}\right)^{\top }  \in K \subset \mathbb{R}^{d}$ for a known QCO set $K$ (i.e., the prior knowledge of the constraint). $\xi _{ij} \sim \mathcal{N}(0,1)$ are i.i.d. standard Gaussian noises. Instead of directly observing the original samples, an unknown adaptive adversary $\mathcal{C}$ from the strong $\epsilon $-contamination model (see Definition \ref{def_epsilon_contamination_model} for the formal definition) tampers with the data. The adversary swaps a fraction of the samples with arbitrary unknown values, where the contaminated fraction of the data is no more than a known, predetermined constant $\epsilon \in [0,c_0)$ for some small $c_0 < \frac{1}{2}$. The objective is to test the following null and alternative hypotheses with controlled type-I and type-II errors from the contaminated samples $\mathbf{X}=\left\{X _{i},\dots,X _{N}\right\}:=\mathcal{C}(\tilde{\mathbf{X}})$
\begin{equation}\label{eq: testing problem}
    \begin{aligned}
        & H _{0}: \mu =0,\\ 
        & H _{1}: \left\lVert \mu \right\rVert_{2}^{}\ge \rho ,\mu \in K.
    \end{aligned}
\end{equation}
Intuitively, if $\left\lVert \mu \right\rVert_{2}^{}$ is overly close to zero or the fraction of the contaminated samples is overly large, we are not able to detect the difference in the worst case. Therefore, we are interested in characterizing the asymptotically smallest value of $\rho$ (depending on other parameters) for which the testing problem becomes feasible. In detail, for a given desired type-I and type-II errors $\alpha $, define the set of acceptable tests under $H _{0}$ as,
\begin{equation}\label{eq: def_A_s}
    A _{s}(N,d,K,\epsilon ,\alpha ,\sigma ):=\left\{\phi :\sup\limits_{\mathcal{C}}\mathbb{P}_{0}\left(\phi (\mathcal{C}(\tilde{\mathbf{X}}))=1 \right)\le \alpha\right\},
\end{equation}
where $\mathbb{P}_{0}$ denotes the law of $\tilde{\mathbf{X}}$ when $\mu =\mathbf{0}$, and similar for $\mathbb{P}_{\mu }$. We write $A _{s}$ and omit the dependency hereafter when the context is clear. Let 
\begin{equation}\label{eq: def_critical_rho}
    \rho (N,d,K,\epsilon ,\alpha ,\sigma ):=\inf\limits_{\rho }\left\{\rho :\rho >0,\sup\limits_{\left\lVert \mu \right\rVert_{2}^{}\ge \rho ,\mu \in K}\sup\limits_{\mathcal{C}}\mathbb{P}_{\mu }\left(\phi (\mathcal{C}(\mathbf{X}))=0\right)\le \alpha, \phi \in A _{s} \right\}
\end{equation}
be the smallest possible $\rho $ such that the testing task with both errors less than $\alpha $ is doable. We also write it as $\rho _{\text{critical}}$ for convenience when the context is clear. We are going to study the asymptotic behaviour of $\rho (N,d,K,\epsilon ,\alpha ,\sigma )$. See Figure \ref{figure: explanation of the studied problem} for an illustrative example of the problem.

Our main result shows that $\rho _{\text{critical}}$ is determined by three independently valid lower bounds. Moreover, the lower bounds are tight and optimal with only the exception of logarithmic terms. The \hyperref[theorem: lower bound 1]{first lower bound} captures the geometric nature of the QCO set $K$ via the Kolmogorov $k$-width. This lower bound is independent of the adversary, therefore also works in the uncorrupted setting, i.e., when $\epsilon =0$. The \hyperref[theorem: lower bound 2]{second lower bound} is obtained within the framework of \textit{differential privacy}, which is introduced in Section \ref{subsection: main lower bound 2}. The \hyperref[theorem: lower bound 3]{third lower bound} is derived from a combined argument related to the identical mixture distributions, computation of the total variation and a designed strategy of the adversary. For the upper bounds, we present two algorithms which conclude whether $H _{0}$ holds or not while ensuring type-I and type-II errors both less than $\alpha $. The requirements on $\lVert \mu \rVert_{2}$ imposed by the \hyperref[algorithm: exponential algorithm]{theoretical algorithm} are slightly closer to the lower bounds. However, this improvement comes at the cost of exponential time complexity, rendering the method computationally intractable for large $N$ and $d$. The \hyperref[algorithm: polynomial algorithm]{polynomial-time algorithm}, in contrast, provides a slightly weaker upper bound than the \hyperref[algorithm: exponential algorithm]{theoretical algorithm} (though still nearly optimal) while achieving polynomial-time complexity. Nevertheless, both upper bounds suggest that the three lower bounds previously discussed suffice to fully characterize the behavior of $\rho_{\text{critical}}$. Moreover, such minimax optimal rates also imply a phase transition phenomenon regarding $\epsilon $ for the contamination process. Finally, the results and algorithms derived in the $\ell _{2}$ norm scenario can be naturally extended to the $\ell _{p}$ case for $1 \le p<2$. Section \ref{section: discussion} contains such relevant discussions.

\subsection{Related Literature}

In this section, we discuss some relevant works. We first comment on the literature on minimax testing with uncorrupted data. In their seminal body of work on nonparametric testing, \cite{ingster1982minimax,ingster1993asymptotically1,ingster1993asymptotically2, ingster1993asymptotically3,ingster2000minimax} showed the difference between testing and estimation in multiple problems. See also \cite{ingster2003nonparametric} for a summary of the many existing results in the area. Despite various progress in \textit{goodness-of-fit} testing in Gaussian models, we are not aware of a general result dedicated to QCO sets. For certain special types of ellipsoids, the problem was first studied by \cite{ermakov1991minimax}, where the author derived exact minimax rates. \cite{suslina1996minimax} studied the asymptotics of the problem of signal detection in $\ell_q$ ellipsoids with $\ell_p$ ball removed. Since $\ell_q$ ellipsoid for $q \ge 2$ are QCO sets, when $p=2$ we provide a generalization of these results. For general ellipsoids, \cite{baraud2002non} studied the testing problem in $\mathbb{R}^n$ and $\ell_2$ and obtained the minimax optimal critical radius. Our work can hence be seen as a generalization to part of \cite{baraud2002non}. In an interesting recent work, \cite{wei2020local} studied local testing rates on ellipsoids. Here local means that the point of interest under the null is not necessarily the zero point, but can be an arbitrary point in the ellipsoid. \cite{wei2020local} proved upper and lower bounds based on localized Kolmogorov widths, which are not guaranteed to match in general, but do match for certain examples with Sobolev ellipsoids. In contrast, the present work focuses on more general sets, but the point of interest under the null hypothesis is always anchored at the zero point, and the critical radius is determined based on the global (instead of local) Kolmogorov widths of the set.

Subsequently, we note that the recent literature on robust testing is relatively sparse. A notable exception is the paper \cite{10353143}, which studies a special case of the problem we consider, where $K = \mathbb{R}^{d}$, but considers a broader class of adversaries. The author also developed a polynomial-time algorithm, the framework of which we adapt to develop a more general algorithm suitable for our settings. In addition, \cite{diakonikolas2017statistical} develops a statistical query (SQ) lower bound for robust estimation problem in high-dimensional setting and extends the results to the testing problem.

\subsection{Organization}

This paper is structured as follows. Section \ref{section: preliminary} introduces the notations and several important concepts used through the paper. Section \ref{section: main result} presents the lower and upper bounds for the signal detection problem and contains the main results of the paper. It also contains the discussion of the phase transitions of the problem. Section \ref{section: extension to lp testing problem} is the extension of the main results in Section \ref{section: main result} to the robust $\ell _{p}$ norm testing problem. Section \ref{section: experiments} presents the experiment procedures, including the data generation, constraint and adversary setups, the experimental results of the polynomial-time algorithms related to Theorem \ref{theorem: upper bound 2} and relevant discussion. Section \ref{section: discussion} contains the discussion which compares the minimax rates of robust testing to robust estimation. It also includes some open problems for future works. The pseudo code of the relevant algorithms discussed in the paper and the additional experimental results of Section \ref{section: experiments} are provided in Appendix \ref{section: algorithms}. We attach necessary probabilistic tools and auxiliary lemmas in Appendix \ref{section: concepts and lemmas}. The proofs of the lower and upper bounds and other main theorems, lemmas, and corollaries in Section \ref{section: main result}, Section \ref{section: extension to lp testing problem} and Section \ref{section: discussion} are included in Appendix \ref{section: proofs of the main results}.

\section{Preliminary}\label{section: preliminary}

\subsection{Notations and Definitions}\label{subsection: notations and definitions}

To simplify the analysis, we assume that $K \subset \mathcal{X}=\mathbb{R}^{d}$, and therefore $\mu \in \mathbb{R}^{d}, X _{i}\in \mathbb{R}^{d}$ through the paper. However, one can extend the results to the case when $K \subset \ell _{2}$ with the same techniques and algorithms with minimal efforts given the assumption $\lim\limits_{k\rightarrow \infty }D _{k}(K)=0$. See Definition \ref{def_Kolmogorov_width} for further details. We use $\rho $ to denote the value of the $\ell _{2}$ norm of $\mu $ under $H _{1}$. $\left\lVert \cdot \right\rVert_{p}^{}$ means the $\ell  _{p}$ norm in $\mathbb{R}^{d}$ for $1 \le p \le \infty $. When $A$ is a set, $\left|A\right|$ denotes the cardinality of $A$. We use calligraphic font $\mathcal{P},\mathcal{Q},\mathcal{R}\dots$ to represent probability distributions. $\mathcal{P}\ll \mathcal{Q}$ means that $\mathcal{P}$ is absolutely continuous with respect to $\mathcal{Q}$. $\text{TV}(\mathcal{P},\mathcal{Q})$ denotes the total variation distance and $\text{KL}(\mathcal{P},\mathcal{Q})$ for the Kullback-Leibler divergence if $\mathcal{P}\ll \mathcal{Q}$. When the underlying probability distribution is not explicitly specified, we use $\mathbb{P}$ alternatively to denote its probability law in general. We use $\left\langle \cdot , \cdot  \right\rangle$ and $(\cdot )^{\top } (\cdot )$ alternatively (possibly with a subscript when needed) to represent the inner product in a Hilbert space. $\text{diag}\{a\}$ represents the diagonal matrix with diagonal entries $a _{i},i=1,2,\dots$. The convenient notation $[n]$ denotes the set $\left\{1,2,\dots,n\right\}$. When $I$ is an index set, $v_{I}$ represents the restriction of $v$ on $I$ if $v$ is a vector and $\mathbf{Y}_{I}$ represents the sub-matrix of $\mathbf{Y}$ with rows selected according to $I$ if $\mathbf{Y}$ is a matrix. $\mathcal{O}$, $\varTheta $, $\varOmega $ notations are used to denote the asymptotic order of rates. The lowercase $c _{1},c _{2},\dots$ denote the universal constants and their specific values may differ from line to line.

In addition to these plain notations, we let $\tilde{\mathbf{X}}:=\left\{\tilde{X}_{1},\dots,\tilde{X}_{N}\right\}$ represent the original samples, $\mathbf{X}:=\mathcal{C}(\tilde{\mathbf{X}})=\left\{X _{1},\dots,X _{N}\right\}$ represent the data contaminated by an adversary  $\mathcal{C}$. Let the uppercase $C$ denote the unknown set of indices corresponding to the contaminated samples. By this notation, we have $\tilde{\mathbf{X}}_{[N]\backslash C}=\mathbf{X}_{[N]\backslash C}$. Moreover, the following definitions formally introduce the concepts of QCO set and the Kolmogorov $k$-width used in this paper, which play crucial roles for our results.

\begin{definition}[Quadratically convex orthosymmetric (QCO) set]\label{def_QCO_set}
    Given a set $K \subset \mathbb{R}^{d}$, we say $K$ is a quadratically convex orthosymmetric (QCO) set if it satisfies the following conditions:\\ 
\hspace*{0.5em}(1), $K$ is convex;\\ 
\hspace*{0.5em}(2), $K$ is quadratically convex, which means that $K ^{2}$ is also convex, where $K ^{2}$ is defined as
\begin{equation*}
    K ^{2}:=\left\{\left(\theta _{1}^{2},\dots,\theta _{d}^{2}\right)^{\top } \left\lvert\right. \left(\theta _{1},\dots,\theta _{d}\right)^{\top } \in K\right\};
\end{equation*} 
\hspace*{0.5em}(3), $K$ is orthosymmetric, which means that if $\theta =\left(\theta _{1},\dots,\theta _{d}\right)^{\top } \in K$, then $\theta _{\eta }:=\left(\eta _{1}\theta _{1},\dots,\eta _{d}\theta _{d}\right)^{\top }\in K$, where $\eta _{i}\in \left\{-1,1\right\}, 1 \le i \le d$.
\end{definition}

\begin{remark*}
    The definition of QCO set can be naturally extended when $K \subset \ell _{2}$, as we only need to replace $[d]$ by the positive integers $\mathbb{N}^{+}$.
\end{remark*}

\begin{definition}[Kolmogorov $k$-width, also known as Kolmogorov $N$-width]\label{def_Kolmogorov_width}
    Let $\mathcal{X}$ be a Banach space equipped with the norm $\left\lVert \cdot \right\rVert_{}^{}$, and $K \subset \mathcal{X}$ is a subset. The Kolmogorov $k$-width is defined as 
\begin{equation}\label{eq: definition of Kolmogorov width}
    D _{k}(K) = \inf_{P \in \mathcal{P}_{k}} \sup_{\theta \in K} \left\lVert \theta -P \theta \right\rVert_{}^{},
\end{equation}

where $\mathcal{P}_{k}$ is the set of all projection operators that project a vector onto some subspace of $\mathcal{X}$ with intrinsic dimension $k$.
\end{definition}

\begin{remark*}
    In this work, our analysis will be confined to orthogonal projections. Moreover, though $\mathcal{P}_{k}$ is the set of all orthogonal projection with intrinsic dimension $k$, we discretize this set --- only consider the projections that are aligned with the coordinate axes without loss of generality, as shown by the analysis later. For now, let us first assume that the norm leveraged in the definition is the $\ell _{2}$ norm.
\end{remark*}

Based on the definition, one can derive several basic properties of the Kolmogorov widths. They are non-increasing with respect to $k$, meaning that $D _{k _{1}}(K)\ge D _{k _{2}}(K)$ for any $0 \le k _{1}\le k _{2}$. In addition, for a finite-dimensional space $\mathcal{X}$ with full dimension $d$, we always have $D _{0}(K)=\sup\limits_{\theta \in K}\left\lVert \theta \right\rVert_{2}^{}$ and $D _{d}(K)=0$ for any $K \subset \mathcal{X}$. For the more general case where $K \subset \ell _{2}$, we shall assume $\lim\limits_{d \rightarrow \infty} D _{d}(K)=0$. In this paper, we further define $D _{k}(K)=\infty $ when $k<0$ and $D _{k}(K)=0$ for $k>d$ for the convenience of the marginal cases.

The connection between the Kolmogorov widths and our main results is established through the following definition of the first and second dimensions of the optimal projections, which are designed to identify the optimal $k$ satisfying certain marginal conditions determined by the geometric properties of the constraint set $K$ and the corruption process.

\begin{definition}[First optimal dimension and projection]\label{def_first_optimal_dimension}
    Given the noise scale $\sigma $, the number of samples $N$ and the prior knowledge $K$, the first optimal dimension is an integer $k _{1}^{\star }(K,\sigma ,N)\in \left[0,d\right]$ depending on $K$, $\sigma $, and $N$ where it satisfies that $D _{k _{1}^{\star }-1}(K)>\frac{(k _{1}^{\star })^{1/4}}{\sqrt{N}}\sigma $ but $D _{k _{1}^{\star }}(K)\le \frac{(k _{1}^{\star }+1)^{1/4}}{\sqrt{N}}\sigma $. The projection operator corresponding to $k _{1}^{\star }(K,\sigma ,N)$ is denoted as $P _{1}^{\star }(K,\sigma ,N)$.
\end{definition}

\begin{definition}[Second optimal dimension and projection]\label{def_second_optimal_dimension}
    Given the noise scale $\sigma $, the number of samples $N$, the prior knowledge $K$ and the corruption rate $\epsilon $, the second optimal dimension is an integer $k _{2}^{\star }(K,\sigma ,N,\epsilon )\in [0,d]$ depending on $K$, $\sigma $, $N$, and $\epsilon $ where it satisfies that $D _{k _{2}^{\star }-1}(K)>\frac{(k _{2}^{\star })^{1/4}\sqrt{\epsilon }}{N ^{1/4}}\sigma $ but $D _{k _{2}^{\star }}(K)\le \frac{(k _{2}^{\star }+1)^{1/4}\sqrt{\epsilon }}{N ^{1/4}}\sigma $. The projection operator corresponding to $k _{2}^{\star }(K,\sigma ,N,\epsilon )$ is denoted as $P _{2}^{\star }(K,\sigma ,N,\epsilon )$.
\end{definition}

We note that from the monotone property of the Kolmogorov widths and the boundary values, such $k _{1}^{\star }(K,\sigma ,N)$ and $k _{2}^{\star }(K,\sigma ,N,\epsilon )$ must exist and are unique. In fact, it could be equally defined as 
\begin{equation*}
    \begin{aligned}
        k _{1}^{\star }(K,\sigma ,N)&:=\max\limits \left\{j \left\lvert\right. 0 \le j \le d, D _{j-1}(K)>\frac{j ^{\frac{1}{4}}}{\sqrt{N}}\sigma \right\},\\ 
        k _{2}^{\star }(K,\sigma ,N,\epsilon )&:=\max\limits \left\{j \left\lvert\right. 0 \le j \le d, D _{j-1}(K)>\frac{j ^{\frac{1}{4}}\sqrt{\epsilon }}{N ^{\frac{1}{4}}}\sigma \right\},
    \end{aligned}
\end{equation*}
or
\begin{equation*}
    \begin{aligned}
        k _{1}^{\star }(K,\sigma ,N)&:=\min\limits \left\{j \left\lvert\right. 0 \le j \le d, D _{j}(K)\le \frac{(j+1)^{\frac{1}{4}}}{\sqrt{N}}\sigma \right\},\\ 
        k _{2}^{\star }(K,\sigma ,N,\epsilon )&:=\min\limits \left\{j \left\lvert\right. 0 \le j \le d, D _{j}(K)\le \frac{(j+1) ^{\frac{1}{4}}\sqrt{\epsilon }}{N ^{\frac{1}{4}}}\sigma \right\}.
    \end{aligned}
\end{equation*}

In the general case where $K \subset \ell _{2}$, the requirement that $\lim\limits_{d\rightarrow \infty }D _{d}(K)=0$ ensures $k _{1}^{\star }(K,\sigma ,N)<\infty $,\footnote{$k _{2}^{\star }(K,\sigma ,N,\epsilon )$ is possibly infinite. This is the case for example if $\epsilon =0, D _{k}(K)>0, \forall k \in \mathbb{N}$.} which is important for the upper bounds and the implementation of the algorithms discussed later. When the context is clear, we omit the dependency with $K$, $\sigma $, $N$, and $\epsilon $, only writing them as $k _{1}^{\star }, k _{2}^{\star }$ and $P _{1}^{\star },P _{2}^{\star }$.

Below are several concrete examples of the concepts introduced above, which serve to illustrate the definitions and provide motivation from practical applications.

\begin{example}[$\mathbb{R}^{d}$]\label{example: QCO R^d}
    The whole space $\mathbb{R}^{d}$ is trivially convex, quadratically convex and orthosymmetric. Hence, it is a QCO set. The Kolmogorov $k$-width can be easily calculated as
    \begin{equation*}
        D _{k}(\mathbb{R}^{d})= \left\{
        \begin{tabular}{ll}
            $\infty $, & $(0 \le k \le d-1)$,\\ 
            $0$, & $(k=d)$.
        \end{tabular}
        \right.
    \end{equation*} 
    The Kolmogorov $k$-width is infinite until $k$ reaches $d$, which implies  one cannot uniformly approximate a vector in $\mathbb{R}^{d}$ by a linear subspace well without any additional constraint. The first and second optimal dimensions are $k _{1}^{\star }=k _{2}^{\star }=d, \forall \epsilon \in [0,c _{0})$.
\end{example}

\begin{example}[Hyperrectangles in $\mathbb{R}^{d}$]\label{example: QCO hr}
    The hyperrectangle in $\mathbb{R}^{d}$ is defined as
    \begin{equation*}
        H _{d}(a _{1},\dots,a _{d}):=\left\{x \left\lvert\right.  x=(x _{1},\dots,x _{d}) \in \mathbb{R}^{d},\left|x_{i}\right|\le a _{i},a _{i}\ge 0, 1 \le i \le d\right\}.
    \end{equation*}
    Obviously, $H _{d}(a _{1},\dots,a _{d})$ is convex and orthosymmetric. One can verify that $H _{d}^{2}(a _{1},\dots,a _{d})=\left\{x ^{2}\left\lvert\right. x \in H _{d}(a _{1},\dots,a _{d})\right\}$ is also a hyperrectangle with each side being $[0,a _{i}^{2}], 1 \le i \le d$. Therefore, it is a QCO set. By reindexing if necessary, we may assume that the sequence $\left\{a _{i}\right\}_{i=1}^{d}$ is arranged in decreasing order. The Kolmogorov widths are 
    \begin{equation*}
        D _{k}(H _{d}(a _{1},\dots,a _{d}))=\left\{
        \begin{tabular}{ll}
            $\sqrt{\sum\limits_{i=k+1}^{d}a _{i}^{2}}$, & $(0 \le k \le d-1)$,\\ 
            $0$, & $(k=d)$.
        \end{tabular}
        \right.
    \end{equation*}
    The decay rate of the Kolmogorov widths depends on the specific choice of $a _{i}$. The first optimal dimension is $k _{1}^{\star }=\max\limits \left\{j \left\lvert\right. 1 \le j \le d,\sqrt{\sum\limits_{i=j}^{d}a _{i}^{2}}>\frac{j ^{1/4}}{\sqrt{N}}\sigma \right\}$ when $\sqrt{\sum\limits_{i=1}^{d}a _{i}^{2}}>\frac{1}{\sqrt{N}}\sigma $ and otherwise $0$. The second optimal dimension is $k _{2}^{\star }=\max\limits \left\{j \left\lvert\right. 1 \le j \le d, \sqrt{\sum\limits_{i=j}^{d}a _{i}^{2}}>\frac{j ^{1/4}\sqrt{\epsilon }}{N ^{1/4}}\sigma \right\}$ when $\sqrt{\sum\limits_{i=1}^{d}a _{i}^{2}}>\frac{\sqrt{\epsilon }}{N ^{1/4}}\sigma $ otherwise $0$.
\end{example}

\begin{example}[Ellipses in $\mathbb{R}^{d}$]\label{example: QCO el}
    An ellipse in $\mathbb{R}^{d}$ that is centered at the origin and aligns with the axes can be represented as
    \begin{equation*}
        E _{d}(a _{1},\dots,a _{d})=\left\{x \left\lvert\right. x \in \mathbb{R}^{d},\sum\limits_{i=1}^{d}\frac{x _{i}^{2}}{a _{i}}\le 1\right\}.
    \end{equation*}
    Any such ellipse is convex and orthosymmetric. The square of the ellipse $E _{d}(a _{1},\dots,a _{d})$ is
    \begin{equation*}
        E _{d}^{2}(a _{1},\dots,a _{d})=\left\{x \left\lvert\right. x \in \mathbb{R}_{+}^{d}, \sum\limits_{i=1}^{d}\frac{x _{i}}{a _{i}}\le 1\right\},
    \end{equation*}
    which is part of a simplex (enclosed by several hyperplanes) in $\mathbb{R}^{d}$, and therefore convex. Figure \ref{fig: qco_ellipse} presents a depiction of $E _{d}^{2}(a _{1},\dots,a _{d})$ when $d = 2$.
    \begin{figure}[htbp]
        \centering
        \includegraphics[scale=0.4]{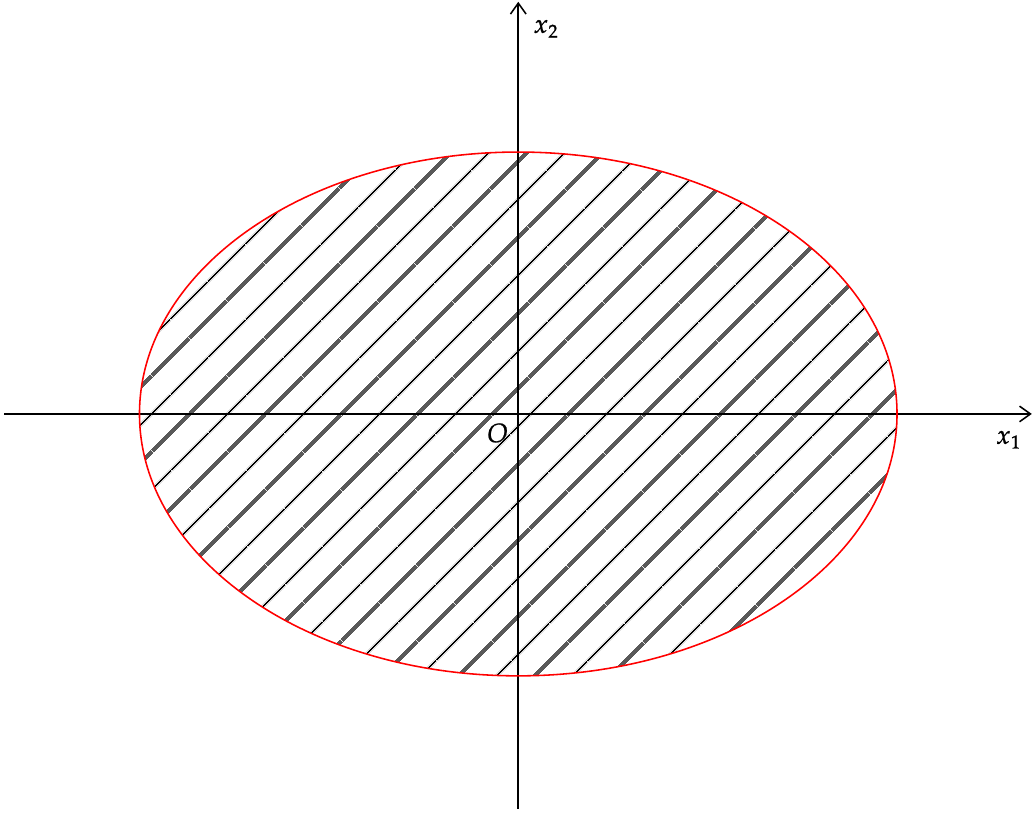}
        \includegraphics[scale=0.4]{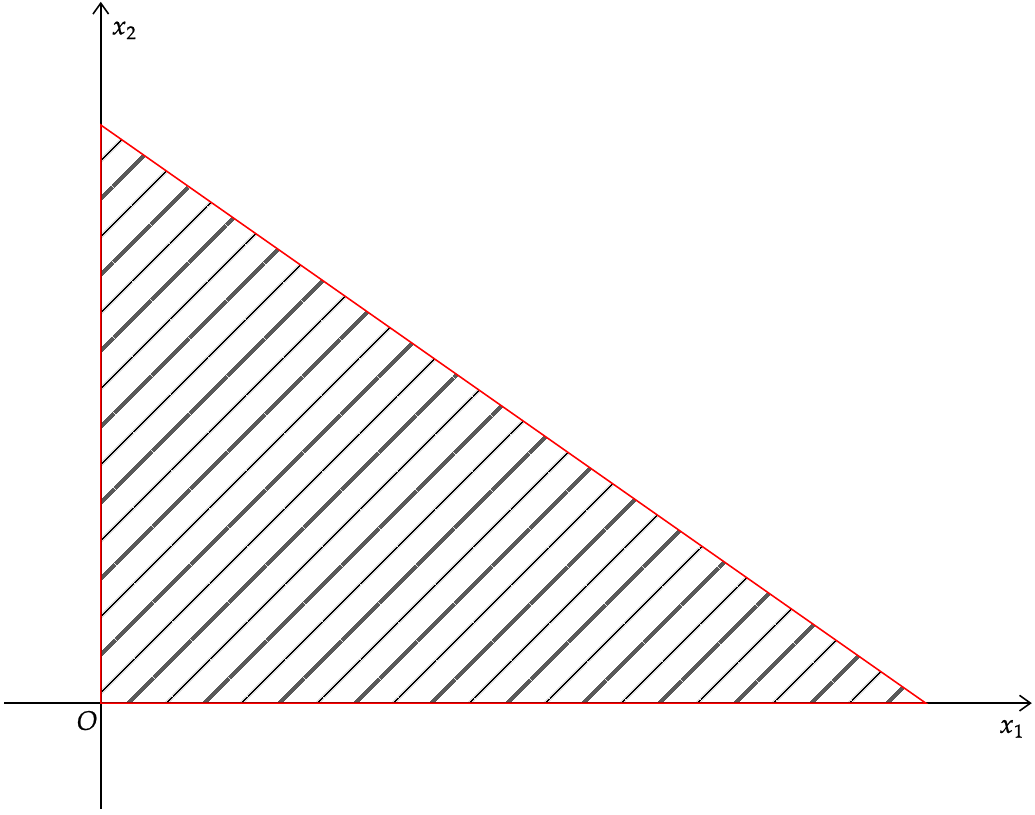}
        \caption{An example of an ellipse and its square-set in $\mathbb{R}^{2}$. The square-set is a triangle in $\mathbb{R}^{2}$, which is convex.}\label{fig: qco_ellipse}
    \end{figure}
    Again, we assume $\left\{a _{i}\right\}_{i=1}^{d}$ is arranged in decreasing order, by which the Kolmogorov widths are
    \begin{equation*}
        D _{k}(E _{d}(a _{1},\dots,a _{d}))=\left\{
        \begin{tabular}{ll}
            $\sqrt{a _{k+1}}$, & $(0 \le k \le d-1)$,\\ 
            $0$, & $(k=d)$.
        \end{tabular}
        \right.
    \end{equation*}
    The first optimal dimension is $k _{1}^{\star }=\max\limits \left\{j \left\lvert\right. 1 \le j \le d, a _{j}>\frac{\sqrt{j}}{N}\sigma ^{2}\right\}$ when $a _{1}>\frac{1}{N}\sigma ^{2}$ and $0$ otherwise. The second optimal dimension $k _{2}^{\star }=\max\limits \left\{j \left\lvert\right. 1 \le j \le d, a _{j}>\frac{\sqrt{j}\epsilon }{\sqrt{N}}\sigma ^{2}\right\}$ when $a _{1}>\frac{\epsilon }{\sqrt{N}}\sigma ^{2}$ and $0$ otherwise.

    When $a _{i}=i ^{-2s}$ for some $s>1$, the corresponding ellipse in $\mathbb{R}^{d}$ can be regarded as a special case of the more general ellipses of Sobolev types, which is defined as $S(s):=\left\{x \left\lvert\right. x \in \mathbb{R}^{\infty },\sum\limits_{i=1}^{\infty }i ^{2s}x _{i}^{2}\le 1\right\}$. This type of sets is closely related to the estimation and density approximation problem with the class of $s $-times differentiable functions in probability space. We refer readers to \cite{wainwright2019high} and \cite{tsybakov2009introduction} for more details.
\end{example}

\begin{remark*}
    When an ellipse is not axes-aligned: (1), the orthosymmetric property is not guaranteed, which can be easily checked; (2), the quadratically convex property is true for $d=2$ and false for $d \ge 3$. To gain the idea of (2), consider an ellipse
    \begin{equation*}
        E _{d}(A)=\left\{x \left\lvert\right. x \in \mathbb{R}^{d},x ^{\top } A ^{-1}x \le 1\right\},
    \end{equation*}
    and its squared set denoted by $E _{d}^{2}(A)$. For any point $y$ in $E ^{2}_{d}(A)$, by definition there exists a sign set $\gamma \in \left\{-1,1\right\}^{d}$ such that $\left(\gamma _{1}\sqrt{y _{1}},\dots,\gamma _{d}\sqrt{y _{d}}\right)\in E _{d}(A)$. Therefore, $E _{d}^{2}(A)$ can be expressed by 
    \begin{equation*}
        E _{d}^{2}(A)=\left\{y \left\lvert\right. y \in \mathbb{R}_{+}^{d},\min\limits _{\gamma \in \left\{-1,1\right\}^{d}}\left\{\sum\limits_{i,j}^{}A _{ij}\gamma _{i}\gamma _{j}\sqrt{y _{i}y _{j}}\right\}\le 1\right\}.
    \end{equation*}
    When $d=2$, by selecting $\gamma _{1}$ and $\gamma _{2}$ properly, $(y _{1},y _{2})\in \mathbb{E}_{d}^{2}(A)$ if and only if 
    \begin{equation*}
        A _{11}y _{1}+A _{22}y _{2}-2 \left|A _{12}\right|\sqrt{y _{1}y _{2}}\le 1.
    \end{equation*}
    Therefore, $E _{2}^{2}(A)$ is convex since it is a sublevel set of a convex function. However, for $d \ge 3$, it is impossible to let all coefficients of $\sqrt{y _{i}y _{j}}$ be negative, which leads to the non-convexity of the LHS.
\end{remark*}

Finally, we conclude this section with the definition of the strong $\epsilon $-contamination model studied in this paper.

\begin{definition}[Strong $\epsilon $-contamination model]\label{def_epsilon_contamination_model}
    Suppose that $\tilde{\mathbf{X}}=\left\{\tilde{X}_{1},\dots,\tilde{X}_{N}\right\}$ are independently drawn according to a distribution $\mathcal{P}$ under either $H _{0}$ or $H _{1}$ with unknown label. Let $\epsilon \in [0,c _{0})$ be a predetermined known constant for some $c _{0}<\frac{1}{2}$. An adversary $\mathcal{C}$ is assumed to have the full knowledge of $\mathcal{P}$. After the sampling step, $\mathcal{C}$ is allowed to inspect $\tilde{\mathbf{X}}$ and arbitrarily replace up to $\epsilon N$ of them before they are available for analysis. The resulting corrupted dataset $\mathbf{X}=\mathcal{C}(\tilde{\mathbf{X}})$ is provided to the statistician while the corruption index set $C$ is unknown.
\end{definition}

\section{Main Results}\label{section: main result}

\subsection{Lower Bounds}\label{subsection: lower bounds}

In this section, we introduce three independent lower bounds, i.e., Theorem \ref{theorem: lower bound 1}, \ref{theorem: lower bound 2}, and \ref{theorem: lower bound 3}, which jointly and fully determine the lower bound of $\rho _{\text{critical}}$.

\subsubsection{Based on a Geometric Property}\label{subsection: main lower bound 1}

The following lemma serves as an important basis for the theoretical lower bounds in Theorems \ref{theorem: lower bound 1} and \ref{theorem: lower bound 2}. The underlying intuition is that the orthosymmetry and convexity of the prior constraint $K$ imply the existence of vectors with small individual coordinates that nonetheless yield a significant aggregate norm. Such ``spread-out'' vectors, and the maximal norm they can attain, are effectively characterized by the Kolmogorov widths of $K$. Specifically, Lemma \ref{lemma: existence of the optimal vector of a QCO set} guarantees the existence of a vector with a large $\ell_2$ norm but a small $\ell_\infty$ norm by leveraging the properties of \hyperref[def_QCO_set]{QCO sets} and the \hyperref[def_Kolmogorov_width]{definition} of the Kolmogorov widths. 

\begin{lemma}\label{lemma: existence of the optimal vector of a QCO set}
    Let $K \subset \mathcal{X}$ be a QCO set and $c$ be a positive constant. Suppose that $D_{k-1}(K) > c \sigma $ for some $k \ge 1$. Then there exists a vector $\theta \in K$ such that $\left\lVert \theta \right\rVert_{2}^{} = c \sigma $ while $\left\lVert \theta \right\rVert_{\infty }^{} \le \frac{c}{\sqrt{k}}\sigma $.
\end{lemma}

The proof of Lemma \ref{lemma: existence of the optimal vector of a QCO set} is derived from the key properties of the QCO set $K$, combined with an application of the minimax theorem to a suitably truncated vector in $K ^{2}$. The formal proof is provided in Appendix \ref{subsection: proof of the existence lemma}. A similar argument also appears in Section 3.3.2 of \cite{neykov2022minimax}.

Such existence in Lemma \ref{lemma: existence of the optimal vector of a QCO set} is fundamental for the construction of the least favorable prior, where we define a uniform mixture over the coordinate-wise sign flips of this vector to induce a small total variation distance between the mixture and the null distribution $H _{0}$, thereby establishing one of the minimax lower bounds. This standard technique (see, e.g., \cite{10.1214/aos/1193342380}) can be adapted to our current problem with the prior constraint $K$ via Lemma \ref{lemma: existence of the optimal vector of a QCO set}.

\begin{theorem}[First lower bound]\label{theorem: lower bound 1}
    If $\rho \le c(\alpha ) \cdot \frac{(k _{1}^{\star })^{1/4}}{\sqrt{N}}\sigma $, where $c$ is a constant that only depends on $\alpha $, or $k _{1}^{\star }=0$, then
    \begin{equation*}
        \inf_{\psi: \mathbb{P}_{0}(\psi = 1)\le \alpha}\sup_{\theta \in K, \|\theta\| \ge \rho} \mathbb{P}_{\theta }(\psi = 0) \ge \alpha .
    \end{equation*}
    Moreover, one possible explicit expressions of $c(\alpha )$ is
    \begin{equation*}
        c(\alpha )=\left[2 \ln \left(1+4 \left(1-2 \alpha \right)^{2}\right)\right]^{\frac{1}{4}}.
    \end{equation*} 
\end{theorem}

For a fixed level $\alpha $, Theorem \ref{theorem: lower bound 1} establishes that if $\rho \lesssim \frac{(k _{1}^{\star })^{1/4}}{\sqrt{N}}\sigma $, no test can uniformly distinguish $H_1$ from $H_0$ with non-trivial power. As an immediate consequence, distinguishing $H_0$ and $H_1$ remains impossible in the presence of an adversary $\mathcal{C}$ under the same condition, since a valid strategy for $\mathcal{C}$ is simply applying the identity map (i.e., leaving the samples unperturbed).

Moreover, it is worth noting the marginal case when $k _{1}^{\star }=0$. Informally speaking, $k _{1}^{\star }=0$ means that the constraint set $K$ is insufficiently rich to support a signal that is statistically distinguishable from $H _{0}$. Since the definition of $k _{1}^{\star }$ only relies on $K$, $N$, $\sigma $, and does not depend on the original samples $\tilde{\mathbf{X}}$ or the observations $\mathbf{X}$, we shall assume in the following text that $D _{0}(K)=\sup\limits_{\theta \in K}\left\lVert \theta \right\rVert_{2}^{}>\frac{1}{\sqrt{N}}\sigma $, and therefore $k _{1}^{\star }\ge 1$. Otherwise, the test will not be feasible in the first place. We kindly refer readers to Appendix \ref{subsection: proof of the main lower bound 1} for the proof of this argument and Theorem \ref{theorem: lower bound 1}, which relies on the selected vector in Lemma \ref{lemma: existence of the optimal vector of a QCO set} and Le Cam's method of ``fuzzy hypothesis''.

Recall that for uncontaminated observations without constraint, the critical value for $\left\lVert \mu \right\rVert_{2}^{}$ is $\rho _{\text{critical}}\asymp \frac{d ^{1/4}}{\sqrt{N}}\sigma $ (see Section 3.3.6 in \cite{ingster2003nonparametric}). Theorem \ref{theorem: lower bound 1} generalizes this classic result in the sense that additional prior knowledge of the mean is allowed and the lower bound is refined as a consequence. Indeed, Theorem \ref{theorem: lower bound 1} recovers the classic result when we set $\mathcal{X}=\mathbb{R}^{d}$ (where $k _{1}^{\star }=d$ according to Example \ref{example: QCO R^d}). Moreover, the theorem implies that when we have non-trivial constraint $K$ on $\mu $, one of the essential dimensional attributes of the problem is the finite $k _{1}^{\star }$ instead of $d$. (Recall that $k _{1}^{\star }<\infty $ is ensured even if $K \subset \ell _{2}$.)

\subsubsection{Based on the Framework of Differential Privacy}\label{subsection: main lower bound 2}

Differential privacy was first introduced by \cite{10.1007/11681878_14} in the context of private data analysis. Formally, it requires the output distribution of a randomized algorithm to be insensitive to small perturbations in the input dataset. A key insight in this field is the strong connection between robust statistics and differential privacy; specifically, under mild conditions, robust estimators can be converted into differentially private algorithms using the \textit{propose-test-release} framework. Subsequently, differential privacy was applied to the goodness-of-fit testing (also called identity testing) in statistics. Pioneering works in this direction include \cite{pmlr-v48-rogers16, pmlr-v54-rogers17a, pmlr-v70-cai17a}.

\begin{definition}[$(\epsilon ,\delta )$-differential privacy, \cite{10.1007/11681878_14}]\label{def_diff_privacy_init}
    For a pair of non-negative real numbers $(\epsilon ,\delta )$, a randomized algorithm $\mathcal{A}$ defined on the sample space $\Omega $ is a $(\epsilon ,\delta )$-differentially private algorithm if for any $S \subset \text{range}(\mathcal{A})$ and any $X,X ^{\prime }\in \Omega $ with the Hamming distance $d _{H}(X, X ^{\prime })\le 1$, we have
\begin{equation}
    \mathbb{P}\left(\mathcal{A}(X)\in S\right)\le e ^{\epsilon }\mathbb{P}\left(\mathcal{A}(X ^{\prime })\in S\right)+\delta .
\end{equation}
When $\delta = 0$, the notion is called \textit{pure differential privacy}; when $\delta > 0$, it is referred to as \textit{approximate differential privacy}.
\end{definition}

Our second lower bound relies on a fundamental coupling lemma: if two datasets $\mathbf{X}$ and $\mathbf{X}^{\prime }$, drawn from $\mathcal{U}$ and $\mathcal{V}$ respectively, can be coupled to overlap on a significant fraction of indices (i.e., $X _{i}=X ^{\prime }_{i}$ for most $i$), then distinguishing between $\mathcal{U}$ and $\mathcal{V}$ is statistically hard.

\begin{lemma}[Core lemma]\label{lemma: main lemma in the proof of lower bound 2 init}
    Suppose $\mathbf{X}=(X _{1},X _{2},\dots,X _{N})^{\top }$ and $\mathbf{X} ^{\prime }=(X ^{\prime }_{1},X _{2}^{\prime },\dots,X _{N}^{\prime })^{\top }$ are random matrices with rows $X _{i},X ^{\prime }_{i}\in \mathbb{R}^{d}$. Assume the laws of $\mathbf{X}$ and $\mathbf{X} ^{\prime }$ are $\mathcal{U}$ and $\mathcal{V}$, respectively. Additionally, assume that there exists a coupling between $\mathcal{U}$ and $\mathcal{V}$ under which the expectation of the Hamming distance between $\mathbf{X}$ and $\mathbf{X} ^{\prime }$ is $\mathcal{O}(\epsilon N)$, i.e.,
\begin{equation}\label{eq: condition in the core lemma of the second lower bound}
    \mathbb{E}_{(\mathbf{X},\mathbf{X} ^{\prime })\sim (\mathcal{U},\mathcal{V})}d _{H}(\mathbf{X},\mathbf{X} ^{\prime })\lesssim \epsilon N.
\end{equation}
Then, there exists no robust test capable of distinguishing between $\mathcal{U}$ and $\mathcal{V}$ under contamination level of $\varOmega (\epsilon )$ while achieving type-I and type-II error probabilities below $0.1$ simultaneously and uniformly.
\end{lemma}

The proof of Lemma \ref{lemma: main lemma in the proof of lower bound 2 init} relies on a result established by \cite{NEURIPS2018_7de32147}, which states that if \eqref{eq: condition in the core lemma of the second lower bound} holds for a $(0,\delta)$-differentially private algorithm, then necessarily $\epsilon N=\Omega(1/\delta)$. Furthermore, for any robust test attempting to distinguish between $\mathcal{U}$ and $\mathcal{V}$, one can construct a $(0,\delta)$-differentially private algorithm via the \textit{black-box robustness-to-privacy transformation} proposed in \cite{10353143} and \cite{10.1145/3564246.3585115}. Combining these two arguments yields a necessary condition for the expected Hamming distance between $\mathbf{X}$ and $\mathbf{X}^{\prime}$.

The following theorem, which establishes the second lower bound for the testing problem \eqref{eq: testing problem}, builds upon Lemma \ref{lemma: main lemma in the proof of lower bound 2 init} and the selected vector $\theta \in K$ associated with the second optimal dimension $k _{2}^{\star }$ (see Lemma \ref{lemma: existence of the optimal vector of a QCO set}). Specifically, let $\mathcal{U}=\mathbb{P}_{0}^{\otimes N}$ and $\mathcal{V}=\mathbb{P}_{1}^{\otimes N}$, where $\mathbb{P}_{0}$ denotes $\mathcal{N}(0,\sigma ^{2})$ and $\mathbb{P}_{1}$ is the mixture of $\mathcal{N}(v,\sigma ^{2})$ with $v$ drawn uniformly from the set
\begin{equation*}
    V :=\left\{v \left\lvert\right. v _{i}=\gamma _{i}\theta _{i}, \gamma _{i}\in \left\{-1,1\right\}, i \ge 1, i \in I\right\}.
\end{equation*}
We construct a chain of distributions interpolating between $\mathcal{U}$ and $\mathcal{V}$ based on the optimal coupling (Lemma \ref{lemma: total variation and coupling}). Consequently, the computation of the expected Hamming distance is reduced to the sum of the total variation distances between adjacent distributions in the chain. We kindly refer readers to Appendix \ref{subsection: proof of the main lower bound 2} for the formal proofs of Lemma \ref{lemma: main lemma in the proof of lower bound 2 init} and Theorem \ref{theorem: lower bound 2}.

\begin{theorem}[Second lower bound]\label{theorem: lower bound 2}
    Consider the testing problem \eqref{eq: testing problem} under the \hyperref[def_epsilon_contamination_model]{strong $\epsilon$-contamination model} with an adversary $\mathcal{C}$. If the signal strength satisfies $\rho \lesssim \frac{(k _{2}^{\star }) ^{1/4}\sqrt{\epsilon }\sigma }{N ^{1/4}}$, then no test can simultaneously achieve type-I and type-II errors bounded by $0.1$. Formally, under this condition, we have
    \begin{equation*}
        \inf\limits_{\phi \in A _{s}(0.1)}\sup\limits_{\mu \in K}\sup\limits_{\mathcal{C}}\mathbb{P}_{\mu }(\phi (\mathcal{C}(\mathbf{X}^{\prime }))=0)\ge 0.1,
    \end{equation*}
    where $A _{s}(0.1)$ denotes the set of tests with type-I error at most $0.1$, as defined in \eqref{eq: def_A_s}.
\end{theorem}

In Lemma \ref{lemma: main lemma in the proof of lower bound 2 init} and Theorem \ref{theorem: lower bound 2}, we fix the constant at $0.1$ for simplicity, though it can be replaced by any sufficiently small universal constant. In the high-dimensional setting where $d$ (and possibly $k_2^\star$) may substantially exceed $N$, one must ensure that
\begin{equation*}
    \frac{(k _{2}^{\star })^{1/4}\sqrt{\epsilon }\sigma }{N ^{1/4}}\le D _{0}(K)=\sup\limits_{\theta \in K}\left\lVert \theta \right\rVert_{2}^{}.
\end{equation*}
If this condition is violated, the testing problem becomes statistically impossible.

\subsubsection{Based on the Estimation Lower Bounds}\label{subsection: main lower bound 3}

The third lower bound, as established in Theorem \ref{theorem: lower bound 3} below, contrasts with Theorems \ref{theorem: lower bound 1} and \ref{theorem: lower bound 2} in that it is determined solely by the diameter of $K$ (denoted as $D _{0}(K)$), the fraction of contamination $\epsilon$, and the noise scale $\sigma$, while remaining independent of $N$.

\begin{theorem}[Third lower bound]\label{theorem: lower bound 3}
    For the signal detection problem \eqref{eq: testing problem} under the \hyperref[def_epsilon_contamination_model]{strong $\epsilon$-contamination model} with an adversary $\mathcal{C}$, if $\epsilon \gtrsim \frac{1}{\sqrt{N}}$ and $\rho \lesssim \epsilon \sigma $, we always have 
\begin{equation}
    \inf\limits_{\phi \in A _{s}}\sup\limits_{\mu \in K,\left\lVert \mu \right\rVert_{2}^{}\ge \rho }\sup\limits_{\mathcal{C}}\mathbb{P}_{\mu }(\phi (\mathcal{C}(\tilde{\mathbf{X}}))=0)= \varTheta (1),
    \label{eq: lower_bound_3_1}
\end{equation}
where $A _{s}$ is the set of acceptable tests defined as \eqref{eq: def_A_s}.
\end{theorem}

\begin{remark*}
    (1). This lower bound is notable as it coincides with the minimax rate for the estimation problem under the strong $\epsilon$-contamination model (see, e.g., \cite{prasadan2025informationtheoreticlimitsrobust}). This finding stands in contrast to the conventional wisdom that testing is typically statistically easier than estimation (see \cite{ingster1993asymptotically1, baraud2002non}). However, our subsequent analysis confirms that in this setting, the minimax rates for estimation and testing indeed match (or are of the same order) when certain conditions are met.

    (2). The condition $\epsilon \gtrsim \frac{1}{\sqrt{N}}$ can be safely removed. See our discussions in Section \ref{subsection: main result}.
\end{remark*}

Recall that $D _{0}(K)=\sup_{\theta \in K}\lVert \theta \rVert_{2}$ represents the radius of $K$ (i.e., the maximal distance from the origin). Theorem \ref{theorem: lower bound 3} implies that a necessary condition for the existence of a test with type-I and type-II errors uniformly bounded by $\alpha$ is that either $\rho \gtrsim \epsilon \sigma$ or $\epsilon =\mathcal{O}(1/\sqrt{N})$. Noting that the signal strength is inherently bounded by $\lVert \mu \rVert_{2} \le D _{0}(K)$, it follows that if $\epsilon \gtrsim 1/\sqrt{N}$ and $\epsilon \sigma \gtrsim D _{0}(K)$, consistent testing is statistically impossible.

The proof of Theorem \ref{theorem: lower bound 3} is inspired by a similar lower bound established for the estimation problem under star-shaped constraints in \cite{prasadan2025informationtheoreticlimitsrobust}. It relies on the construction of a specific distribution that admits two distinct mixture representations, with components corresponding to $H_{0}$ and $H_{1}$ respectively. Crucially, the adversary $\mathcal{C}$ can exploit this structural ambiguity to design a contamination strategy that renders the resulting distributions statistically indistinguishable when the lower bound condition is violated. The full proof is detailed in Appendix \ref{subsection: proof of the main lower bound 3}.

\subsection{Upper Bounds}\label{subsection: upper bounds}

The upper bounds are established constructively. Specifically, we design explicit testing procedures and identify the required signal strength $\rho =\lVert \mu \rVert_{2}$ that guarantees test validity (i.e., uniformly bounded type-I and type-II errors).

\subsubsection{Theoretical Algorithm}\label{subsection: theoretical algorithm}

Since the adversary $\mathcal{C}$ is allowed to modify only a limited number of samples, a natural strategy to mitigate the contamination is to identify subsets of $\mathbf{X}$ that exhibit a \textit{consistency property}.

\begin{definition}[Consistent subset]\label{def_consistent_subset_init}
    A subset $S \subset \mathbf{X}=\left\{X _{1},\dots,X _{N}\right\}$ is called a consistent subset regarding the contamination fraction $\epsilon $ and a test $\phi :2 ^{S}\mapsto \left\{0,1\right\}$ if\\ 
    \hspace*{0.5em}(i), $\left|S\right|\ge (1-\epsilon )N$,\\ 
    \hspace*{0.5em}(ii), $\phi (S ^{\prime })=\phi (S)$ for any $S ^{\prime }\subset S$ with $\left|S ^{\prime }\right|\ge (1-2\epsilon )N$.
\end{definition}

As shown in Appendix \ref{subsection: proof of the main upper bound 1}, such subsets exist with high probability and behave as if they consisted entirely of authentic samples. To provide some intuition, consider a specific $\chi ^{2}$-type test $\phi _{e}$. Assume that the original (unobserved) data $\tilde{\mathbf{X}}$ is inherently consistent. This guarantees the existence of a consistent subset within $\mathbf{X}$ (at least for the set of uncorrupted samples $\mathbf{X}_{[N]\backslash C}$). Let $\mathbf{X}_{S}$ denote any such consistent subset. Consequently, we can recover the test result of the unobserved $\tilde{\mathbf{X}}$ by establishing the following equivalence by the definition of the \hyperref[def_consistent_subset_init]{consistent subset}:
\begin{equation*}
    \phi _{e}(\mathbf{X}_{S})=\phi _{e}(\mathbf{X}_{S}\cap \mathbf{X}_{[N]\backslash C})=\phi _{e}(\mathbf{X}_{[N]\backslash C})=\phi _{e}(\tilde{\mathbf{X}}).
\end{equation*}

The following theorem formally establishes the existence of a valid testing procedure under the condition \eqref{eq: required condition of the main upper bound 1}, and consequently provides an upper bound for the testing problem \eqref{eq: testing problem}. Specifically, we demonstrate that, when \eqref{eq: required condition of the main upper bound 1} holds, there exists a desired test $\phi _{e}$ and at least one consistent subset $\mathbf{X}_{S}$ with respect to $\phi _{e}$, upon which the entire testing task becomes feasible. To achieve this, we explicitly define an event $E _{e}$ of $\chi ^{2}$-type, on which $\phi _{e}$ is determined, and prove that the original samples $\tilde{\mathbf{X}}$ is consistent with high probability. The detailed proof is provided in Appendix \ref{subsection: proof of the main upper bound 1}.

\begin{theorem}[Theoretical upper bound]\label{theorem: upper bound 1}
    Given a dataset $\mathbf{X}$ from the \hyperref[def_epsilon_contamination_model]{strong $\epsilon $-contamination model} and a fixed error tolerance $\alpha $, for any $k \in [1,d]$\footnote{When $K \subset \ell _{2}$, we allow $1 \le k< \infty $.}, there exists a test $\phi _{e}$ for \eqref{eq: testing problem} such that the type-I and type-II errors are both uniformly less than $\alpha $ if $\rho ^{2}\gtrsim E _{\text{raw}}^{2}(\epsilon ,K,N,\sigma ^{2},k)$ where $E _{\text{raw}}^{2}(\epsilon ,K,N,\sigma ^{2},k)$ is defined as
\begin{equation}\label{eq: required condition of the main upper bound 1}
    E _{\text{raw}}^{2}(\epsilon ,K,N,\sigma ^{2},k):=D _{k}^{2}(K)+\sigma ^{2}\max\limits \left\{\frac{\sqrt{k}}{N},\epsilon ^{2}\ln \left(\frac{1}{\epsilon }\right),\sqrt{\frac{\epsilon ^{2}\ln \left(\frac{1}{\epsilon }\right)k}{N}}\right\}.
\end{equation}
In other words, we have
\begin{equation*}
    \rho _{\text{critical}}\lesssim \min\limits _{k}E _{\text{raw}}(\epsilon ,K,N,\sigma ^{2},k).
\end{equation*}
\end{theorem}

We can further leverage the results of Theorem \ref{theorem: upper bound 1} by optimizing over $k$ and applying the definition of $k _{1}^{\star }$ and $k _{2}^{\star }$.

\begin{corollary}[]\label{corollary: upper bound 1}
    Given the same setting as Theorem \ref{theorem: upper bound 1}, we can alternatively define $E(\epsilon ,K,N,\sigma ^{2})$, which is also a valid upper bound of $\rho _{\text{critical}}$, as
    \begin{equation}\label{eq: alt required condition of the main upper bound 1}
        E ^{2}(\epsilon ,K,N,\sigma ^{2}):=\sigma ^{2}\max\limits \left\{\frac{\sqrt{\min\limits \left\{k _{1}^{\star },k _{2}^{\star }\right\}}}{N},\epsilon ^{2}\ln \left(\frac{1}{\epsilon }\right),\sqrt{\frac{\epsilon ^{2}\ln \left(\frac{1}{\epsilon }\right)\min\limits \left\{k _{1}^{\star },k _{2}^{\star }\right\}}{N}}\right\}
    \end{equation}
\end{corollary}
\begin{proof}
    By the definition of $k _{1}^{\star }$, we know $D _{k _{1}^{\star }}(K)\le \frac{\sqrt{k _{1}^{\star }+1}}{N}\sigma ^{2}$. Therefore, if we set $k=k _{1}^{\star }$ in \eqref{eq: required condition of the main upper bound 1}, $D _{k _{1}^{\star }}(K)$ can be absorbed into the other term involving the maximum. The logic also holds for $k _{2}^{\star }$. Finally, since \eqref{eq: required condition of the main upper bound 1} is always a valid upper bound, we select the smaller one, which leads to \eqref{eq: alt required condition of the main upper bound 1}.
\end{proof}

Corollary \ref{corollary: upper bound 1} effectively selects the smaller of the optimal dimensions $k _{1}^{\star }$ and $k _{2}^{\star }$, and consequently determines the optimal projection on which $\phi _{e}$ is based. We show in Section \ref{subsection: main result} that the upper bound \eqref{eq: alt required condition of the main upper bound 1} matches the theoretical lower bounds up to a logarithmic factor of $\ln(1/\epsilon)$.

A significant improvement implied by Theorem \ref{theorem: upper bound 1} is that it does not require the asymptotic condition $\epsilon \rightarrow 0$, which is an essential assumption in many previous works. However, the test $\phi _{e}$ is limited by its reliance on an exhaustive search over all subsets of $\mathbf{X}$, resulting in exponential computational complexity with respect to $N$. This prohibitive cost renders the test computationally infeasible for large sample sizes. In the next section, motivated by the filtering method in \cite{10353143}, we propose a polynomial-time algorithm designed to handle the prior constraint $K$ in the presence of the adversary $\mathcal{C}$.

\subsubsection{Polynomial-Time Algorithm}\label{subsection: polynomial-time algorithm}

In this section, we develop a polynomial-time algorithm\footnote{By polynomial-time'', we assume that computing the optimal projections (determining the optimal dimensions $k _{1}^{\star }$ and $k _{2}^{\star }$) can be performed in polynomial time, as is the case for any ellipsoid or hyperrectangle.} to detect whether $\mu =\mathbf{0}$ under a slightly stronger condition than the bound \eqref{eq: required condition of the main upper bound 1}. In contrast to the theoretical algorithm which requires an exhaustive search, the proposed algorithm leverages weights to adaptively filter potentially corrupted samples. Specifically, this filtering strategy and the algorithm are built and developed upon the framework of \hyperref[def_omega_regularity_init]{$\omega $-regularity}, established in prior works (e.g., \cite{10353143, dong2019quantumentropyscoringfast}).

\begin{definition}[$\omega $-regularity]\label{def_omega_regularity_init}
    Given a weight vector $\omega =(\omega _{1},\dots,\omega _{N})^{\top } $, $\mathbf{Y}$ is said to be $(\epsilon ,\beta _{1},\beta _{2})$-regular if for all subsets $S \subset [N]$ with $\left|S\right|\le \epsilon N$, we have the following properties:\\ 
    \begin{tabular}{rl}
        (i), & $\left|\sum\limits_{i \in S}^{}\left\lVert Y _{i}\right\rVert_{2}^{2}-\left|S\right|k ^{\star }\right|\le c \beta _{1}$,\\ 
        (ii), & $\left|\left\lVert \sum\limits_{i \in S}^{}\sqrt{\omega }_{i} Y _{i}\right\rVert_{2}^{2}-\left\lVert \omega _{S}\right\rVert_{1}^{}k ^{\star }\right|\le c \beta _{2}$, and\\ 
        (iii), & $\left|\left\langle \sum\limits_{i \in S}^{}\sqrt{\omega }_{i}Y _{i}, \sum\limits_{j \in [N]}^{}\sqrt{\omega }_{j}Y _{j} \right\rangle-\left\lVert \omega _{S}\right\rVert_{1}^{} k ^{\star }\right|\le c \sqrt{N}\beta _{1}$,
    \end{tabular}

    where we recall that $\omega _{S} \in \mathbb{R}^{N}$ means the restriction of $\omega $ on the set $S$.
\end{definition}

The intuition behind the polynomial-time algorithm parallels that of the theoretical construction, although establishing this connection relies on more advanced random matrix theory concerning (sub-)Gaussian random variables (see Appendix \ref{section: concepts and lemmas} for technical preliminaries). We first establish that the authentic samples satisfy the \hyperref[def_omega_regularity_init]{$\omega $-regularity} condition for specified parameters $\beta _{1},\beta _{2}$ with high probability (setting the initial weights $\omega =\mathbf{1}_{[N]}$). Our analysis demonstrates that, under sufficient conditions on the signal strength $\rho =\lVert \mu \rVert_{2}$, the re-weighted samples effectively mimic the behavior of the unobserved authentic data $\tilde{\mathbf{X}}$. Consequently, we are able to recover the test result corresponding to $\tilde{\mathbf{X}}$. Such conditions are formalized in the following theorem.

\begin{theorem}[Polynomial upper bound]\label{theorem: upper bound 2}
    Assuming the same conditions as in Theorem \ref{theorem: upper bound 1}, for any $1 \le k \le d$\footnote{When $K \subset \ell _{2}$, we allow $1 \le k <\infty $.}, as long as the following condition is satisfied, there exists a test $\phi _{p}$ to distinguish between $H _{0}$ and $H _{1}$ with type-I and type-II errors both uniformly less than $\alpha $.
    \begin{equation}\label{eq: required condition of the main upper bound 2}
        \rho ^{2}\gtrsim \underbrace{D _{k}^{2}(K)+\sigma ^{2}\max\limits \left\{\frac{\epsilon \ln \left(\frac{N}{\alpha }\right)}{\sqrt{N}},\epsilon ^{2}\ln \left(\frac{N}{\alpha }\right),\sqrt{\frac{\epsilon ^{2}k \ln \left(\frac{N}{\alpha }\right)}{N}}, \frac{\sqrt{k}\ln \left(\frac{1}{\alpha }\right)}{N}\right\}}_{:=P _{\text{raw}}^{2}(\epsilon ,K,N,\sigma ^{2},k)}.
    \end{equation}
    In other words, we have 
    \begin{equation*}
        \rho _{\text{critical}}\lesssim \min\limits _{k}P _{\text{raw}}(\epsilon ,K,N,\sigma ^{2},k)
    \end{equation*}
    In addition, the time and space complexity of the algorithm is at most polynomial in terms of the parameters $N,k,\frac{1}{\epsilon } ,\frac{1}{\alpha }$.
\end{theorem}

Similar to Theorem \ref{theorem: upper bound 1}, the following corollary characterizes the upper bound \eqref{eq: required condition of the main upper bound 2} with $k _{1}^{\star }$ and $k _{2}^{\star }$. The proof is omitted since it follows from Corollary \ref{corollary: upper bound 1} with minor modifications.

\begin{corollary}[]\label{corollary: upper bound 2}
    Given the same setting as Theorem \ref{theorem: upper bound 2}, we can alternatively define $P(\epsilon ,K,N,\sigma ^{2})$, which is also a valid upper bound of $\rho _{\text{critical}}$, as 
    \begin{equation}\label{eq: alt required condition of the main upper bound 2}
        P ^{2}(\epsilon ,K,N,\sigma ^{2})=\sigma ^{2}\max\limits \left\{\frac{\epsilon \ln \left(\frac{N}{\alpha }\right)}{\sqrt{N}},\epsilon ^{2}\ln \left(\frac{N}{\alpha }\right),\sqrt{\frac{\epsilon ^{2}\min\limits \left\{k _{1}^{\star },k _{2}^{\star }\right\}\ln \left(\frac{N}{\alpha }\right)}{N}},\frac{\sqrt{\min\limits \left\{k _{1}^{\star },k _{2}^{\star }\right\}}\ln \left(\frac{1}{\alpha }\right)}{N}\right\}.
    \end{equation}
\end{corollary}

Analogous to the upper bound \eqref{eq: alt required condition of the main upper bound 1}, \eqref{eq: alt required condition of the main upper bound 2} also nearly matches the theoretical lower bounds with only the exception of logarithmic terms. However, one can easily check that the upper bound \eqref{eq: alt required condition of the main upper bound 1} is generally tighter than \eqref{eq: alt required condition of the main upper bound 2}, i.e., $E(\epsilon ,K,N,\sigma ^{2})\le P(\epsilon ,K,N,\sigma ^{2})$.

We acknowledge that the proof of Theorem \ref{theorem: upper bound 2} and the associated algorithm build upon the framework established in \cite{10353143}. However, the results in \cite{10353143} are restricted to the regime where $\lVert \mu \rVert_{2} \lesssim 1$ in the unconstrained setting (i.e., $K = \mathbb{R}^{d}$). Consequently, we have significantly extended the algorithm to accommodate arbitrary signal magnitudes and to incorporate the prior constraint $K$ via the optimal projection $P^{\star}$. Furthermore, our analysis refines the arguments of \cite{10353143} to address certain technical issues. We kindly refer readers to Appendix \ref{subsection: proof of the main upper bound 2} for the complete proof of the properties regarding the \hyperref[def_omega_regularity_init]{$\omega $-regularity} and Theorem \ref{theorem: upper bound 2} and Appendix \ref{section: algorithms} for the pseudo code of the polynomial-time algorithm.

\subsection{Main Results}\label{subsection: main result}

In this section, we integrate the preceding theorems and bounds and summarize the main results. The lower bound is synthesized from Theorems \ref{theorem: lower bound 1}, \ref{theorem: lower bound 2}, \ref{theorem: lower bound 3}. Although the bounds are derived through different approaches, an interesting phase transition phenomenon emerges, determined by the ranges of $\epsilon $ over which each term dominates others.

\begin{theorem}[Main lower bound]\label{theorem: main lower bound}
    Given the mean testing problem \eqref{eq: testing problem} in the Gaussian sequence model, with potential contamination from a strong adversary $\mathcal{C}$ and a prior constraint on the mean $\mu$ induced by a QCO set $K$, the following condition is necessary to ensure the existence of a valid test whose type-I and type-II errors are both uniformly below a prescribed constant $\alpha $
    \begin{equation}\label{eq: main theoretical lower bound}
        \rho ^{2}\gtrsim \sigma ^{2}\max\limits \left\{\frac{\sqrt{k _{1}^{\star }}}{N},\epsilon \sqrt{\frac{k _{2}^{\star }}{N}},\epsilon ^{2}\right\},
    \end{equation}
    where $N$ is the sample size, $\sigma ^{2}$ is the variance of the noise, $\epsilon $ is the fraction of contamination, and $k _{1}^{\star },k _{2}^{\star }$ are the optimal dimensions exclusively determined by $N$, $\sigma $, $\epsilon $ and the set $K$ as in Definition \ref{def_first_optimal_dimension}, \ref{def_second_optimal_dimension}.
\end{theorem}

Roughly speaking, Theorem \ref{theorem: main lower bound} demonstrates that the minimax lower bound of the testing problem \eqref{eq: testing problem} is jointly determined by three parts --- the geometric properties of the set $K$, which is represented by the factors $k _{1}^{\star },k _{2}^{\star }$; the properties of the original data, which is implied by $N$ and $\sigma $; and the corruption process, which is implied by $\epsilon $.

The corruption fraction $\epsilon $ connects the three components in the lower bound. In the \hyperref[subsection: proof of the main corollary]{proof} of Corollary \ref{corollary: matching between bounds} below, we show that the first term, which is independent with $\epsilon $\footnote{We should recall that $k _{1}^{\star }$ is independent with $\epsilon $}, dominates the other two if and only if $\epsilon \le \frac{1}{\sqrt{N}}$. Similarly, the last term, which is independent with $N$, is the dominant if and only if $\epsilon \gtrsim \sqrt{\frac{k _{2}^{\star }}{N}}$. Therefore, the set $\left\{\epsilon \left\lvert\right. 0 \le \epsilon <c _{0},\epsilon \in \left[\frac{1}{\sqrt{N}},\sqrt{\frac{k _{2}^{\star }}{N}}\right]\right\}$\footnote{We express in this form because $k _{2}^{\star }$ also depends on $\epsilon $.} serves as a threshold belt. When $\epsilon $ lies below this belt, it is too small for the adversary $\mathcal{C}$ to exert its influence, and the optimal minimax rate becomes essentially independent of $\epsilon $. In other words, it is free to contaminate the dataset in this region. In contrast, if $\epsilon $ lies above the belt (for example, $\epsilon \asymp N ^{-\frac{1}{4}}$), then the power of $\mathcal{C}$ can be fully utilized, to the extent that increasing the sample size $N$ does not improve the problem's difficulty. When $\epsilon $ falls within the belt, the minimax rate is the result of the delicate balance between the adversary $\mathcal{C}$ and the testing procedures $\phi _{e},\phi _{p}$. Based on this argument, \eqref{eq: main theoretical lower bound} can be equivalently written as 
\begin{equation}\label{eq: alter main theoretical lower bound}
    \rho ^{2}\gtrsim \left\{
    \begin{tabular}{lll}
        $\frac{\sqrt{k _{1}^{\star }}}{N}\sigma ^{2}$, & $\left(0 \le \epsilon \lesssim \frac{1}{\sqrt{N}}\right)$ & (Free contamination phase),\\ 
        $\epsilon \sqrt{\frac{k _{2}^{\star }}{N}}\sigma ^{2}$, & $\left(\frac{1}{\sqrt{N}}\lesssim \epsilon \lesssim \sqrt{\frac{k _{2}^{\star }}{N}}\right)$ & (Balance phase),\\ 
        $\epsilon ^{2}\sigma ^{2}$, & $\left(\epsilon \gtrsim \sqrt{\frac{k _{2}^{\star }}{N}}\right)$, & (Overly contaminated phase).
    \end{tabular}
    \right.
\end{equation}

Informally, a message implied by \eqref{eq: alter main theoretical lower bound} is that after the ``safe threshold'' $\frac{1}{\sqrt{N}}$, every single contaminated sample in the strong $\epsilon $-contamination model invalidates the statistical testing use of $\sqrt{\frac{N}{k _{2}^{\star }}}$ uncontaminated samples on average, until the overly contaminated phase.

For the upper bound, Theorem \ref{theorem: main upper bound} below is a natural downstream conclusion of Theorem \ref{theorem: upper bound 1} and \ref{theorem: upper bound 2}.

\begin{theorem}[Main upper bound]\label{theorem: main upper bound}
    Given the same setting as in Theorem \ref{theorem: main lower bound}, if $\rho \gtrsim E(\epsilon ,K,N,\sigma ^{2})$, there exists a test $\phi _{e}$ achieving type-I and type-II errors uniformly less than $\alpha $. Furthermore, if $\rho \gtrsim P(\epsilon ,K,N,\sigma ^{2})$ and the determination of $k _{1}^{\star }$ and $k _{2}^{\star }$ can be completed within polynomial-time of $\left(N,d,\frac{1}{\epsilon }\right)$, there exists a test $\phi _{p}$ that is able to finish the testing process within polynomial-time of $\left(N,k _{1}^{\star },k _{2}^{\star },\frac{1}{\epsilon },\frac{1}{\alpha }\right)$ while still achieving type-I and type-II errors uniformly less than $\alpha $. $E(\epsilon ,K,N,\sigma ^{2})$ and $P(\epsilon ,K,N,\sigma ^{2})$ are defined in \eqref{eq: alt required condition of the main upper bound 1} and \eqref{eq: alt required condition of the main upper bound 2} respectively.
\end{theorem}

We further establish that the main lower bound in Theorem \ref{theorem: main lower bound} is minimax optimal and there is only minor difference between the conditions required by the two algorithms in Theorem \ref{theorem: main upper bound}. This is summarized in the following corollary.

\begin{corollary}[Nearly matching between the lower bounds and the upper bounds]\label{corollary: matching between bounds}
    The upper bound $E(\epsilon ,K,N,\sigma ^{2})$ nearly matches the theoretical lower bound \eqref{eq: main theoretical lower bound} with only an exception of the logarithmic factor $\ln \left(\frac{1}{\epsilon }\right)$. The upper bound $P(\epsilon ,K,N,\sigma ^{2})$ of $\phi _{p}$ nearly matches $E(\epsilon ,K,N,\sigma ^{2})$ with only exceptions of the logarithmic factors $\ln N$, $\ln \alpha $, and $\ln \left(\frac{1}{\epsilon }\right)$.
\end{corollary}

For Corollary \ref{corollary: matching between bounds}, it involves some refined analysis with the main lower bound \eqref{eq: main theoretical lower bound} and $E(\epsilon ,K,N,\sigma ^{2})$. We kindly refer readers to Appendix \ref{subsection: proof of the main corollary} for the detailed arguments.

\section{Extension to \texorpdfstring{$\ell_{p}$}{lp} Norm Separation for \texorpdfstring{$1 \le p <2$}{1<=p<2}}\label{section: extension to lp testing problem}

As an alternative, with the same assumptions on the samples $\tilde{\mathbf{X}}$ and the adversary $\mathcal{C}$, we may consider the following $\ell _{p}$ robust testing problem.
\begin{equation}\label{eq: lp testing problem}
    \begin{aligned}
        & H _{0}: \mu =0,\\ 
        & H _{1}: \left\lVert \mu \right\rVert_{p}^{}\ge \rho ,\mu \in K.
    \end{aligned}
\end{equation}
The classic results for the testing problem \eqref{eq: testing problem} can be naturally generalized to \eqref{eq: lp testing problem} when $1 \le p<2$ via the H\"older's inequality by \cite{ingster2003nonparametric}. While Ingster's approach cannot be directly transferred to our settings, we establish similar generalized results with the $p$-convex condition and accordingly adjusted optimal dimensions and projections. Note that in this section, we shall substitute the $\ell _{2}$ norm with the $\ell _{p}$ norm for the Kolmogorov $k$-width. Furthermore, to simplify the analysis, we only consider the discretized version (see the remark after Definition \ref{def_Kolmogorov_width}).

\begin{definition}[Discretized Kolmogorov $k$-Width]\label{def_discretized_Kolmogorov_width}
    Let $\mathcal{X}$ be a Banach space equipped with the norm $\left\lVert \cdot \right\rVert_{}^{}$, and $K \subset \mathcal{X}$ is a subset. The discretized Kolmogorov $k$-width is defined as 
    \begin{equation}\label{eq: definition of discretized Kolmogorov width}
        \tilde{D} _{k}(K) = \inf_{P \in \tilde{\mathcal{P}}_{k}} \sup_{\theta \in K} \left\lVert \theta -P \theta \right\rVert_{}^{},
    \end{equation}

    where $\tilde{\mathcal{P}}_{k}$ is the set of all projection operators that align with the axes and project a vector onto some subspace of $\mathcal{X}$ with intrinsic dimension $k$.
\end{definition}

Though discretized, it inherits most of important properties from the original definition. When the norm is $\left\lVert \cdot \right\rVert_{p}^{}$, let us denote the corresponding Kolmogorov $k$-width for $K$ as $\tilde{D}_{p,k}(K)$.

\begin{definition}[$p$-convex orthosymmetric (PCO) set]\label{definition: p-convex set}
    For $p \ge 1$, a set $K \subset \mathcal{X}$ is a $p$-convex orthosymmetric set if\\
    \hspace*{0.5em}(1). $K$ is convex;\\ 
    \hspace*{0.5em}(2). $K ^{p}$ is convex, where $K ^{p}$ is defined as
    \begin{equation*}
        K ^{p}:=\left\{\left(\left|\theta _{1}\right|^{p},\dots,\left|\theta _{d}\right|^{p}\right)^{\top } \left\lvert\right. \left(\theta _{1},\dots,\theta _{d}\right)^{\top } \in K\right\};
    \end{equation*}
    \hspace*{0.5em}(3), $K$ is orthosymmetric, which means that if $\theta =\left(\theta _{1},\dots,\theta _{d}\right)^{\top } \in K$, then $\theta _{\eta }:=\left(\eta _{1}\theta _{1},\dots,\eta _{d}\theta _{d}\right)^{\top }\in K$, where $\eta _{i}\in \left\{-1,1\right\}, 1 \le i \le d$.
\end{definition}
It is not difficult to verify that the whole space, hyperrectangles and ellipsoids are still all PCO sets for $1 \le p <2$.

Analogous to previous analysis with the $\ell _{2}$ norm, the new lower bound relies on the following revised version of Lemma \ref{lemma: existence of the optimal vector of a QCO set}. It is not difficult to check that the proof is not essentially different.

\begin{lemma}[]\label{lemma: general existence of the optimal vector in a QCO set}
    Let $K \subset \mathcal{X}$ be a $p$-convex orthosymmetric set, and suppose that $\tilde{D}_{p,k-1}(K) > c \sigma $. Then there exists a vector $\theta \in K$ such that $\left\lVert \theta \right\rVert_{p}^{} = c \sigma $ and $\left\lVert \theta \right\rVert_{\infty }^{}\le \frac{c}{k ^{1/p}}\sigma $.
\end{lemma}

We also require the following new optimal dimensions and projections for $\ell _{p}$ robust testing problem.

\begin{definition}[$\ell _{p}$ optimal dimensions and projections]\label{general_optimal_dimension}
    Given the noise scale $\sigma $, the number of samples $N$, the prior knowledge $K$ and the corruption rate $\epsilon $, define the following optimal dimensions
    \begin{equation*}
    \begin{aligned}
        & k _{p,1}^{\star }:=\max\limits \left\{j \left\lvert\right. 0 \le j \le d, \tilde{D} _{p,j}(K)> \frac{j ^{\frac{1}{p}-\frac{1}{4}}}{\sqrt{N}}\sigma \right\},\\ 
        & k _{p,2}^{\star }:=\max\limits \left\{j \left\lvert\right. 0 \le j \le d, \tilde{D} _{p,j}(K)> \frac{\sqrt{\epsilon }j ^{\frac{1}{p}-\frac{1}{4}}}{N ^{\frac{1}{4}}}\sigma \right\},\\ 
        & k _{p,3}^{\star }:=\max\limits \left\{j \left\lvert\right. 0 \le j \le d, \tilde{D} _{p,j}(K)> \epsilon \sigma j ^{\frac{1}{p}-\frac{1}{2}}\right\}.
    \end{aligned}
    \end{equation*}
    Analogous to Definition \ref{def_first_optimal_dimension} and \ref{def_second_optimal_dimension}, the corresponding optimal projections are defined accordingly and denoted as $P _{p,1}^{\star },P _{p,2}^{\star },P _{p,3}^{\star }$.
\end{definition}

Given the modified lemma and definitions above, we present the lower bounds for the $\ell _{p}$ robust testing problem \eqref{eq: lp testing problem}.

\begin{theorem}[$\ell _{p}$ main lower bound]\label{theorem: lp main lower bound}
    Given the mean testing problem \eqref{eq: lp testing problem} in the Gaussian sequence model, with potential contamination from a strong adversary $\mathcal{C}$ and a prior constraint on the mean $\mu$ induced by a p-convex orthosymmetric set $K$, the following condition is necessary to ensure the existence of a valid test whose type-I and type-II errors are both uniformly below a prescribed constant $\alpha $
    \begin{equation}\label{eq: lp main lower bound}
    \left\lVert \theta \right\rVert_{p}^{}\gtrsim \sigma \max\limits \left\{\frac{(k _{p,1}^{\star })^{\frac{1}{p}-\frac{1}{4}}}{\sqrt{N}},\frac{\sqrt{\epsilon }(k _{p,2}^{\star })^{\frac{1}{p}-\frac{1}{4}}}{N ^{\frac{1}{4}}},\epsilon (k _{p,3}^{\star })^{\frac{1}{p}-\frac{1}{2}}\right\}.
    \end{equation}
    where $N$ is the sample size, $\sigma ^{2}$ is the variance of the noise, $\epsilon $ is the fraction of contamination, and $k _{p,1}^{\star },k _{p,2}^{\star },k _{p,3}^{\star }$ are the optimal dimensions exclusively determined by $N$, $\sigma $, $\epsilon $ and the set $K$ as in Definition \ref{general_optimal_dimension}.
\end{theorem}

Undoubtedly, the lower bound \eqref{eq: lp main lower bound} recovers \eqref{eq: main theoretical lower bound} as expected when $p=2$. A significant difference from the $\ell _{2}$ norm separation is that we have another optimal dimension regarding the last term in \eqref{eq: lp main lower bound}. The intrinsic reason is that when covariance matrices are known and identical, the KL divergence between Gaussian distributions --- an essential quantity we leverage in the \hyperref[subsection: proof of the main lower bound 3]{proof} of Theorem \ref{theorem: lower bound 3}, depends on the $\ell _{2}$ norm between the means. When considering the $\ell _{p}$ norm separation, we require $k _{p,3}^{\star }$ to adapt to this discrepancy. We defer the proof of Theorem \ref{theorem: lp main lower bound} to Appendix \ref{subsection: proof of the lp main lower bound}.

In $\ell _{2}$ norm separation, the upper bound \eqref{eq: alt required condition of the main upper bound 1} matches the lower bound \eqref{eq: main theoretical lower bound} except for logarithmic factors. We identify similar results in $\ell _{p}$ norm separation, which is summarized as Theorem \ref{theorem: lp main upper bound} below.

\begin{theorem}[$\ell _{p}$ main upper bound]\label{theorem: lp main upper bound}
    If $\rho \gtrsim \sigma \max\limits \left\{\frac{(k _{p,1}^{\star })^{\frac{1}{p}-\frac{1}{4}}}{\sqrt{N}},\frac{\sqrt{\epsilon }(k _{p,2}^{\star })^{\frac{1}{p}-\frac{1}{4}}}{N ^{\frac{1}{4}}},\epsilon (k _{p,3}^{\star })^{\frac{1}{p}-\frac{1}{2}}\right\}$, then there exists a test to distinguish between $H _{0}$ and $H _{1}$ with uniformly small type-I and type-II errors. Moreover, the time complexity of the algorithm is polynomial-time in $\left(N,d,\frac{1}{\epsilon },\frac{1}{\alpha }\right)$. Note that in the expression above, we omit the polylogarithmic factors in $\left(N,k _{p,1}^{\star },k _{p,2}^{\star },k _{p,3}^{\star },\frac{1}{\alpha },\frac{1}{\epsilon }\right)$.
\end{theorem}

At first examination, Theorem \ref{theorem: lp main upper bound} appears to be a direct and natural corollary of Theorem \ref{theorem: main upper bound}. However, a closer inspection reveals subtleties in its derivation. A crucial observation is that, though one can select the minimum $k _{p}^{\star }:=\min\limits \left\{k _{p,1}^{\star }k _{p,2}^{\star },k _{p,3}^{\star }\right\}$ and its optimal projection $P _{p}^{\star }$, and is guaranteed that there exist three vectors $\theta _{p,1},\theta _{p,2},\theta _{p,3}$ such that $\left\lVert P _{p}^{\star }\theta _{p,1}\right\rVert_{2}^{2}\gtrsim \frac{\sqrt{k ^{\star }}}{N}\sigma ^{2}, \left\lVert P _{p}^{\star }\theta _{p,2}\right\rVert_{2}^{2}\gtrsim \epsilon \sqrt{\frac{k _{p}^{\star }}{N}}\sigma ^{2},\left\lVert P _{p}^{\star }\theta _{p,3}\right\rVert_{2}^{2}\gtrsim \epsilon ^{2}\sigma ^{2}$, $\theta _{p,1}, \theta _{p,2},\theta _{p,3}$ are not necessary to be identical. In Theorem \ref{theorem: upper bound 1} or Theorem \ref{theorem: upper bound 2}, however, we require a single vector satisfying all the conditions simultaneously. In fact, similar to the $\ell _{2}$ norm case, we do need to select the minimal optimal dimension. To overcome the difficulty specified above, we again need to exploit the properties of the $p$-convex orthosymmetric set and the Kolmogorov widths. The proof of Theorem \ref{theorem: lp main upper bound} is provided in Appendix \ref{subsection: proof of lp main upper bound}.

\section{Experiments}\label{section: experiments}

In this section, we present the results of simulation to verify the correctness of the theorems above and algorithms in Appendix \ref{section: algorithms} on the application side. We generate artificial data, design different strategies for the adversary to contaminate the data, and run the (polynomial-time) algorithm to check its effectiveness. The pseudo code is attached in Appendix \ref{section: algorithms}.

\subsection{Data Generation}\label{section: setting of generation}

We generate artificial data using the function \texttt{numpy.random.multivariate} in the \texttt{Numpy} package in \texttt{Python}. To fully verify the correctness of the algorithms, we consider both the classical and high-dimensional regimes separately, and explore a variety of parameter configurations. Specifically, in the classical scenario, we let the parameter tuple $\left(N,d,\epsilon\right)$ iterate over all possible combinations of $\left\{200,2000\right\}\times\left\{5,10,20,50\right\}\times\left\{0.01,0.015,0.02,0.025\right\}$; in the high-dimensional case, we let the tuple iterate over all possible combinations of $\left\{100,500\right\}\times \left\{100,500\right\}\times \left\{0.01,0.015,0.02,0.025\right\}$. We set $\sigma =1$ in all generation procedures.

\subsection{Choice of Constraints}\label{subsection: setting of constraints}

For the choice of constraints, i.e., $K$, we select the common ellipsoidal constraint:
\begin{equation}\label{eq: ellipsoidal constraint}
    \mu \in K _{e}:=\left\{v \left\lvert\right. v \in \mathbb{R}^{d},\sum\limits_{i=1}^{n}\frac{v _{i}^{2}}{a _{i}}\le 1,a _{i}=\sqrt{d}\cdot i ^{-2}, 1 \le i \le d\right\}.
\end{equation}
Note that our choice of the axes of the ellipse ensures that the maximal $\ell _{2}$ norm of the vectors in $K _{e}$ is greater than the theoretical lower bound \eqref{eq: main theoretical lower bound}. For a given positive value $c \le \sqrt{a _{1}}$ as the $\ell _{2}$ norm of $\mu $ under $H _{1}$, the following approach uniquely determines the specific vector $\mu _{c}$ in $K _{e}$ such that $\left\lVert \mu _{c}\right\rVert_{2}^{}=c$.
\begin{equation}\label{eq: the specific choice of vector in the ellipse}
    \mu _{c}=\left\{
    \begin{tabular}{ll}
        $(\sqrt{\frac{c ^{2}/a _{k ^{\star }+1}-1}{1/a _{k ^{\star }+1}-1/a _{1}}},\underbrace{0,\dots,0}_{k ^{\star }-1 \text{ zeros}} ,\sqrt{\frac{1-c ^{2}/a _{1}}{1/a _{k ^{\star }+1}-1/a _{1}}},0,\dots,0)$, & if $c ^{2}>a _{k ^{\star }+1}$ and $k ^{\star }<d$,\\ 
        $(\underbrace{0,\dots,0}_{k ^{\star } \text{ zeros}},c,0,\dots,0)$, & if $c ^{2}\le a _{k ^{\star }+1}$ and $k ^{\star }<d$,\\ 
        $\left(c,0,\dots,0\right)$, & if $k ^{\star }=d$,
    \end{tabular}
    \right.
\end{equation}
where $k ^{\star }=\min\limits \left\{k _{1}^{\star },k _{2}^{\star }\right\}$. We verify in Appendix \ref{subsection: verification of the property of minimal} that the approach \eqref{eq: the specific choice of vector in the ellipse} ensures that the resulting vector $\mu _{c}$ achieves the minimal value of $\left\lVert P ^{\star }\mu \right\rVert_{2}^{}$ among the set $\left\{v \left\lvert\right. v \in \mathbb{R}^{d},v \in K _{e},\left\lVert v\right\rVert_{2}^{}=c\right\}$, where $P ^{\star }$ is the optimal projection operator corresponding to $k ^{\star }$. We believe that this property places $\mu _{c}$ in a near-worst-case position for the algorithm, and thus provides an appropriate setting to test its robustness.

\subsection{Configurations of Adversary}\label{section: setting of adversary}
For the adversary $\mathcal{C}$ in the strong $\epsilon $-contamination model, we design the following two strategies. Note that in the strategies, $\mu _{0}$ is an unknown fake mean vector determined by $\mathcal{C}$ that can depend on the constraint $K$.

\begin{algorithm}[htbp]
    \caption{strategy 1 of the adversary $\mathcal{C}$}\label{algorithm: strategy 1 of the adversary}
    
    \Import $\tilde{\mathbf{X}}, N,d,\sigma =1$

    \vspace*{1em}

    Compute $\texttt{arr}\in \mathbb{R}^{N}$, $\texttt{arr}_{i}=\left\lVert \tilde{X}_{i}\right\rVert_{2}^{},1 \le i \le N$

    \If {$H _{1}$ is true}{
        $C=\text{indices of the largest } \epsilon N \text{ entries in } \texttt{arr}$
        
        $\mathbf{X}_{C}=\epsilon N \text{ i.i.d. fake data from } \mathcal{N}\left(\mathbf{0},\mathbf{I}_{d}\right)$

        $\tilde{\mathbf{X}}_{C}=\mathbf{X}_{C}$
    }

    \Else {
        \Import $\mu _{0}$

        \vspace*{1em}

        $C=\text{indices of the smallest } \epsilon N \text{ entries in } \texttt{arr}$

        $\mathbf{X}_{C}=\epsilon N \text{ i.i.d. fake data from } \mathcal{N}\left(\mu _{0},\mathbf{I}_{d}\right)$

        $\tilde{\mathbf{X}}_{C}=\mathbf{X}_{C}$
    }

    \Return $\mathbf{X}=\mathcal{C}(\tilde{\mathbf{X}})$
\end{algorithm}

\begin{algorithm}[htbp]
    \caption{strategy 2 of the adversary $\mathcal{C}$}\label{algorithm: strategy 2 of the adversary}
    
    \Import $\tilde{\mathbf{X}}, N,d,\sigma =1$

    \vspace*{1em}

    \If {$H _{1}$ is true}{
        Compute $\texttt{dist}\in \mathbb{R}^{N}$, $\texttt{dist}_{i}=\left\lVert \tilde{X}_{i}\right\rVert_{2}^{},1 \le i \le N$

        $C=\text{indices of the largest } \epsilon N \text{ entries in } \texttt{dist}$
        
        $\mathbf{X}_{C}=\epsilon N \text{ i.i.d. fake data from } \mathcal{N}\left(\mathbf{0},\mathbf{I}_{d}\right)$

        $\tilde{\mathbf{X}}_{C}=\mathbf{X}_{C}$
    }

    \Else {
        Compute $\mu _{0}$ according to the constraint type

        Compute $\texttt{dist}\in \mathbb{R}^{N}$, $\texttt{dist}_{i}=\left\lVert \tilde{X}_{i}-\mu _{0}\right\rVert_{2}^{}$

        $C=\text{indices of the largest } \epsilon N \text{ entries in } \texttt{dist}$

        $\mathbf{X}_{C}=\epsilon N \text{ i.i.d. fake data from } \mathcal{N}\left(\mu _{0},\mathbf{I}_{d}\right)$

        $\tilde{\mathbf{X}}_{C}=\mathbf{X}_{C}$
    }

    \Return $\mathbf{X}=\mathcal{C}(\tilde{\mathbf{X}})$
\end{algorithm}

In strategies \ref{algorithm: strategy 1 of the adversary} and \ref{algorithm: strategy 2 of the adversary}, the adversary $\mathcal{C}$ and the resulting corruption are largely determined on the \texttt{arr} and \texttt{dist} arrays, which compute the distance between the original data and the origin, and between the original data and an imagined mean $\mu _{0}$, respectively. Both strategies were implemented in our experiments, and the results indicate that the choice has minimal impact on the outcomes.

\subsection{Experimental Results}\label{section: experimental results}

Experimental results for $N=200$ in the classical case and $N=100$ in the high-dimensional case are presented in the following figures. Results for other parameter configurations are similar and are deferred to Appendix \ref{section: algorithms}. In our experiments, we maintained the type-I error near zero while identifying the empirical $\tilde{\rho} _{\text{critical}}$ under $H _{1}$ as the critical value of $\rho $ when the corresponding type-II error exactly crosses $\beta =0.05$.

\begin{figure}[htbp]
    \centering
    \includegraphics[width=0.84\textwidth]{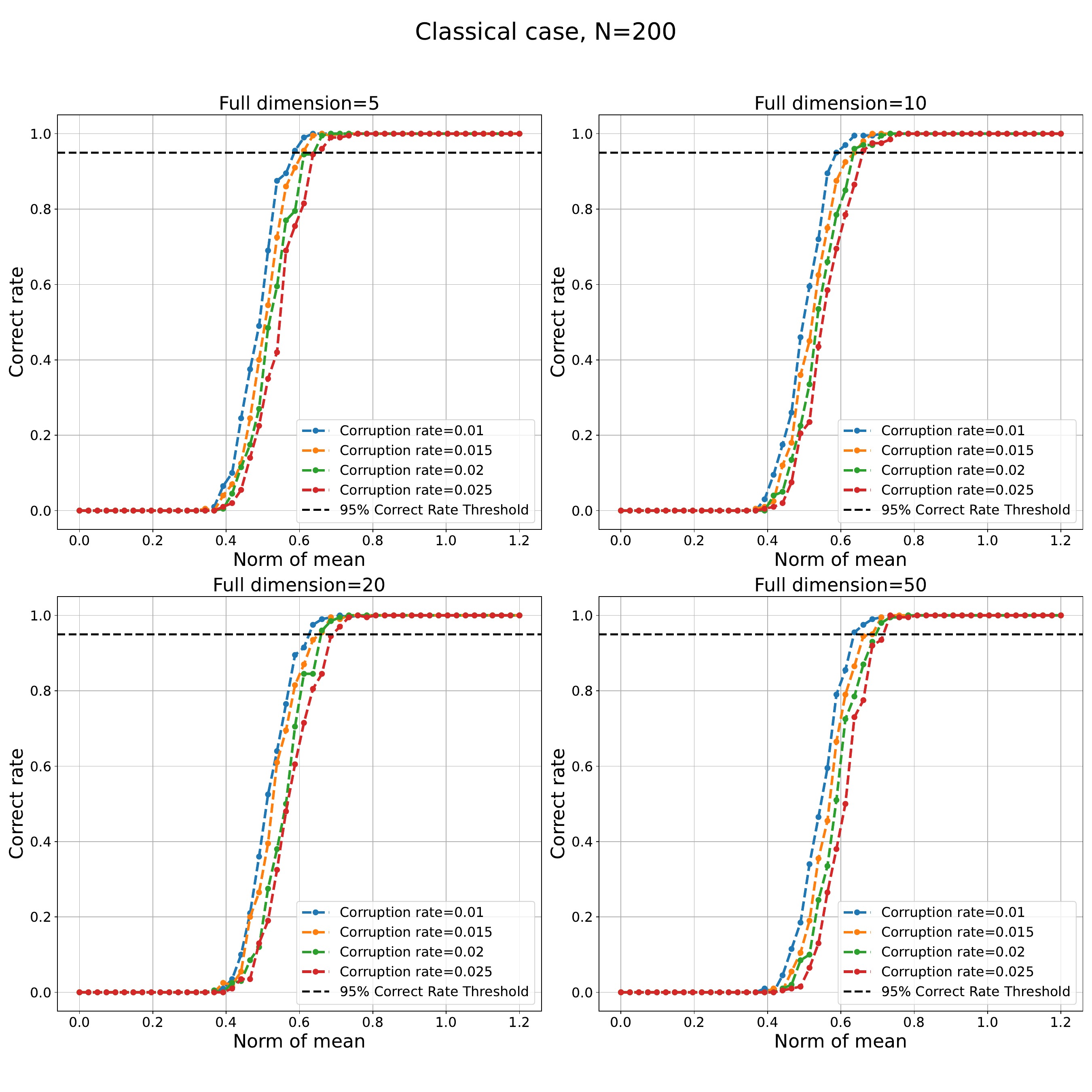}
    \includegraphics[width=0.84\textwidth]{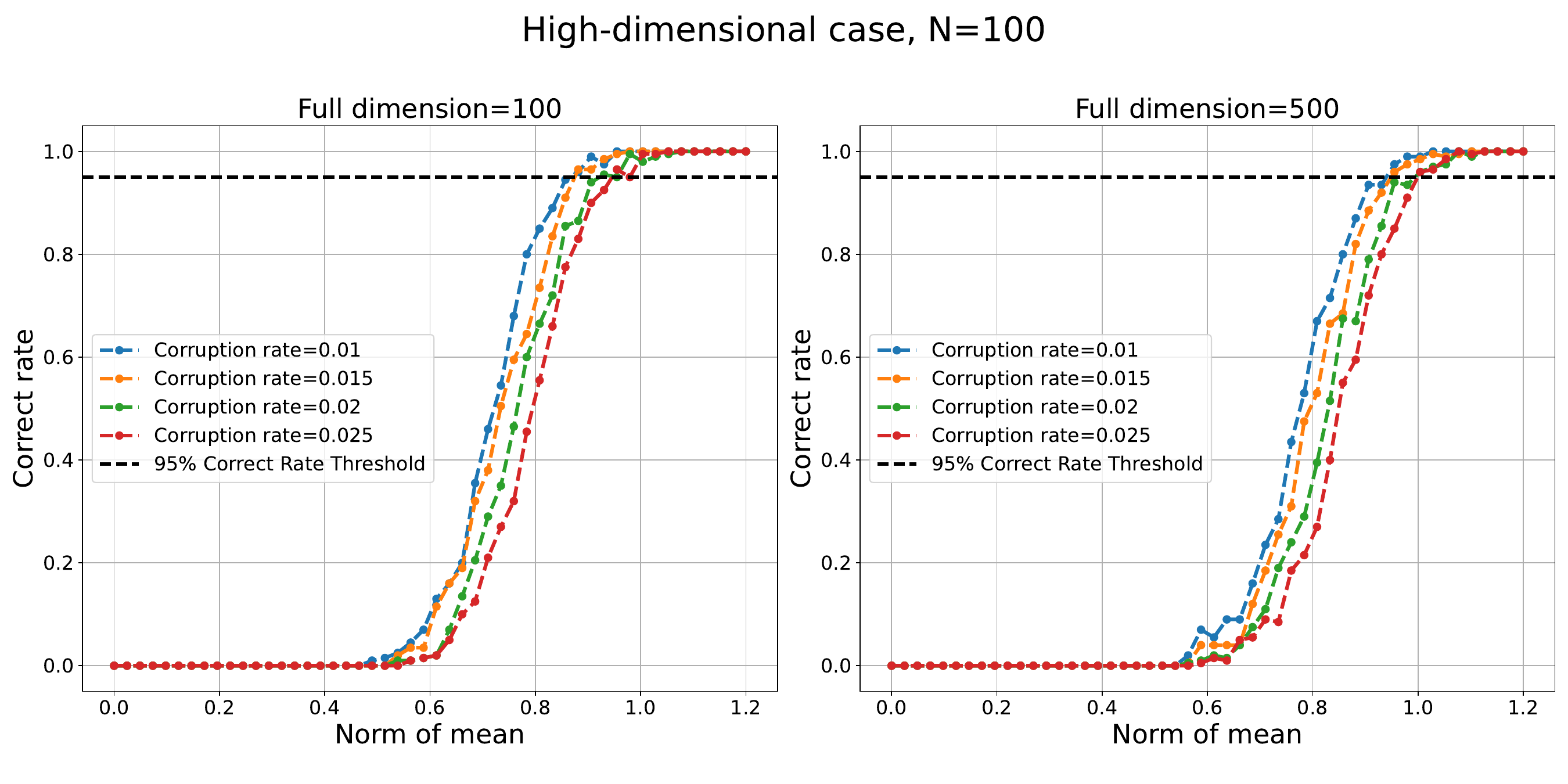}
    \caption{\textbf{Top}: Rejection rate (or Empirical power) versus the $\ell _{2}$ norm of the mean under $H _{1}$ in the classic setting ($N=200$). \textbf{Bottom}: Rejection rate versus the $\ell _{2}$ norm of the mean in the high-dimensional setting ($N=100$).}
\end{figure}

We observed several notable phenomena from the simulation results. First, the algorithm performs reasonably well across all parameter configurations, even for the moderate sample sizes: it successfully rejects $H _{0}$ when $\left\lVert \mu \right\rVert_{2}^{}$ is sufficiently large and fails to reject when the signal strength is too weak, in both the classical and the high-dimensional regimes. Second, focusing on $\rho _{\text{critical}}$ --- the critical norm of $\mu $ that yields a type-II error of $\alpha =0.05$ --- we find that $\rho _{\text{critical}}$ decreases as $N$ increases, and increases as $d$ grows. This behavior is consistent with theoretical intuition: larger sample sizes improve detectability, whereas higher dimensionality makes the testing problem more challenging. However, a more refined observation is that the influence of $d$ diminishes once the dimension becomes sufficiently large (e.g., around $d=50$ in our experiments). The reason behind this phenomenon is that, although $d$ increases to relatively large values, the optimal projection dimension $k ^{\star }$ becomes increasingly insensitive to $d$ under our setting of $N=1,\sigma =1$, and the ellipsoidal constraint. A direct calculation shows that when $N=2000$, the optimal dimensions corresponding to full dimensions $[5,10,20,50]$ are $[5,9,10,12]$, respectively. This demonstrates that the gaps between the optimal dimensions can be substantially smaller than those between the actual dimensions, which leads to the phenomenon. Finally, the plots indicate that as $N$ increases, the performance of the algorithm becomes less sensitive to variations in $\epsilon $. Though not fully verified, we believe that this is attributed to the fact that for larger $N$, the eigensubspace structure of $\mathbf{X}\mathbf{X}^{\top } $ (or $\mathbf{X}^{\top } \mathbf{X}$) becomes more stable with respect to the corrupted observations.

We also remark that we computed the ratio between $\tilde{\rho }_{\text{critical}}$ and the theoretical lower bound \eqref{eq: main theoretical lower bound} in an attempt to assess the tightness of the bound. Although the preliminary results appeared encouraging, we chose not to report them or draw conclusions. The primary reason is that both \eqref{eq: main theoretical lower bound} and its corresponding upper bound are derived under idealized minimax conditions that depend on the precise forms of $K$, $\mu $, and the adversary $\mathcal{C}$. In practical settings, however, these worst case quantities ($\mu $ and $\mathcal{C}$) are inherently unknown to us. Consequently, despite our simulations incorporating designed mean structures and adversarial strategies, they cannot guarantee that the conditions required for demonstrating optimality are satisfied. Identifying these conditions poses a fundamental challenge, and we leave this direction to future work.

For the hyperrectangle-type constraint, we conducted the experiments with the same settings and procedures. The results exhibited behavior similar to those obtained under the ellipsoidal constraint. Therefore, we present only the results under the ellipsoidal case.

\section{Discussion}\label{section: discussion}

\subsection{Rotated QCO Sets}\label{subsection: rotated QCO sets}

For completeness, we note that if the set $K$ is an orthogonal rotation of a QCO set --- i.e., $K = U K'$ for some orthogonal matrix $U$ --- then one may transform the data by considering $U ^{\top } X$ and subsequently solve the corresponding problem on the QCO set $K ^{\prime }$. This follows from the fact that $U ^{\top } \tilde{X}_{i} = U ^{\top }\mu + U ^{\top }\xi _{i}$, where $U ^{\top }\xi _{i}\sim N(0, \sigma^2 \mathbf{I})$ and $U ^{\top } \mu \in K'$. Moreover, the computational cost of determining the optimal dimensions for $K$ and $K ^{\prime }$ is essentially the same in most cases.

\subsection{Equivalence in Optimal Dimensions}\label{subsection: equivalence in optimal dimensions}

An alternative but equivalent way to define the optimal dimension is to simply select $k _{\text{alt}}^{\star }$ satisfying $D _{k _{\text{alt}}^{\star }-1}(K)>\left(\frac{k _{\text{alt}}^{\star }}{N}\right)^{1/4}\max\limits \left\{\sqrt{\epsilon },N ^{-\frac{1}{4}}\right\}\sigma $ and $D _{k _{\text{alt}}^{\star }}(K)\le \left(\frac{k _{\text{alt}}^{\star }+1}{N}\right)^{1/4}\max\limits \left\{\sqrt{\epsilon },N ^{-\frac{1}{4}}\right\}\sigma $. It also has the equivalent definitions as the discussion following Definition \ref{def_first_optimal_dimension} and \ref{def_second_optimal_dimension}. Once $k _{\text{alt}}^{\star }$ is determined, we simply project the data using $P ^{\star }_{\text{alt}}$, the optimal projection associated with $k _{\text{alt}}^{\star }$. The remainder of the testing procedure is identical to that in Algorithms \ref{algorithm: exponential algorithm} and \ref{algorithm: polynomial algorithm}. By this new definition, the lower bound will be adjusted accordingly as 
\begin{equation}\label{eq: alt main theoretical lower bound}
    \rho ^{2}\gtrsim \sigma ^{2}\max\limits \left\{\sqrt{\frac{k _{\text{alt}}^{\star }}{N}}\max\limits \left\{\epsilon ,\frac{1}{\sqrt{N}}\right\},\epsilon ^{2}\right\}.
\end{equation}
Similarly, the upper bounds $E _{\text{alt}}^{2}(\epsilon ,K,N,\sigma ^{2})$ and $P _{\text{alt}}^{2}(\epsilon ,K,N,\sigma ^{2})$ are obtained by substituting all $k _{1}^{\star }$ and $k _{2}^{\star }$ in $E ^{2}(\epsilon ,K,N,\sigma ^{2})$ and $P ^{2}(\epsilon ,K,N,\sigma ^{2})$ by $k _{\text{alt}}^{\star }$. With the proof in Appendix \ref{subsection: proof of the main corollary}, one can verify that $k _{\text{alt}}^{\star }, P _{\text{alt}}^{\star }$ and its corresponding testing process are equivalent to that of $k _{1}^{\star }$, $k _{2}^{\star }$, $P _{1}^{\star }$, and $P _{2}^{\star }$. The new lower bound \eqref{eq: alt main theoretical lower bound} and the upper bounds $E _{\text{alt}}^{2}(\epsilon ,K,N,\sigma ^{2})$ and $P ^{2}(\epsilon ,K,N,\sigma ^{2})$ do match the previous lower bound \eqref{eq: main theoretical lower bound} and the upper bounds \eqref{eq: alt required condition of the main upper bound 1}, \eqref{eq: alt required condition of the main upper bound 2} as well.

\subsection{Testing-Estimation Relationship}\label{subsection: testing estimation relationship}

The results developed in Section \ref{section: main result} should be compared and contrasted with the existing results for the corresponding estimation problem. Specifically, in the estimation setting --- where the same distribution assumptions and corruption model are imposed, but the goal is to estimate the mean vector $\mu $ --- the minimax rates are established by \cite{donoho1990minimax} and \cite{neykov2022minimax}. Although the authors do not explicitly employ the concepts of Kolmogorov widths and optimal dimensions, their results can nonetheless be interpreted within our framework. Notably, the Kolmogorov widths and the associated optimal dimensions again play a central role in these results.

\begin{definition}[Optimal dimension and projection in estimation]
    Given the same assumptions on the distribution, the constraint $K$, and the corruption model as \ref{def_first_optimal_dimension}, \ref{def_second_optimal_dimension}, and \ref{def_epsilon_contamination_model}, $k _{e}^{\star }$ is defined as the dimension that satisfies $D _{k _{e}^{\star }-1}(K)>\sqrt{\frac{k _{e}^{\star }}{N}}\sigma $ but $D _{k _{e}^{\star }}\le \sqrt{\frac{k _{e}^{\star }+1}{N}}\sigma $. It can be equivalently defined as 
    \begin{equation}\label{def_estimation_optimal_dimension}
        k _{e}^{\star }:=\max\limits _{k}\left\{k \left\lvert\right. k \ge 0, D _{k-1}(K)\ge \sqrt{\frac{k}{N}}\sigma \right\}.
    \end{equation}
\end{definition}

\begin{theorem}[Robust Estimation Rate, \cite{prasadan2025informationtheoreticlimitsrobust}]\label{theorem: lower bound of estimation}
    Given the same assumptions on the distribution, the constraint $K$, and the corruption model as \eqref{eq: testing problem}, the minimax rate for the estimation problem is given by
\begin{equation}\label{eq: lower bound of estimation}
    \inf\limits_{\hat{\nu }} \sup\limits_{\mu \in K} \mathbb{E}_\mu \left\lVert \hat{\nu }(\mathbf{X})-\mu \right\rVert_{2}^{2}\asymp \sigma ^{2}\max\limits \left\{\frac{k _{e}^{\star }}{N},\epsilon ^{2}\right\},
\end{equation}
where $\inf $ is taken with respect to all measurable functions $\hat \nu$ of the data.
\end{theorem}
When $K=\mathbb{R}^{d}$ and there is no corruption, we have $k _{e}^{\star }=d$ and $\epsilon =0$, in which case the rate above recovers the classic result that $\inf\limits_{\hat{\nu }}\sup\limits_{\mu }\mathbb{E}_{\mu }\left\lVert \hat{\nu }(\tilde{\mathbf{X}})-\mu \right\rVert_{2}^{2}\asymp \frac{d}{N}\sigma ^{2}$. It is also well known that, in this setting, the estimation problem is generally more difficult than the testing problem. The reason is simple and intuitive --- one is able to derive a valid test with small error probabilities whenever $\left\lVert \mu \right\rVert_{2}^{2}$ is asymptotically larger than the minimax lower bound for the estimation task. This relationship, unsurprisingly, is inherited by our more general problem.

\begin{lemma}[Testing-estimation relationship]\label{lemma: testing estimation relationship}
    Given the same assumptions on the distribution, the constraint $K$, and the corruption model above, we have
    \begin{equation}
        \sigma ^{2}\max\limits \left\{\frac{\sqrt{k _{1}^{\star }}}{N},\epsilon \sqrt{\frac{k _{2}^{\star }}{N}},\epsilon ^{2}\right\}\lesssim \sigma ^{2}\max\limits \left\{\frac{k _{e}^{\star }}{N},\epsilon ^{2}\right\}
    \end{equation}
\end{lemma}

We refer readers to Appendix \ref{subsection: proof of the testing estimation relationship} for a straightforward computational verification of Lemma \ref{lemma: testing estimation relationship}. Following the discussion after Corollary \ref{corollary: matching between bounds}, we can also analyze the phase transition of the estimation problem. It occurs around $\epsilon \asymp \sqrt{\frac{k _{e}^{\star }}{N}}$. If $\epsilon \lesssim \sqrt{\frac{k _{e}^{\star }}{N}}$, then $\epsilon $ is too small to influence the minimax rate of the estimation error. Otherwise (for example, $\epsilon \asymp N ^{-\frac{1}{4}}$), the adversary is strong enough to the extent that $N$ is excluded and the minimax rate only relies on $\epsilon $ and $\sigma $.

Another interesting observation from Lemma \ref{lemma: testing estimation relationship} is that although estimation is generally more challenging than testing, the gap disappears as the corruption rate increases. In fact, both quantities stabilize at $\epsilon ^{2}\sigma ^{2}$ when $\epsilon $ is sufficient large. Similar phenomenon has also been observed in a recent work \cite{kania2025testingimprecisehypotheses}. See Section \ref{subsection: open problems and future work} for additional details.

\subsection{Open Problems and Future Work}\label{subsection: open problems and future work}

An intriguing open direction is to investigate the localized minimax testing rates for QCO sets. A preliminary step in this direction was taken by \cite{wei2020local}, who studied local testing rates in the setting of ellipsoids. It would be interesting to understand the extent to which their results can be generalized to QCO sets. Moreover, from an intuitive standpoint, the origin appears to be the most challenging point for testing under QCO constraints, as it may correspond to the largest local volume relative to other points. This raises a natural question: does the global minimax testing rate coincide with the localized rate at the origin, and thus with the rates established in this paper?

A fundamental assumption of this paper is that the covariance matrix is fully known (and taken to be $\mathbf{I}_{d}$). This, however, is rarely guaranteed in real-world applications. An important direction for future work is therefore to understand whether --- and to what extent --- we can still distinguish between $H _{0}$ and $H _{1}$ when $\bm{\Sigma }$ must be estimated rather than treated as prior knowledge. A recent relevant work in robust covariance matrix estimation is \cite{robustgaussiancovarianceestimation}, where the authors proposed an algorithm achieving a near-optimal error rate $\mathcal{O}\left(\epsilon \ln \left(\frac{1}{\epsilon }\right)\right)$ in the Mahalanobis norm under $\epsilon $-strong contamination model while remaining computationally efficient. Nevertheless, this algorithm is designed for zero-mean distributions and unconstrained setting. It remains unknown whether reliable testing is achieveble and how the constraint $K$ can facilitate the problem when $\bm{\Sigma }$ is unknown.
  
In \cite{kania2025testingimprecisehypotheses}, the authors consider a more general problem than \eqref{eq: original testing problem}:
\begin{equation}\label{eq:general_testing_problem}
    \begin{aligned}
        & H _{0}: \left\lVert \mu \right\rVert_{2}^{}\le \rho _{0},\\ 
        & H _{1}: \left\lVert \mu \right\rVert_{2}^{}\ge \rho _{1}.
    \end{aligned}
\end{equation}
In other words, rather than being fixed at the origin, the mean vector under the null is allowed to vary within a neighborhood of the origin. An interesting corollary of their work is that as the neighborhood expands, the minimax testing rate converges to the minimax estimation rate --- an effect that parallels a key phenomenon in our setting. Such match in estimation and testing rate is also observed by other work (e.g., \cite{DBLP:conf/alt/Ndaoud19}). Given the framework of their problem, it is natural to ask whether one can establish the minimax rate with constrained and adversarial setting of \eqref{eq:general_testing_problem}. From the perspective of the testing-estimation relationship, an even more important question is: under what conditions does the gap between testing and estimation disappear?

Another important avenue for future research will be to study $\ell _{p}$ norm separation for $p>2$. We believe there are fundamental changes compared to the problem discussed in this work. The first reason is that $\chi ^{2}$-test is no longer suitable for the testing problem when $p>2$. Instead, $\chi ^{p}$-test is a useful tool in this scenario, see the results by \cite{ingster2003nonparametric} when there is no constraint and adversary. The second reason is that the geometric properties of the $\ell _{p}$ balls when $p>2$ are known to be dramatically different compared to those when $1 \le p \le 2$. This affects the selection of the extreme vector as in Lemma \ref{lemma: existence of the optimal vector of a QCO set}. It is also unclear at the moment whether the Kolmogorov widths will still drive the minimax rates in that case.

Finally, a more ambitious open problem is to characterize the minimax testing rates over general convex sets that contain the origin. It would be particularly exciting to determine whether an optimal test exists for every such convex set $K$, analogous to the result of \cite{neykov2022minimax}, which established the existence of a universally optimal estimation procedure for arbitrary convex sets in the corresponding estimation problem.

\section*{Acknowledgments}
The authors are grateful to Siva Balakrishnan for discussions on local testing on ellipsoids and related minimax testing problems. The authors also would like to thank Cl\'ement Canonne for a helpful email exchange in which Cl\'ement suggested to us that in the unconstrained problem, $\left\lVert \mu \right\rVert_{2}^{}\gtrsim \epsilon \sigma $ is required, which we extended to show Theorem \ref{theorem: lower bound 3}.

\newpage

\bibliographystyle{apalike}
\bibliography{biblio}

\newpage

\appendix

\section{Theoretical and Polynomial-Time Algorithms}\label{section: algorithms}

\subsection{Theoretical Algorithm}

According to Theorem \ref{theorem: upper bound 1} and Section \ref{subsection: proof of the main upper bound 1}, we first attempt to search for a consistent set $T$ with respect to $\phi _{e}$. $\phi _{e}(\mathbf{X}_{T})$ then determines whether we are to accept or reject $H _{0}$.

\begin{algorithm}[H]
    \caption{theoretical algorithm when $\rho \gtrsim E(\epsilon ,K,N,\sigma ^{2})$}
    \label{algorithm: exponential algorithm}

    \Import $\mathbf{X},N,K,\sigma $

    \vspace*{1em}

    Compute $k _{1}^{\star }, P _{1}^{\star }, k _{2}^{\star }, P _{2}^{\star }$ from $N,K,\sigma ,\epsilon $ as Definition \ref{def_first_optimal_dimension} and \ref{def_second_optimal_dimension}

    \If{$k _{1}^{\star }\le k _{2}^{\star }$}{
        $k ^{\star }= k _{1}^{\star }, P ^{\star }=P _{1}^{\star }$
    }

    \Else{
        $k ^{\star }=k _{2}^{\star }, P ^{\star }=P _{2}^{\star }$
    }

    $\mathbf{Y}=\mathbf{X}P ^{\star }$

    \For{$T$ in $\left\{\text{all subsets with cardinality greater than }(1-\epsilon )N\right\}$}{

        \texttt{consistent\_flag} $=$ \texttt{True}

        \If{$\left\lVert \sum\limits_{i \in T}^{}Y _{i}\right\rVert_{2}^{2}-k ^{\star }\left|T\right|\sigma ^{2}\ge c \left|T\right|^{2}E ^{2}(\epsilon ,K,N,\sigma )$}{
            \texttt{subset\_test} $=$ \texttt{True}
        }
        \Else{
            \texttt{subset\_test} $=$ \texttt{False}
        }
        \For{$T ^{\prime }$ in $\left\{\text{all subsets of $T$ with cardinality greater than }(1-2\epsilon )N\right\}$}{

            \If{$\left\lVert \sum\limits_{i \in T ^{\prime }}^{}Y _{i}\right\rVert_{2}^{2}-k ^{\star }\left|T ^{\prime }\right|\sigma ^{2}\ge c \left|T ^{\prime }\right|^{2}E ^{2}(\epsilon ,K,N,\sigma )$}{
                \texttt{subsubset\_test} $=$ \texttt{True}
            }
            \Else{
                \texttt{subsubset\_test} $=$ \texttt{False}
            }
            \If{\texttt{subset\_test} XOR \texttt{subsubset\_test}}{
                \texttt{consistent\_flag} $=$ \texttt{False}

                \Break
            }
        }
        \If{\texttt{consistent\_flag}}{
            \If{\texttt{subset\_test}}{
                Reject $H _{0}$

                \Return
            }
            \Else{
                Accept $H _{0}$

                \Return
            }
        }
    }
    Accept $H _{0}$

    \Return
\end{algorithm}

\subsection{Polynomial-Time Algorithm}

The polynomial-time algorithm combines the three sub-algorithms in Section \ref{subsection: proof of the main upper bound 2} together to filter the observations $\mathbf{X}$ and through the weight $\omega $. The theoretical guarantee is provided in Section \ref{subsection: proof of the main upper bound 2}.

\begin{algorithm*}[htbp]
    \caption{Prefiltering.}\label{algorithm: prefiltering init}
    Set $\gamma _{1}=c\left[\sqrt{k ^{\star }\ln \left(\frac{N}{\alpha }\right)}+\ln \left(\frac{N}{\alpha }\right)\right]$, $\texttt{count}=0,\texttt{i}=0$

    \While {$\texttt{i} < N$}{
        \If {$\left|\left\lVert Y _{i}\right\rVert_{2}^{2}-k ^{\star }\right|> \gamma _{1}$}{

            $\texttt{count}=\texttt{count}+1$

            \If {$\texttt{count}>\epsilon N$}{
                \Return \texttt{None}
            } 
            Delete $Y _{i}$ from $\mathbf{Y}$      
        }
        $\texttt{i}=\texttt{i}+1$
    }

    \Return $\mathbf{Y}$
\end{algorithm*}

\begin{algorithm}[htbp]
    \caption{Sample filtering when $N>k ^{\star }$.}\label{algorithm n > k init}
    Set $\gamma _{2}$ as \eqref{eq: definition of gamma_2}, and $\lambda =\left\lVert Y ^{\top } D(\omega )Y-N \mathbf{I}_{k ^{\star }}\right\rVert_{2}^{}$. ($\omega $ is initialized as $\mathbf{1}$.)

    \While {$\lambda \ge \gamma _{2}$}{
        Set $v$ to be the unit singular vector associated with $\lambda $

        Compute $\tau _{i}=\left\langle v, Y _{i} \right\rangle ^{2}\mathbf{1}_{\left\{\omega _{i}>0\right\}}$ for $1 \le i \le N$

        Set $I$ be the smallest index such that $\sum\limits_{i=1}^{I}\omega _{i}\ge 2\epsilon N$

        Update $w$ according to \eqref{eq: update of omega}

        \If {$\left\lVert \omega \right\rVert_{1}^{}<N(1-2\epsilon )$} {

            \Return \texttt{None}
        }

        Set $\lambda =\left\lVert Y ^{\top } D(\omega )Y-N \mathbf{I}_{k ^{\star }}\right\rVert_{2}^{}$
    }

    \Return $\omega $
\end{algorithm}

\begin{algorithm}[htbp]
    \caption{Sample filtering when $N \le k ^{\star }$.}\label{algorithm n <= k init}
    Set $\gamma _{3}$ as \eqref{eq: definition of gamma_3}, and $\lambda =\left\lVert \sqrt{D(\omega)}\mathbf{Y}\mathbf{Y} ^{\top } \sqrt{D(\omega)}-D(\omega )\right\rVert_{2}^{}$. ($\omega $ is initialized as $\mathbf{1}$.)

    \While {$\lambda \ge \gamma _{3}$}{
        Set $v$ to be the unit singular vector associated with $\lambda $

        Compute $\tau _{i}=\frac{v _{i}^{2}}{\omega _{i}}\mathbf{1}_{\left\{w _{i}>0\right\}}$

        Update $w$ according to \eqref{eq: update of omega 2}

        \If {$\left\lVert \omega \right\rVert_{1}^{}<N(1-6\epsilon )$} {

            \Return \texttt{None}
        }

        Set $\lambda =\left\lVert \sqrt{D(\omega)}\mathbf{Y}\mathbf{Y} ^{\top } \sqrt{D(\omega)}-D(\omega )\right\rVert_{2}^{}$
    }

    \Return $\omega $
\end{algorithm}

\begin{algorithm}[htbp]
    \caption{Weight filtering.}\label{algorithm: weight filtering init}
    Compute $\tau _{i}=\left|\left\langle \sqrt{\omega _{i}}Y _{i}, \sum\limits_{j=1}^{N}\sqrt{\omega _{j}}Y _{j} \right\rangle-\omega _{i}k ^{\star }\right|\mathbf{1}_{\left\{\omega _{i}>0\right\}}, 1 \le i \le N$

    Sort $\tau _{i}$ by desreasing order and find the indices $\left\{i _{1},i _{2},\dots,i _{\epsilon N}\right\}$ corresponding to the first $\epsilon N$ maximal $\tau _{i}$

    Set $\omega _{i _{1}},\omega _{i _{2}},\dots,\omega _{i _{\epsilon N}}$ to zero

    \Return $\omega $
\end{algorithm}

\begin{algorithm}[htbp]
    \caption{polynomial-time algorithm when $\rho \gtrsim P(\epsilon ,k ^{\star },N,1)$.}\label{algorithm: polynomial algorithm}
    
    \Import $\mathbf{X}, N,K,\sigma =1$

    \Import Algorithm \ref{algorithm: prefiltering init}, \ref{algorithm n > k init}, \ref{algorithm n <= k init}, \ref{algorithm: weight filtering init}

    \vspace*{1em}

    $\mathbf{X}=\mathbf{X}/\sigma $

    Compute $k _{1}^{\star }, P _{1}^{\star }, k _{2}^{\star }, P _{2}^{\star }$ from $N,K,\sigma =1,\epsilon $ as Definition \ref{def_first_optimal_dimension} and \ref{def_second_optimal_dimension}

    \If{$k _{1}^{\star }\le k _{2}^{\star }$}{
        $k ^{\star }= k _{1}^{\star }, P ^{\star }=P _{1}^{\star }$
    }

    \Else{
        $k ^{\star }=k _{2}^{\star }, P ^{\star }=P _{2}^{\star }$
    }

    $\mathbf{Y}=\mathbf{X}P ^{\star }$

    Execute Algorithm \ref{algorithm: prefiltering init}, and record the returning value $\texttt{R} _{1}$

    \If {$\texttt{R} _{1}$ is \texttt{None}}{
        Reject $H _{0}$
        
        \Return
    }

    \Else {
        $\mathbf{Y}=\texttt{R} _{1}$
    }

    Set $\omega =\mathbf{1}_{[k ^{\star }]}$

    \If {$N>k ^{\star }$}{
        Execute Algorithm \ref{algorithm n > k init}, and record the returning value $\texttt{R}_{2}$
    }
    \Else {
        Execute Algorithm \ref{algorithm n <= k init}, and record the returning value $\texttt{R} _{2}$
    }
    \If {$\texttt{R}_{2}$ is \texttt{None}} {
        Reject $H _{0}$

        \Return
    }

    \Else {
        $\omega =\texttt{R} _{2}$
    }

    Execute Algorithm \ref{algorithm: weight filtering init}, and record the returning value $\texttt{R} _{3}$

    $\omega =\texttt{R} _{3}$

    \If {$\left|\left\lVert \sum\limits_{i=1}^{N}\sqrt{\omega _{i}}Y _{i}\right\rVert_{2}^{2}-k ^{\star } \left\lVert \omega \right\rVert_{1}^{}\right|\ge c _{2} N ^{2}P(\epsilon ,k ^{\star },N,1)$}{
        Reject $H _{0}$
    }

    \Else {
        Accept $H _{0}$
    }

    \Return
\end{algorithm}

\newpage

\subsection{Deferred Experiment Results}

The following figures display additional results deferred from Section \ref{section: experimental results}.

\begin{figure}[htbp]
    \centering
    \includegraphics[width=0.84\textwidth]{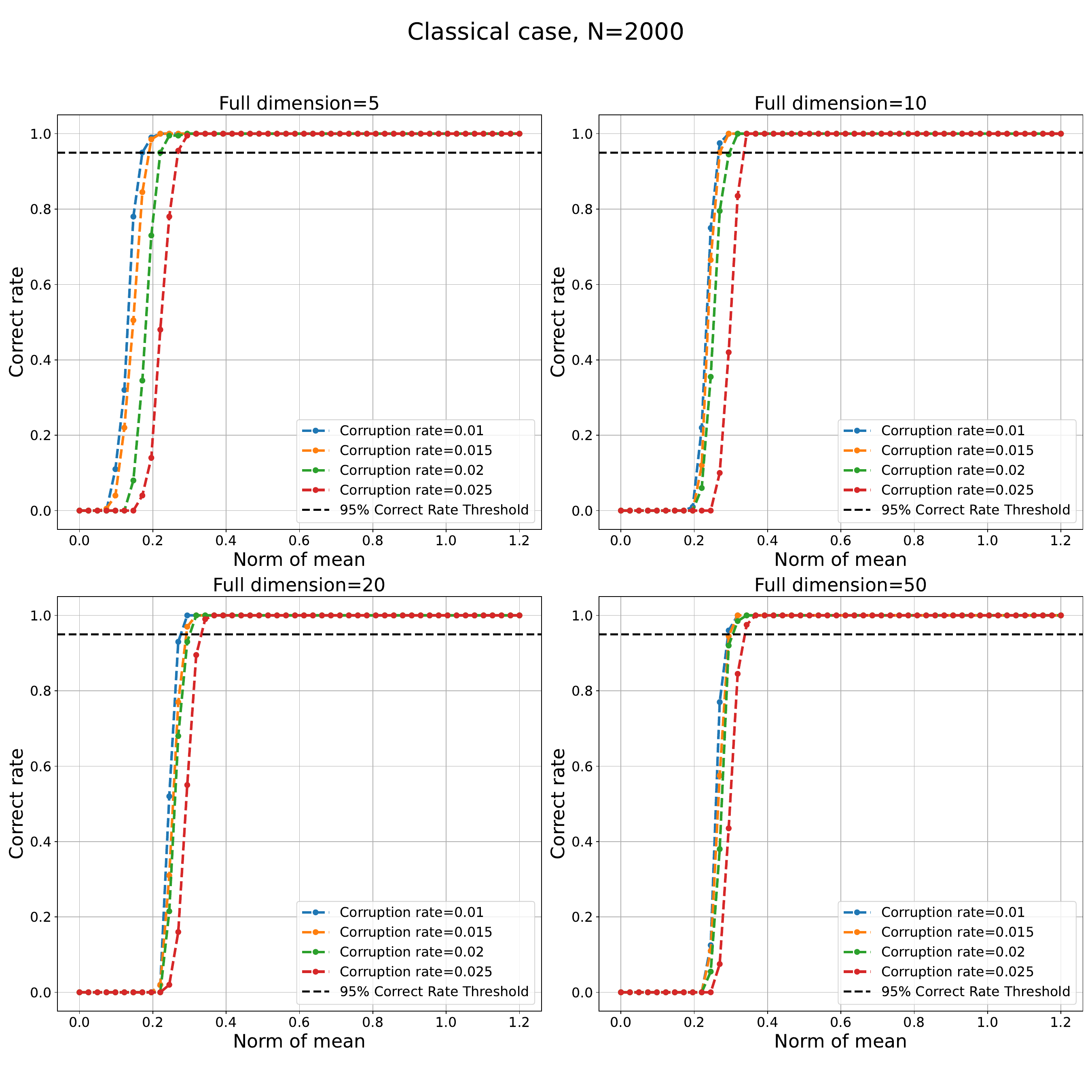}
    \includegraphics[width=0.84\textwidth]{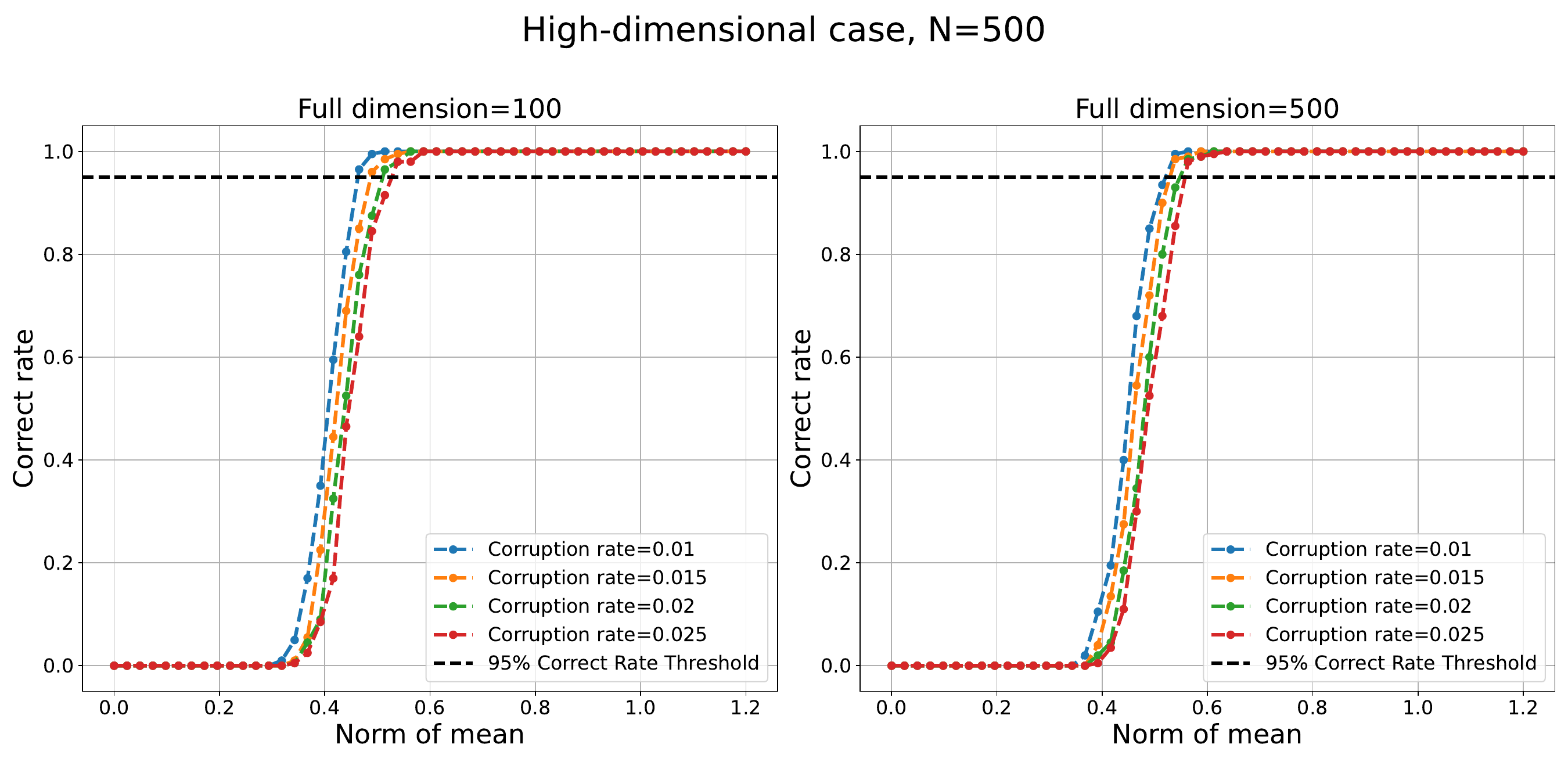}
    \caption{\textbf{Top}: Rejection rate (or Empirical power) versus the $\ell _{2}$ norm of the mean under $H _{1}$ in the classic setting ($N=2000$). \textbf{Bottom}: Rejection rate versus the $\ell _{2}$ norm of the mean in the high-dimensional setting ($N=500$).}
\end{figure}

\newpage

\section{Some Basic Probabilistic Concepts and Useful Lemmas}\label{section: concepts and lemmas}

In this section, we recall some basic probabilistic concepts and useful lemmas that are frequently mentioned and leveraged in this paper. Most of them are related to concentrations and the measures of distance between distributions.

\begin{definition}[Total variation distance]\label{def_total_variance}
    For a $\sigma $-algebra $(\Omega ,\mathscr{F})$ and the probabilities $\mathcal{P},\mathcal{Q}$ defined on $\mathscr{F}$, the total variation distance between $\mathcal{P}$ and $\mathcal{Q}$ is defined as 
    \begin{equation*}
        \text{TV}(\mathcal{P},\mathcal{Q})=\sup\limits_{E \in \mathscr{F}}\left|\mathcal{P}(E)-\mathcal{Q}(E)\right|.
    \end{equation*}
\end{definition}

From the definition, we immediately know that $\text{TV}(\mathcal{P},\mathcal{Q})=\text{TV}(\mathcal{Q},\mathcal{P})$ and $0 \le \text{TV}(\mathcal{P},\mathcal{Q})\le 1$. $\text{TV}(\mathcal{P},\mathcal{Q})=0$ if and only if $\mathcal{P}=\mathcal{Q}$ with only exception of zero-measure sets. When $\mathcal{P},\mathcal{Q}$ has densities $p,q$ with respect to some measure $\nu $ and space $\Omega _{\nu }$, it can be also easily shown that 
\begin{equation}\label{eq: property of tv distance}
    \text{TV}(\mathcal{P},\mathcal{Q})=\displaystyle\int_{\Omega _{\nu }}^{}\min\limits \left\{p,q\right\}d\nu =\frac{1}{2}\displaystyle\int_{\Omega _{\nu }}^{}\left|p-q\right|d\nu .
\end{equation}

\begin{definition}[the Kullback-Leibler divergence (KL divergence)]\label{def_KL_divergence}
    Given the same $\sigma $-algebra and probabilities, additionally assume that $\mathcal{P}\ll \mathcal{Q}$. The KL divergence between $\mathcal{P}$ and $\mathcal{Q}$ is defined as 
    \begin{equation*}
        \text{KL}(\mathcal{P},\mathcal{Q})=\displaystyle\int_{\Omega }^{}\ln \frac{d \mathcal{P}}{d \mathcal{Q}}d \mathcal{P}.
    \end{equation*}
\end{definition}

KL divergence is always non-negative, with $\text{KL}(\mathcal{P},\mathcal{Q})=0$ implying $\mathcal{P}=\mathcal{Q}$ with only exception of zero-measure sets. Contrasted with total variation distance, KL divergence is not symmetric, meaning that $\text{KL}(\mathcal{P},\mathcal{Q})\neq \text{KL}(\mathcal{Q},\mathcal{P})$ in general, even if $\frac{d \mathcal{P}}{d \mathcal{Q}}$ and $\frac{d \mathcal{Q}}{d \mathcal{P}}$ both exist.

\begin{lemma}[Pinsker's inequality]\label{lemma: Pinsker's inequality}
    For the probability distributions $\mathcal{P},\mathcal{Q}$ with $\mathcal{P}\ll \mathcal{Q}$, we have 
    \begin{equation*}
        \text{TV}(\mathcal{P},\mathcal{Q})\le \sqrt{\frac{1}{2}\text{KL}(\mathcal{P},\mathcal{Q})}.
    \end{equation*}
    When $\mathcal{P}\ll \mathcal{Q}$ and $\mathcal{Q}\ll \mathcal{P}$, it can be refined as 
    \begin{equation*}
        \text{TV}(\mathcal{P},\mathcal{Q})\le \sqrt{\frac{1}{2} \min\limits \left\{\text{KL}(\mathcal{P},\mathcal{Q}),\text{KL}(\mathcal{Q},\mathcal{P})\right\}}.
    \end{equation*}
\end{lemma}

There are many proofs for this lemma. One can refer \cite{tsybakov2009introduction} for a quick proof.

\begin{definition}[Coupling]\label{def_coupling}
    Given probability distributions $\mathcal{P},\mathcal{Q}$, a valid coupling of $\mathcal{P}$ and $\mathcal{Q}$ is any pair of random variables $(X,Y)$ where the marginal distribution of $X$ ang $Y$ is $\mathcal{P},\mathcal{Q}$ respectively.
\end{definition}

The following classic lemma help us connect the total variation distance with the Hamming distance in the \hyperref[subsection: proof of the main lower bound 2]{proof} of Theorem \ref{theorem: lower bound 2}.

\begin{lemma}[]\label{lemma: total variation and coupling}
    For any probability distributions $\mathcal{P}$ and $\mathcal{Q}$, their total variation distance equals to the infimum over the probability that $X$ is not equal to $Y$, where $(X,Y)$ is a valid coupling of $\mathcal{P}$ and $\mathcal{Q}$.
    \begin{equation}
        \text{TV}(P,Q)=\inf\limits_{\text{$(X,Y)$ is a valid coupling of P and Q}}\mathbb{P}(X\ne Y).
    \end{equation}
    The coupling that achieves the equality above is referred as the \textit{optimal coupling} between $\mathcal{P}$ and $\mathcal{Q}$.
\end{lemma}

The following definitions and lemmas concern the sub-Gaussian and sub-exponential random variables, which are closely related to the concentration bounds in the proof of the main theorems.

\begin{definition}[Sub-Gaussian random variables]\label{def_sub_Gaussian_var}
    For a univariate random variable $X$, define its sub-Gaussian norm as $\left\lVert X\right\rVert_{\psi _{2}}^{}:=\inf\limits_{}\left\{t:t>0,\mathbb{E}\exp\left\{\frac{X ^{2}}{t ^{2}}\right\}\le 2\right\}$. If such $\left\lVert X\right\rVert_{\psi _{2}}^{}$ does not exist, then define $\left\lVert X\right\rVert_{\psi _{2}}^{}=\infty $. $X$ is said to be sub-Gaussian if $\left\lVert X\right\rVert_{\psi _{2}}^{}<\infty $.
\end{definition}

There are many equivalent definitions for sub-Gaussian random variables. The most frequently mentioned one is $\mathbb{P}(\left|X\right|\ge t)\le 2 \exp\left\{-\frac{t ^{2}}{c \left\lVert X\right\rVert_{\psi _{2}}^{2}}\right\}$.

\begin{definition}[Sub-exponential random variables]\label{def_sub_exp_var}
    For a univariate random variable $X$, define its sub-exponential norm as $\left\lVert X\right\rVert_{\psi _{1}}^{}:=\inf\limits_{}\left\{t:t>0,\mathbb{E}\exp\left\{\frac{\left|X\right|}{t}\right\}\le 2\right\}$. If such $\left\lVert X\right\rVert_{\psi _{1}}^{}$ does not exist, then define $\left\lVert X\right\rVert_{\psi _{1}}^{}=\infty $. $X$ is said to be sub-exponential if $\left\lVert X\right\rVert_{\psi _{1}}^{}<\infty $.
\end{definition}

There are also many equivalent definitions for sub-exponential random variables. A popular one is $\mathbb{P}(\left|X-\mathbb{E}X\right|\ge t)\le 2 \exp\left\{-c\min\limits \left\{\frac{t ^{2}}{\left\lVert X\right\rVert_{\psi _{1}}^{2}},\frac{t}{\left\lVert X\right\rVert_{\psi _{1}}^{}}\right\}\right\}$.

The requirement of sub-exponential random variables is weaker than that of sub-Gaussian in the sense that if $X$ is sub-Gaussian, it must be sub-exponential. Moreover, if $X$ is sub-Gaussian, then by definition $X ^{2}$ is sub-exponential.

The following is the famous Bernstein's inequality, which handles multiple independent sub-Gaussian (sub-exponential) random variables. It also generalizes the equivalent definitions mentioned above.

\begin{lemma}[Bernstein's inequality]\label{lemma: Bernstein's inequality}
    Assume $X _{1},\dots,X _{n}$ are independent, zero-mean random variables. If $X _{1},\dots,X _{n}$ are sub-Gaussian, then for any $t \ge 0$, we have 
    \begin{equation}
        \mathbb{P}\left(\sum\limits_{i=1}^{n}X _{i}\ge t\right)\le \exp\left\{-c t ^{2}/\sum\limits_{i=1}^{n}\left\lVert X _{i}\right\rVert_{\psi _{2}}^{2}\right\}.
    \end{equation}
\end{lemma}

If $X _{1},\dots,X _{n}$ are sub-exponential, then we have 
\begin{equation}\label{eq: tail bound of sub-exp var}
    \mathbb{P}\left(\sum\limits_{i=1}^{n}X _{i}\ge t\right)\le \exp\left\{-c \min\limits \left\{t ^{2}/\sum\limits_{i=1}^{n}\left\lVert X _{i}\right\rVert_{\psi _{1}}^{2},t/\max\limits _{i}\left\{\left\lVert X _{i}\right\rVert_{\psi _{1}}^{}\right\}\right\}\right\}.
\end{equation}

The next lemma is particularly useful when we control the $\ell _{2}$ norm of a random vector with independent sub-Gaussian coordinates.

\begin{lemma}[]\label{lemma: sub_Gaussian_norm_of_l2_norm}
    Assume $X=(X _{1},\dots,X _{n})\in \mathbb{R}^{d}$ where $X _{1},\dots,X _{n}$ are independent, $\mathbb{E}X _{i}=0, \text{Var}\left(X _{i}\right)=1$, and each $X _{i}$ is sub-Gaussian. Then $\left\lVert X\right\rVert_{2}^{}-\sqrt{n}$ is sub-Gaussian with $\left\lVert \left\lVert X\right\rVert_{2}^{}-\sqrt{n}\right\rVert_{\psi _{2}}^{}\le c \max\limits _{i}\left\{\left\lVert X _{i}\right\rVert_{\psi _{2}}^{2}\right\}$.
\end{lemma}

The definitions and lemmas above can be found in many modern statistics textbook. For an example, see \cite{hdp}.

The famous Hanson--Wright inequality is a powerful tool to bound the tail probability of quadratic forms of the sub-Gaussian random variables. It can also be viewed as a generalization of the Bernstein's inequality \eqref{eq: tail bound of sub-exp var}. See \cite{10.1214/aoms/1177693335} for the initial paper and \cite{rudelson2013hansonwrightinequalitysubgaussianconcentration} for a modern proof.

\begin{lemma}[Hanson--Wright inequality]\label{lemma: hw_inequality}
    Assume that $X=(X _{1},\dots,X _{n})\in \mathbb{R}^{n}$ is a random vector with independent, zero-mean and sub-Gaussian coordinates. Let $A$ be an $n \times n$ matrix and $K=\max \limits _{i}\left\lVert X _{i}\right\rVert_{\psi _{2}}^{}$. Then for any $t \ge 0$, we have 
     \begin{equation}
        \mathbb{P}\left(|X ^{\top } AX-\mathbb{E}X ^{\top } AX|>t\right)\le 2 \exp\left\{-c \min\limits \left\{\frac{t ^{2}}{K ^{4}\left\lVert A\right\rVert_{F}^{2}},\frac{t}{K ^{2}\left\lVert A\right\rVert_{2}^{}}\right\}\right\},
     \end{equation}
     where $\left\lVert \cdot \right\rVert_{F}^{}$ is the Frobenius norm of a matrix, and $\left\lVert \cdot \right\rVert_{2}^{}$ is the $\ell _{2}$ operator norm of a matrix.
\end{lemma}

We now turn to the Wishart matrix, which plays a central role in covariance estimation in high-dimensional Gaussian models. Additional properties of the Wishart distribution and its applications in estimation can be found in \cite{KUBOKAWA20081906} and \cite{Bishop_2018}.

\begin{definition}[Wishart matrix]
    Given a matrix $X \in \mathbb{R}^{n \times d}$ where each row is independently drawn from $\mathcal{N}(0,\bm{\Sigma })$ for some $\bm{\Sigma }\in \mathbb{R}^{d \times d}$. Then, the Wishart matrix with scale matrix $\bm{\Sigma }$ and dimension $(d,n)$, denoted as $W _{d}(\bm{\Sigma },n)$, is defined as 
    \begin{equation*}
        W _{d}(\bm{\Sigma },n)=\mathbf{X} ^{\top }\mathbf{X}.
    \end{equation*}
    Specifically, when $\bm{\Sigma }=\mathbf{I}_{d}$ is the identity matrix, the notation is simplified to $W _{d}(n)$. The distribution of $W _{d}(\bm{\Sigma },n)$ is called Wishart distribution with parameter $d,n,\bm{\Sigma }$. When the meaning is clear from the context, we also use notation $W _{d}(\bm{\Sigma },n)$ (and also $W _{d}(n)$) as Wishart distribution.
\end{definition}

For the Gaussian random matrix $X$, the following lemma can be used to bound its maximum and minimum singular values, which can be further used to establish the bound for the $\ell _{2}$ operator norm of $W _{d}(\bm{\Sigma },n)$. One can find a proof of this lemma in \cite{wainwright2019high} by the concentration bound for Lipschitz functions and the Sudakov-Fernique comparison inequality.

\begin{lemma}[]\label{lemma: bounds of singular values of X}
    Assume $X \in \mathbb{R}^{n \times d}$, and each row is independently drawn from $\mathcal{N}(0,\bm{\Sigma })$, then for any $\epsilon >0$, we have
    \begin{equation*}
        \mathbb{P}\left(\frac{\sigma _{\text{max}}(X)}{\sqrt{n}}\ge \gamma _{\text{max}}(\sqrt{\bm{\Sigma }})(1+\epsilon )+\sqrt{\frac{\text{tr}(\bm{\Sigma })}{n}}\right)\le \exp\left\{-\frac{n\epsilon ^{2}}{2}\right\}.
    \end{equation*}
    When $n \ge d$, we further have the following control on minimum singular value of $X$:
    \begin{equation*}
        \mathbb{P}\left(\frac{\sigma _{\text{min}}(X)}{\sqrt{n}}\le \gamma _{\text{min}}(\sqrt{\bm{\Sigma }})(1-\epsilon  )-\sqrt{\frac{\text{tr}(\bm{\Sigma })}{n}}\right)\le \exp\left\{-\frac{n\epsilon ^{2}}{2}\right\}.
    \end{equation*}
    Here $\sigma _{\text{max}}(X)$ (or $\sigma _{\text{min}}(X)$) is the max (or min) singular value of $X$ and $\gamma _{\text{max}}(\sqrt{\bm{\Sigma }})$ (or $\gamma _{\text{min}}(\sqrt{\bm{\Sigma }})$) is the max (or min) eigenvalue of $\sqrt{\bm{\Sigma }}$. 
\end{lemma}

We can leverage Lemma \ref{lemma: bounds of singular values of X} to bound the $\ell _{2}$ operator norm of $W _{d}(n)$.

\begin{corollary}[]\label{corollary: bound on the l2-operator norm on Wishart matrix}
    Assume $W=\mathbf{X} ^{\top }\mathbf{X} \in \mathbb{R}^{d \times d}$ is drawn from $W _{d}(n)$, and $n \ge d$, then for any $\delta >0$, with probability higher than $1-2\delta $, we have
    \begin{equation*}
        \left\lVert W-n \mathbf{I} _{d}\right\rVert_{2}^{}\le C \left(\sqrt{nd}+\sqrt{n \ln \left(\frac{1}{\delta }\right)}+\ln \left(\frac{1}{\delta }\right)\right)
    \end{equation*}
\end{corollary}

\begin{proof}
    Set $\bm{\Sigma }=\mathbf{I}_{d}$ in the first inequality of Lemma \ref{lemma: bounds of singular values of X}, we have:
    \begin{equation*}
        \mathbb{P}\left(\frac{\sigma _{\text{max}}(X)}{\sqrt{n}}\ge (1+\epsilon  )+\sqrt{\frac{d}{n}}\right)\le \exp\left\{-\frac{n\epsilon  ^{2}}{2}\right\}.
    \end{equation*}
    Using $\sigma _{\text{max}}(X)=\sqrt{\gamma _{\text{max}}(W _{d}(n))}$ and after some calculation, we have 
    \begin{equation*}
        \mathbb{P}\left(\gamma _{\text{max}}(W _{d}(n))-n \ge d+2 \sqrt{nd}+n\epsilon ^{2}+2n\epsilon +2 \sqrt{nd}\epsilon \right)\le \exp\left\{-\frac{n\epsilon ^{2}}{2}\right\}.
    \end{equation*}
    Setting $\delta =\sqrt{\frac{\ln \left(\frac{1}{\epsilon }\right)}{n}}$ yields 
    \begin{equation*}
        \mathbb{P}\left(\gamma _{\text{max}}(W _{d}(n))-n \ge d+2\sqrt{nd}+\ln \left(\frac{1}{\delta }\right)+2\sqrt{n \ln \left(\frac{1}{\delta }\right)}+2 \sqrt{d \ln \left(\frac{1}{\delta }\right)}\right)\le \delta .
    \end{equation*}
    The lower bound follows from the condition that $n \ge d$. The upper bound is derived with a similar argument. The proof is completed.
\end{proof}

Under our problem setting, the following corollary is also required. It allows the mean vector $\mu $ of the normal distribution to be nonzero, but at the expense of introducing a dependence on $\left\lVert \mu \right\rVert_{2}^{}$.

\begin{corollary}[]\label{corollary: bounds on the l2-operator norm of X^TX and XX^T}
    Assume $\mathbf{X} \in \mathbb{R}^{n \times d}$ and each row is independently drawn from $\mathcal{N}(\mu ,\mathbf{I}_{d})$ with $\rho :=\left\lVert \mu \right\rVert_{2}^{}$, then for $\delta >0$, with probability higher than $1-2\delta $, we have:\\ 
    \begin{tabular}{rl}
        (i), $n \ge d$: & $\left\lVert \mathbf{X} ^{\top }\mathbf{X}-n \mathbf{I}_{d}\right\rVert_{2}^{}\le C _{1}\left((1+\rho )\left(\sqrt{nd}+\sqrt{n \ln \left(\frac{1}{\delta }\right)}\right)+\ln \left(\frac{1}{\delta }\right)+n \rho ^{2}\right)$, \\ 
        (ii), $n <d$: & $\left\lVert \mathbf{X}\mathbf{X} ^{\top } -d \mathbf{I}_{n}\right\rVert_{2}^{}\le C _{2}\left((1+\rho )\left(\sqrt{nd}+\sqrt{d \ln \left(\frac{1}{\delta }\right)}\right)+\ln \left(\frac{1}{\delta }\right)+n \rho ^{2}\right)$.
    \end{tabular}
\end{corollary}

\begin{proof}
    When $n \ge d$, let $Z _{i}=X _{i}-\mu $, we have 
    \begin{equation*}
        \left\lVert  \mathbf{X} ^{\top }\mathbf{X}-n \mathbf{I}_{d}\right\rVert_{2}^{}=\left\lVert \sum\limits_{i=1}^{n}X _{i}X _{i}^{\top }-n \mathbf{I}_{d}\right\rVert_{2}^{}\le \left\lVert \sum\limits_{i=1}^{n}Z _{i} Z _{i}^{\top } -n \mathbf{I}_{d}\right\rVert_{2}^{}+2 \left\lVert \sum\limits_{i=1}^{n}\mu Z _{i}^{\top } \right\rVert_{2}^{}+n\left\lVert \mu \mu ^{\top }\right\rVert_{2}^{}.
    \end{equation*}
    By Corollary \ref{corollary: bound on the l2-operator norm on Wishart matrix}, we already have $\left\lVert \sum\limits_{i=1}^{n}Z _{i}Z _{i}^{\top } -n \mathbf{I}_{d}\right\rVert_{2}^{}\le C _{1}\left(\sqrt{nd}+\sqrt{n \ln \left(\frac{1}{\delta }\right)}+\ln \left(\frac{1}{\delta }\right)\right)$. We also have $n\left\lVert \mu \mu ^{\top } \right\rVert_{2}^{}=n \alpha ^{2}$. It remains to bound $\left\lVert \sum\limits_{i=1}^{n}\mu Z _{i}^{\top } \right\rVert_{2}^{}$.
    \begin{equation*}
        \left\lVert \sum\limits_{i=1}^{n}\mu Z _{i}^{\top } \right\rVert_{2}^{}\le \left\lVert \mu \right\rVert_{2}^{}\cdot \left\lVert \sum\limits_{i=1}^{n}Z _{i}\right\rVert_{2}^{}.
    \end{equation*}
    Since $Z _{i}\sim \mathcal{N}(0,\mathbf{I}_{d})$ for each $i$, from \ref{lemma: sub_Gaussian_norm_of_l2_norm} we know that $\left\lVert \sum\limits_{i=1}^{n}\frac{Z _{i}}{\sqrt{n}}-d\right\rVert_{2}^{}$ is a sub-Gaussian random variable and therefore with probability higher than $1-\delta $ we have:
    \begin{equation*}
        \left\lVert \sum\limits_{i=1}^{n}\frac{Z _{i}}{\sqrt{n}}-d\right\rVert_{2}^{}\le C \sqrt{\ln \left(\frac{1}{\delta }\right)},
    \end{equation*}
    which leads to $\left\lVert \sum\limits_{i=1}^{n}Y _{i}\right\rVert_{2}^{}\le C \left(\sqrt{nd}+\sqrt{n \ln \left(\frac{1}{\delta }\right)}\right)$. Combined with previous results on $\left\lVert \sum\limits_{i=1}^{n}Z _{i}Z _{i}^{\top } -n \mathbf{I}_{d}\right\rVert_{2}^{}$ and $n\left\lVert \mu \mu ^{\top } \right\rVert_{2}^{}$, we complete the proof of (i).\\ 
    To bound $\left\lVert \mathbf{X}\mathbf{X} ^{\top } -d \mathbf{I}_{n}\right\rVert_{2}^{}$, we notice that $\mathbf{X}\mathbf{X} ^{\top } $ has the same non-trivial eigenvalues with $\mathbf{X} ^{\top }\mathbf{X}=\mathbf{Z} ^{\top } \mathbf{Z}+2 \sum\limits_{i=1}^{n}\mu Z _{i}^{\top } +n\mu \mu ^{\top } $, in this case, we only need to address the $\ell _{2}$ operator norm of $\mathbf{Z} ^{\top } \mathbf{Z}$, which again equals to $\mathbf{Z}\mathbf{Z} ^{\top } $. Since $\mathbf{Z}\mathbf{Z} ^{\top } $ actually follows the distribution of $W _{n}(d)$, the results follow from the proof of (i).
\end{proof}

The next lemma is from Fact 4.2 in \cite{dong2019quantumentropyscoringfast}. It is an analog of Corollary \ref{corollary: bound on the l2-operator norm on Wishart matrix}, but within the field when the $\ell _{1}$ norm of a weight $\omega $ is small compared to $n$.

\begin{lemma}[]\label{lemma: bound on the max eigenvalue of the covariance matrix}
    Assume $\mathbf{X} \in \mathbb{R}^{n \times d}$ and each row is independently drawn from $\mathcal{N}(0,\mathbf{I}_{d})$. Then, with some constant $C _{1}>0$ and another sufficiently large constant $C _{2}$, as long as $n \ge C _{2}\cdot \frac{1}{\epsilon }$, the following holds with probability at least $1-\delta $ for any $v,\omega $ such that $\left\lVert v\right\rVert_{2}^{}=1$, $\left\lVert \omega \right\rVert_{1}^{}\le C _{2}\cdot \epsilon n$:
    \begin{equation*}
        \left\lVert \sum\limits_{i=1}^{n}\omega _{i}X _{i}X _{i}^{\top } \right\rVert_{2}^{}\le C _{1}\left(\epsilon n \ln \left(\frac{1}{\epsilon }\right)+d+\ln \left(\frac{1}{\delta }\right)\right).
    \end{equation*}
\end{lemma}

\begin{corollary}[]\label{corollary: bound on the max eigenvalue of the covariance matrix in general}
    Assuming the same conditions as Lemma \ref{lemma: bound on the max eigenvalue of the covariance matrix}, we further allow $X _{i}$ from $\mathcal{N}(\mu ,\mathbf{I}_{d})$ with $\left\lVert \mu \right\rVert_{2}^{}\le \alpha $. As a consequence, the inequality will track the dependence on $\alpha $:
    \begin{equation*}
        \left\lVert \sum\limits_{i=1}^{n}\omega _{i}X _{i}X _{i}^{\top } \right\rVert_{2}^{}\le C \left[\epsilon n \left(\ln \left(\frac{1}{\epsilon }\right)+\alpha ^{2}\right)+\sqrt{\epsilon n}\alpha + d+\ln \left(\frac{1}{\delta }\right)\right].
    \end{equation*}
\end{corollary}

\begin{proof}
    The proof is basically from Lemma \ref{lemma: bound on the max eigenvalue of the covariance matrix} the standard argument of concentration. For any $v$ such that $\left\lVert v\right\rVert_{2}^{}=1$, we have:
    \begin{equation*}
        v ^{\top } \left(\sum\limits_{i=1}^{n}\omega _{i}X _{i}X _{i}^{\top }\right)v=v ^{\top } \left(\sum\limits_{i=1}^{n}(X _{i}-\mu )(X _{i}-\mu )^{\top } \right)v+\sum\limits_{i=1}^{n}\omega _{i}(v ^{\top } \mu )^{2}+2 (v ^{\top } \mu )\sum\limits_{i=1}^{n}\omega _{i}(v ^{\top }\left( X _{i}-\mu\right))
    \end{equation*}
    Bound of the first term can be obtained from Lemma \ref{lemma: bound on the max eigenvalue of the covariance matrix}. For the second term, by the Cauchy--Schwarz inequality we can bound it by $2\left\lVert \omega \right\rVert_{1}^{} \cdot \alpha ^{2}\le C \cdot \epsilon n \alpha ^{2}$. For the third term, since $\sum\limits_{i=1}^{n}\omega _{i}v ^{\top } (X _{i}-\mu )$ follows from $\mathcal{N}(0,\left\lVert \omega \right\rVert_{2}^{2})$, we know we can bound it by $C \cdot \alpha \left\lVert \omega \right\rVert_{2}^{}\le C \cdot \alpha \sqrt{\epsilon n}$. Combining them together yields the result.
\end{proof}

H\"older's inequality is a fundamentally important tool in probability, mathematical analysis and functional analysis. We briefly review it here.

\begin{lemma}[H\"older's inequality]\label{lemma: holder's inequality}
    Let $(\Omega ,\mathscr{F})$ be a $\sigma $-algebra equipped with a measure $\mu $. For a pair of conjugate numbers $p,q$: $\frac{1}{p}+\frac{1}{q}, 1 \le p,q \le \infty $, and $f \in L _{p}(\Omega ,\mu )$, $g \in L _{q}(\Omega ,\mu )$, we have 
    \begin{equation*}
        \left\lVert fg\right\rVert_{1}^{}=\displaystyle\int_{\Omega }^{}\left|fg\right|d\mu \le \left(\displaystyle\int_{\Omega }^{}\left|f\right|^{p}d\mu \right)^{\frac{1}{p}}\left(\displaystyle\int_{\Omega }^{}\left|g\right|^{q}d\mu \right)^{\frac{1}{q}}=\left\lVert f\right\rVert_{p}^{}\left\lVert g\right\rVert_{q}^{}.
    \end{equation*}
    The equality holds if and only if $a \left|f\right|^{p}=b \left|g\right|^{q}$ a.s. for some constants $a, b$.
\end{lemma}

The next lemma is a famous result for the combination number.

\begin{lemma}[]\label{lemma: combination number}
    For any $0 \le k \le n$, we have $\sum\limits_{i=0}^{k}\binom{n}{i}\le \left(\frac{en}{k}\right)^{k}$.
\end{lemma}

\begin{proof}
    By the Taylor's expansion of $e ^{k}$, we know $e ^{k}=\sum\limits_{i=0}^{\infty }\frac{k ^{i}}{i!}\ge \sum\limits_{i=0}^{k}\frac{k ^{i}}{i!}$. Therefore, it suffices to show that $\frac{k ^{k}}{n ^{k}}\binom{n}{i}\le \frac{k ^{i}}{i!}$ for any $0 \le i \le k$, and this is equivalent to showing that $\frac{k ^{k}}{n ^{k}}\prod\limits_{j=0}^{i-1}(n-j)\le k ^{i}$. However, the last statement can be directly obtained from the fact that $\frac{k ^{k-i}}{n ^{k-i}}\le 1 \le \frac{n ^{i}}{\prod\limits_{j=0}^{i-1}(n-j)}$.
\end{proof}

Finally, we present an intuitive and useful result at the intersection of convexity and combinatorics. While this result was originally introduced in \cite{10353143}, their proof contains a technical gap which we address below.

\begin{lemma}[]\label{lemma: convex combination}
    Denote $[0,1]^{N}$ to be the set $\left\{(x _{1},x _{2},\dots,x _{N}):0 \le x _{i}\le 1,\forall i \in [N]\right\}$. Given any nonnegative integer $k \le N$ and any vector $x \in [0,1]^{N}$ with condition $\left\lVert x \right\rVert_{1}^{}\le k$, there exists a positive integer $M$, sets $T _{1},T _{2},\dots,T _{M}\subset [N]$, and positive numbers $a _{1},a _{2},\dots,a _{M}$ such that the following statements hold:\\
    \begin{tabular}{rl}
        (i), & $\left|T _{i}\right|\le k, \forall i \in [M]$;\\ 
        (ii), & $0 \le a _{i} \le 1, \forall i \in [M],\sum\limits_{i=1}^{M}a _{i}\le 1$;\\ 
        (iii), & $x =\sum\limits_{i=1}^{M}a _{i}\mathbf{1}_{T _{i}}$,
    \end{tabular}\\
    where same as before $\mathbf{1}_{T _{i}}$ is the indicator vector of the set $T _{i}$.
\end{lemma}

\begin{proof}

We first show the following result: for any vector $(x_1,\ldots, x_N) \in [0,1]^{n}$ and an given integer $k \le N$ such that $\left\lVert x\right\rVert_{1}^{} = k$, and $0\le x_i \le 1$ there exist sets $T_i$ of cardinality $|T_i| = k$ and weights $a _{i} \ge 0$ with $\sum\limits_{i=1}^{M}a _{i}=1$ such that $x = \sum\limits_{i=1}^{M}a _{i} \mathbf{1}_{T_i}$.

The proof of this statement uses Birkhoff-von Neumann's theorem. Construct the following $N \times N$ matrix

\begin{equation*}
    W :=\begin{pmatrix}
        \frac{x _{1}}{k} & \frac{x _{2}}{k} & \dots & \frac{x _{N}}{k}\\ 
        \frac{x _{1}}{k} & \frac{x _{2}}{k} & \dots & \frac{x _{N}}{k}\\ 
        \vdots & \vdots & \dots & \vdots \\
        \frac{1-x _{1}}{N-k} & \frac{1-x _{2}}{N-k} & \dots & \frac{1-x _{N}}{N-k}\\ 
        \vdots & \vdots & \dots & \vdots \\ 
        \frac{1-x _{1}}{N-k} & \frac{1-x _{2}}{N-k} & \dots & \frac{1-X _{N}}{N-k}
    \end{pmatrix},
\end{equation*}
where the $W _{ij}$ is defined as
\begin{equation*}
    W _{ij}=\left\{
    \begin{tabular}{ll}
        $\frac{x _{j}}{k}$, & $1 \le i \le k$,\\ 
        $\frac{1-x _{j}}{N-k}$, & $k+1 \le i \le N$.
    \end{tabular}  
    \right.
\end{equation*}

It is easy to check that this is a doubly stochastic matrix. By Birkhoff-von Neumann's theorem, there exists a set of permutation matrices $P _{1},\dots,P _{M}$ such that 
\begin{equation*}
    W=\sum\limits_{i=1}^{M}\omega _{i}P _{i}.
\end{equation*}
Since the sum of the first $k$ rows of $W$ is $x$, now summing all the first $k$ rows of the matrices $P_i$ shows the desired result.

For the general case, suppose that $k-1 <\left\lVert x\right\rVert_{1}^{}<k$. Consider the augmented vector $\tilde x = (x_1, \ldots, x_N, k - \sum\limits_{i=1}^{N} x_i)$, and apply the logic above to that vector. This gives us a representation $\tilde x = \sum\limits_{i=1}^{\tilde{M}} w_i v_i$, where $v_i$ consists of $k$-ones and $(N + 1 - k)$-zeros. Erasing the last entry of the $v_i$'s, we obtain a representation of $x$ with desired vectors and weights. The proof is completed.
\end{proof}

\newpage

\section{Deferred Proofs of the Lemmas and Main Theorems}\label{section: proofs of the main results}

\subsection{Proof of Lemma \ref{lemma: existence of the optimal vector of a QCO set}}\label{subsection: proof of the existence lemma}

\begin{proof}[Proof of Lemma \ref{lemma: existence of the optimal vector of a QCO set}]
    By the conditions and the definition of the Kolmogorov $k$-width, the following inequality holds:
    \begin{equation*}
        \inf_{P \in \mathcal{P}_{k-1}} \sup_{\theta \in K} \left\lVert \theta -P \theta \right\rVert_{2}^{2} > c ^{2}\sigma ^{2}.
    \end{equation*}
    If we only consider the projections aligned with the coordinate axes, there are $n_k := \binom{d}{k}$ such projections (and countably many when $\mathcal{X}\subset \ell _{2}$). Denote the index set of the whole space $\mathcal{X}$ by $I$ --- i.e., $I=[d]$ when $\mathcal{X}=\mathbb{R}^{d}$ and $\mathbb{N}$ when $\mathcal{X}\subset \ell _{2}$ --- and the index sets of such projections by $J$. We have
    \begin{equation*}
        \inf_{P \in \mathcal{P}_{k-1}} \sup_{\theta \in K} =\left\lVert \theta -P \theta \right\rVert_{2}^{2} \le \inf_{S \in J} \sup_{\theta \in K} \sum_{i \in I} \theta^2_i - \sum_{i \in S}\theta_i^2,
    \end{equation*}
    where the minimum over $S$ is taken with respect to all elements in $J$. Since $K$ is quadratically convex, the RHS above can be written as
    \begin{equation*}
        \inf_{S \in J} \sup_{\theta \in K} \sum_{i \in I} \theta^2_i - \sum_{i \in S}\theta_i^2 \le \inf _{S \in J} \sup_{t \in K^2} \sum_{i \in I} t_i - \sum_{i \in S}t_i = \inf_{w \in S_k} \sup_{t \in K^2} \left\lVert t\right\rVert_{1}^{} - \langle w, t\rangle,
    \end{equation*}
    where $w$ ranges in the set $S_k := \left\{e \left\lvert\right. e \in \mathbb{R}^{d} \text{ (or $\mathbb{R}^{\infty }$)}, e _{i}\in \left\{0,1\right\}, \left\lVert e\right\rVert_{1}^{}=k-1\right\}$. It follows that for each small enough $\varepsilon > 0$ and each $w_i \in S_k$ there exists a $t_i$ such that $\left\lVert t _{i}\right\rVert_{1}^{} - \langle w_i, t_i\rangle \ge c ^{2}\sigma ^{2} - \varepsilon$. Since $K$ is convex and orthosymmetric, we may assume without loss of generality that $t_i \ge 0$ for all $i \in I$, and $t _{i}=0$ on the support of $w_i$. In this way, after appropriately scaling down the coordinates of $t _{i}$ (which is permissible since $K$ is convex), we obtain 
    \begin{equation*}
        \left\lVert t _{i}\right\rVert_{1}^{}=\left\lVert t _{i}\right\rVert_{1}^{}-\left\langle w _{i}, t _{i} \right\rangle=c ^{2}\sigma ^{2}-\epsilon ,
    \end{equation*}
    where the second equality follows from the fact that $\mathbf{0} \in K$. We now argue that there exists a convex combination $t_{\alpha} := \sum\limits_{i \in J}^{} \alpha _{i}t _{i}$ such that $\left\lVert t _{\alpha }\right\rVert_{\infty}^{} \le \frac{c ^{2}\sigma ^{2}}{k}$. First, we observe that since all $t_i$ have nonnegative entries, we have $\left\lVert t _{\alpha }\right\rVert_{\infty }^{} = \sup\limits_{r: r \ge 0, \left\lVert r\right\rVert_{1}^{} = 1} \langle r, t_{\alpha}\rangle$. Hence, it suffices to show that
    \begin{equation*}
        \inf_{\alpha: \alpha \ge 0, \left\lVert \alpha \right\rVert_{1}^{} = 1} \sup_{r: r \ge 0, \|r\|_1 = 1} \langle r, t_{\alpha}\rangle \le \frac{c ^{2}\sigma ^{2}}{k}.
    \end{equation*}
    Since both sets over which the optimization is performed are convex, and the function $ \langle r, t_{\alpha}\rangle = \sum\limits_{i \in J} \alpha_i \langle r, t_i \rangle$ is convex-concave (in fact, it is linear in both arguments), by the minimax theorem, we have
    \begin{equation*}
        \inf_{\alpha: \alpha \ge 0, \left\lVert \alpha \right\rVert_{1}^{} = 1} \sup_{r: r \ge 0, \left\lVert r\right\rVert_{1}^{} = 1} \langle r, t_{\alpha}\rangle = \sup_{r: r \ge 0, \left\lVert r\right\rVert_{1}^{} = 1} \inf_{\alpha: \alpha \ge 0, \left\lVert \alpha \right\rVert_{1}^{} = 1} \langle r, t_{\alpha}\rangle = \sup_{r: r \ge 0, \left\lVert r\right\rVert_{1}^{} = 1} \min_{i \in J} \langle r, t_{i}\rangle.
    \end{equation*}
    Let $U$ denote the index set of the largest $k-1$ entries of $r$. Observe that $\min\limits_{i \in J} \langle r, t_{i}\rangle \le \langle r, t_{r}\rangle$, where $t _{r}$ is chosen such that $t _{r,U}=\mathbf{0}$. Consequently,
    \begin{equation*}
        \langle r, t_r\rangle = \sum\limits_{i \not \in U} r_{i} t_{r,i} \le \left\lVert r _{\backslash U}\right\rVert_{\infty }^{}\left\lVert t _{r}\right\rVert_{1}^{} = \left\lVert r _{\backslash U}\right\rVert_{\infty }^{} (c ^{2}\sigma ^{2} - \varepsilon),
    \end{equation*}
    where $r_{\setminus U}$ denotes the vector obtained by removing the largest $k-1$ coordinates of $r$. Finally, observe that $\left\lVert r _{\backslash U}\right\rVert_{\infty }^{} \le \frac{1}{k}$ since $1 = \left\lVert r\right\rVert_{1}^{} \ge k \left\lVert r _{\backslash U}\right\rVert_{\infty }^{}$. Hence, we conclude that there exists $t ^{\star } = t_{\alpha ^{\star }}$ such that
    \begin{equation*}
        \left\lVert t ^{\star }\right\rVert_{\infty }^{} \le \frac{c ^{2}\sigma ^{2}-\epsilon }{k}< \frac{c ^{2}\sigma ^{2}}{k}.
    \end{equation*}
    In addition, since $t^*$ is a convex combination of vectors $t_i$, we must have $\left\lVert t ^{\star }\right\rVert_{1}^{} = c ^{2}\sigma^2 - \varepsilon$. Taking $\varepsilon \rightarrow 0$ and noticing $\left\lVert t ^{\star }\right\rVert_{\infty }^{}=\left\lVert \theta ^{\star }\right\rVert_{\infty }^{2}, \left\lVert t ^{\star }\right\rVert_{1}^{}=\left\lVert \theta ^{\star }\right\rVert_{2}^{2}$ complete the proof.
\end{proof}

\subsection{Proof of Theorem \ref{theorem: lower bound 1}}\label{subsection: proof of the main lower bound 1}

\begin{proof}[Proof of Theorem \ref{theorem: lower bound 1}]
    We first assume $k ^{\star }\ge 1$, and start by a classical treatment of minimax testing lower bounds which we take from \cite{baraud2002non}. Let $\nu$ denote an arbitrary prior on the set $\left\{\mu \left\lvert\right. \mu \in K, \left\lVert \mu \right\rVert_{2}^{} \ge \rho\right\}$, and $\mathbb{P}_{\nu }$ denote the marginal joint distribution of $X _{1},X _{2},\dots,X _{N}$ given the prior $\nu $. We have that
    \begin{equation*}
        \begin{aligned}
            \inf_{\psi: \mathbb{P}_{0}(\psi = 1) \le \alpha} \sup_{\mu \in K, \left\lVert \mu \right\rVert_{2}^{} \ge \rho} \mathbb{P}_{\mu }(\psi = 0) & \ge
            \inf_{\psi: \mathbb{P}_{0}(\psi = 1) \le \alpha} \mathbb{P}_{\nu }(\psi = 0)\\
            & \ge \inf_{\psi: \mathbb{P}_{0}(\psi = 1)\le \alpha }1-\mathbb{P}_{\nu }(\psi = 1)\\
            & \ge \inf_{\psi : \mathbb{P}_{0}(\psi = 1)\le \alpha }1-\mathbb{P}_{0}(\psi =1)+\left(\mathbb{P}_{0}(\psi =1)-\mathbb{P} _{\nu }(\psi =1)\right) \\
            & \ge 1-\alpha -\text{TV}\left(\mathbb{P}_{0},\mathbb{P}_{\nu }\right)\\
            & \overset{\text{(i)}}{=} 1 - \alpha - \frac{1}{2}\left\lVert \mathbb{P}_{0}-\mathbb{P}_{\nu }\right\rVert_{1}^{},
        \end{aligned}
    \end{equation*}
    where (i) uses \eqref{eq: property of tv distance}. Now assuming $\mathbb{P}_\nu \ll \mathbb{P}_0$, and let $L_\nu = \frac{d \mathbb{P}_\nu}{ d \mathbb{P}_0}$ be the Radon--Nikodym derivative of $\mathbb{P}_{\nu }$ to $\mathbb{P}_{0}$. By Cauchy--Schwartz inequality, we have
    \begin{equation*}
        \left\lVert \mathbb{P}_{0}-\mathbb{P}_{\nu }\right\rVert_{1}^{} = \mathbb{E}_{0}\left|L _{\nu }-1\right| \le \sqrt{\mathbb{E}_{0}\left(L _{\nu }-1\right)^{2}}=\sqrt{\mathbb{E}_{0}\left[L _{\nu }^{2}\right]-1}.
    \end{equation*}
    Hence we conclude
    \begin{equation*}
        \inf_{\psi: \mathbb{P}_0(\psi = 1) \le \alpha} \sup_{\mu \in K, \left\lVert \mu \right\rVert_{2}^{} \ge \rho} \mathbb{P}_{\mu }(\psi = 0) \ge 1 - \alpha - \frac{1}{2}\sqrt{\mathbb{E}_{0}\left[L _{\nu }^{2}\right] - 1}.
    \end{equation*}
    Let $\theta$ be the vector from Lemma \ref{lemma: existence of the optimal vector of a QCO set} when $k=k _{1}^{\star }$, which means $\left\lVert \theta \right\rVert_{2}^{}=\frac{\left(k _{1}^{\star }\right)^{1/4}}{\sqrt{N}}\sigma $ and $\left\lVert \theta \right\rVert_{\infty }^{}\le \frac{1}{\sqrt{N}\left(k _{1}^{\star }\right)^{1/4}}\sigma $. We will use the prior $\left(\kappa \gamma_{i} \theta _{i}\right)_{i \in I} \in K$, where $\gamma_i$ are i.i.d. Rademacher random variables for $i \in I$\footnote{$I=[d]$ when $K \subset \mathbb{R}^{d}$ and $I=\mathbb{N}$ when $K \subset \ell _{2}$.} and $0 < \kappa < 1$ is some small constant. Observe that all these vectors have squared Euclidean norm equal to $\rho^2 = \frac{\kappa ^{2}\sqrt{k _{1}^{\star }}}{N} \sigma^2$.

    By the definition of $\mathbb{P}_{\nu }$, calculating the $\chi^2$-divergence of this mixture yields:
    \begin{equation*}
        \begin{aligned}
        \mathbb{E}_{0} (L_\nu)^2 - 1 & = \mathbb{E}_{X _{1},X _{2},\dots,X _{N} \sim \text{i.i.d.} \mathcal{N}(0,\sigma ^{2}\mathbf{I}_{d})}\left(\frac{d \mathbb{P}_{\nu }}{d \mathbb{P}_{0}}\right)^{2}-1\\
        & =\mathbb{E}_{\mu ,\mu ^{\prime}}\mathbb{E}_{X _{1},X _{2},\dots,X _{N} \sim \text{i.i.d.} \mathcal{N}(0,\sigma ^{2}\mathbf{I}_{d})}\frac{\prod\limits_{i=1}^{N}\exp \left\{\frac{-\left\lVert X _{i}-\mu \right\rVert_{2}^{2}}{2 \sigma ^{2}}\right\}\prod\limits_{i=1}^{N}\exp \left\{-\frac{\left\lVert X _{i}-\mu ^{\prime }\right\rVert_{2}^{2}}{2\sigma^2}\right\}}{\prod\limits_{i=1}^{N}\exp \left\{-\frac{\left\lVert X _{i}\right\rVert_{2}^{2}}{\sigma^2}\right\}} - 1 \\
        & = \mathbb{E}_{\mu ,\mu ^{\prime}}\exp\left\{\frac{N\mu ^{\top } \mu ^{\prime }}{\sigma ^{2}}\right\} - 1,
        \end{aligned}
    \end{equation*}
    where by $\mathbb{E}_{\mu ,\mu ^{\prime}}$ we mean expectation over two i.i.d. draws from the prior distribution $\nu $ that we specified above. Using the symmetry of Rademacher random variable, we have 
    \begin{equation*}
        \begin{aligned}
            \mathbb{E}_{\mu ,\mu ^{\prime}}\exp \left\{\frac{N\mu ^{\top } \mu ^{\prime }}{\sigma ^{2}}\right\} - 1 & =\mathbb{E}_{\mu ,\mu ^{\prime}}\exp \left\{N\sum\limits_{i \in I}^{}\frac{\mu _{i}\mu _{i}^{\prime}}{\sigma ^{2}}\right\}-1 \\
            & =\prod_{i \in I}\cosh\left(\frac{N\kappa^2\theta_i^2}{\sigma ^2}\right) - 1\\
            & \le \exp \left\{\sum_{i \in I}\frac{N ^{2}\theta_i^4 \kappa^4}{2 \sigma ^{4}}\right\} - 1\\
            & \le \exp \left\{\frac{N ^{2}\left\lVert \theta \right\rVert_{\infty }^{2}\left\lVert \theta \right\rVert_{2}^{2}\kappa^4}{2\sigma^4}\right\}-1\\
            & \le \exp \left\{\frac{\kappa ^{4}}{2} \right\} - 1,
        \end{aligned}
    \end{equation*}
    where we used the well known inequality $\cosh(x) \le \exp \left\{\frac{x ^{2}}{2}\right\}$, which can be checked simply by comparing the Taylor series, and the properties of $\theta$ from Lemma \ref{lemma: existence of the optimal vector of a QCO set}. It follows that the minimax risk for testing is lower bounded by
    \begin{equation*}
        1 - \alpha - \frac{1}{2}\sqrt{\exp \left\{\frac{\kappa ^{4}}{2}\right\}-1}.
    \end{equation*}
    Recall that the $\ell _{2}$ norm of the selected $\mu $ is $\rho =\left\lVert \mu \right\rVert_{2}^{}=\frac{(k _{1}^{\star })^{1/4}\kappa }{\sqrt{N}}\sigma $. If we let the expression above equal to $\alpha $ and solve for $\kappa $, we obtain $\kappa =c(\alpha )$. Consequently, we conclude that if $\kappa \le c(\alpha )$, and therefore $\rho \le c(\alpha )\cdot \frac{(k _{1}^{\star })^{1/4}}{\sqrt{N}}\sigma $, testing with uniform power $\alpha $ is not possible.

    The remaining corner case is when $k _{1}^{\star }=0$, or equivalently $D _{0}(K)=\sup\limits_{\theta \in K}\left\lVert \theta \right\rVert_{2}^{}\le \frac{1}{\sqrt{N}}\sigma $. In this case, any point $\theta \in K$ satisfies $\left\lVert \theta \right\rVert_{2}^{} \le \frac{1}{\sqrt{N}}\sigma $. Taking a point $\theta _{0}$ that (nearly) achieves the equality, and considering the mixture formed by $\gamma \theta _{0}$, where $\gamma $ is a Rademacher random variable, one obtains that the resulting $\chi^2$-divergence is bounded by $\cosh \left(\frac{N\left\lVert \theta \right\rVert_{2}^{2}}{\sigma ^{2}}\right) - 1 \le \exp \left\{\frac{N ^{2}\left\lVert \theta \right\rVert_{2}^{4}}{2 \sigma ^{4}}\right\}- 1 \le \sqrt{e}-1$ and consequently
    \begin{equation*}
        \inf_{\psi: \mathbb{P}_0(\psi = 1) \le \alpha} \sup_{\mu \in K} \mathbb{P}_{\mu }(\psi = 0) \ge 1 - \alpha - \frac{1}{2}\sqrt{\sqrt{e}-1},
    \end{equation*}
    implying that it is impossible to test the problem. Hence, $k _{1}^{\star }$ must be at least positive for sufficiently small $\alpha $.
\end{proof}

\subsection{Proof of Theorem \ref{theorem: lower bound 2}}\label{subsection: proof of the main lower bound 2}

In this section, we formally prove Theorem \ref{theorem: lower bound 2}. We first recall the \hyperref[lemma: main lemma in the proof of lower bound 2]{core lemma}. To establish Lemma \ref{lemma: main lemma in the proof of lower bound 2}, two auxiliary lemmas, Lemma \ref{lemma: requirement of robust diff private algorithm} and Lemma \ref{lemma: construction of a diff private algorithm}, are presented. Based on these tools, we then construct a symmetric set served as the prior for the mean $\mu $ via the QCO properties and compute the Hamming distance as an application of Lemma \ref{lemma: total variation and coupling} via a sequential optimal coupling. Finally, the lower bound is yielded by the \hyperref[lemma: Pinsker's inequality]{Pinsker's inequality} and a concrete calculation of the \hyperref[def_KL_divergence]{KL divergence}.

\begin{lemma}[Core lemma]\label{lemma: main lemma in the proof of lower bound 2}
    Suppose $\mathbf{X}=(X _{1},X _{2},\dots,X _{N})^{\top }$ and $\mathbf{X} ^{\prime }=(X ^{\prime }_{1},X _{2}^{\prime },\dots,X _{N}^{\prime })^{\top }$ are random matrices with rows $X _{i},X ^{\prime }_{i}\in \mathbb{R}^{d}$. Assume the laws of $\mathbf{X}$ and $\mathbf{X} ^{\prime }$ are $\mathcal{U}$ and $\mathcal{V}$, respectively. Additionally, assume that there exists a coupling between $\mathcal{U}$ and $\mathcal{V}$ under which the expectation of the Hamming distance between $\mathbf{X}$ and $\mathbf{X} ^{\prime }$ is $\mathcal{O}(\epsilon N)$, i.e.,
\begin{equation}
    \mathbb{E}_{(\mathbf{X},\mathbf{X} ^{\prime })\sim (\mathcal{U},\mathcal{V})}d _{H}(\mathbf{X},\mathbf{X} ^{\prime })\lesssim \epsilon N.
\end{equation}
Then, there exists no robust test capable of distinguishing between $\mathcal{U}$ and $\mathcal{V}$ under contamination level of $\varOmega (\epsilon )$ while achieving type-I and type-II error probabilities below $0.1$ simultaneously and uniformly.
\end{lemma}

We start with a review of the notion of $(\epsilon ,\delta )$-differential privacy for convenience.

\begin{definition}[$(\epsilon ,\delta )$-differential privacy, \cite{10.1007/11681878_14}]\label{def_diff_privacy}
    For a pair of non-negative real numbers $(\epsilon ,\delta )$, a randomized algorithm $\mathcal{A}$ defined on the sample space $\Omega $ is an $(\epsilon ,\delta )$-differentially private algorithm if for any $S \subset \text{range}(\mathcal{A})$ and any $X,X ^{\prime }\in \Omega $ with the Hamming distance $d _{H}(X, X ^{\prime })\le 1$, we have
\begin{equation}
    \mathbb{P}\left(\mathcal{A}(X)\in S\right)\le e ^{\epsilon }\mathbb{P}\left(\mathcal{A}(X ^{\prime })\in S\right)+\delta .
\end{equation}
When $\delta = 0$, the notion is called \textit{pure differential privacy}; when $\delta > 0$, it is referred to as \textit{approximate differential privacy}.
\end{definition}

Lemma \ref{lemma: requirement of robust diff private algorithm} is motivated by Theorem 11 in \cite{NEURIPS2018_7de32147}, and shares a nearly identical conclusion with the \hyperref[lemma: main lemma in the proof of lower bound 2]{core lemma}.

\begin{lemma}[]\label{lemma: requirement of robust diff private algorithm}
    Assume that there is a coupling between $\mathcal{U}$ and $\mathcal{V}$ under which that the expectation of the Hamming distance between $\mathbf{X}$ and $\mathbf{X} ^{\prime }$ is smaller or equal than $D$, i.e.,
\begin{equation}
    \mathbb{E}_{(\mathbf{X},\mathbf{X} ^{\prime })\sim (\mathcal{U},\mathcal{V})}d _{H}(\mathbf{X},\mathbf{X} ^{\prime })\le D.
\end{equation}
Then, if a test $\mathcal{A}$ to distinguish between $\mathcal{U}$ and $\mathcal{V}$: $\Omega \mapsto \left\{0,1\right\}$ is a $(0, \delta )$-differentially private algorithm, with type-I and type-II error probabilities both below $0.1$, then we must have $D=\varOmega (\frac{1}{\delta })$.
\end{lemma}

\begin{proof}
    By the assumption of the lemma, suppose that $(\mathbf{X},\mathbf{X} ^{\prime })$ is the random variable from the coupling. We have $\mathbb{P}\left(\mathcal{A}(\mathbf{X})=0\right)\ge 0.9$ and $\mathbb{P}\left(\mathcal{A}(\mathbf{X} ^{\prime })=1\right)\ge 0.9$. By Markov's inequality, we have $\mathbb{P}\left(d _{H}(\mathbf{X},\mathbf{X} ^{\prime })>10D\right)\le \frac{\mathbb{E}d _{H}(\mathbf{X},\mathbf{X} ^{\prime })}{10D}\le 0.1$. From the union bound, we obtain that
\begin{equation}
    \mathbb{P}\left(\mathcal{A}(\mathbf{X})=0 \cap \mathcal{A}(\mathbf{X} ^{\prime })=1 \cap d _{H}(\mathbf{X},\mathbf{X} ^{\prime })<10D\right)\ge 0.7.
    \label{eq:basic_probability_inequality}
\end{equation}
By the definition of \hyperref[def_diff_privacy]{differential privacy}, we can construct a chain from $\mathbf{X}$ to $\mathbf{X} ^{\prime }$ such that each neighbour points differ on at most one coordinate, by which we have
\begin{equation}     |\mathbb{P}\left(\mathcal{A}(\mathbf{X})=0\right)-\mathbb{P}\left(\mathcal{A}(\mathbf{X} ^{\prime })=0\right)|\le \delta \cdot d _{H}(\mathbf{X},\mathbf{X} ^{\prime }).
    \label{eq:difference_1}
\end{equation}
It is worth noting that the probability in \eqref{eq:basic_probability_inequality} includes two parts of randomness: the randomness of $\mathbf{X}$ and $\mathbf{X} ^{\prime }$ from $\Omega $, and the randomness of $\mathcal{A}$. Now take a realization of $\mathbf{X}$ and $\mathbf{X} ^{\prime }$, denoted by $\mathbf{x}$ and $\mathbf{x} ^{\prime }$, that satisfy \eqref{eq:basic_probability_inequality}, we have
\begin{equation}
    1-0.7 \ge \mathbb{P}(\mathcal{A}(\mathbf{x} ^{\prime })=0)\ge \mathbb{P}(\mathcal{A}(\mathbf{x})=0)-\delta \cdot d _{H}(\mathbf{x},\mathbf{x} ^{\prime })\ge 0.7-\delta \cdot 10D,
\end{equation}
from which we obtain that $10D \ge 0.4 \frac{1}{\delta }$. The proof is completed.
\end{proof}

Lemma \ref{lemma: construction of a diff private algorithm} shows how to obtain a $(0,\delta )$-differentially private test from a given robust test between $\mathcal{U}$ and $\mathcal{V}$ using the black-box robustness-to-privacy transformation. The techniques are from \cite{10353143} and \cite{10.1145/3564246.3585115}.

\begin{lemma}[]\label{lemma: construction of a diff private algorithm}
    Suppose $\mathcal{A}$ is a robust test for the mean testing problem \eqref{eq: testing problem} achieving type-I and type-II errors at most $0.1$ under a contamination level at most $\epsilon $. Then, for $\delta =\frac{1}{\epsilon N}$, there exists a $(0,\delta )$-differentially private algorithm $\widetilde{\mathcal{A}}$ whose error probabilities are also at most $0.1$ for any dataset with size $N$.
\end{lemma}

\begin{proof}
    For any algorithm $\mathcal{A}$ that attempts to distinguish between $H _{0}$ and $H _{1}$ based on $\mathbf{X}$, define $\mathcal{S}(X,\mathcal{A})$ as the smallest number such that there exists another dataset $\mathbf{X} ^{\prime }$ with $\mathcal{A}(\mathbf{X} ^{\prime })=1$ and $d _{H}(\mathbf{X},\mathbf{X} ^{\prime })=\mathcal{S}(\mathbf{X},\mathcal{A})$. Using this quantity, we construct a new algorithm $\widetilde{\mathcal{A}}$ that outputs $1$ with probability $\max \left(0,1-\delta \cdot \mathcal{S}(X,\mathcal{A})\right)$. We now prove that $\widetilde{\mathcal{A}}$ is the desired algorithm.

First, we notice that for any two datasets $\mathbf{X} ^{\prime }$, $\mathbf{X} ^{\prime \prime }$ that satisfy $d _{H}(\mathbf{X} ^{\prime },\mathbf{X} ^{\prime \prime })\le 1$, we always have $|\mathcal{S}(X^{\prime },\mathcal{A})-\mathcal{S}(X ^{\prime \prime },\mathcal{A})|\le 1$ for any algorithm $\mathcal{A}$. This is foreseen since if we denote $\mathbf{X}^{\dag}$ as the specific dataset that satisfies $\mathcal{A}(\mathbf{X}^{\dag})=1$ and $d _{H}(\mathbf{X}^{\prime },\mathbf{X}^{\dag})=\mathcal{S}(\mathbf{X}^{\prime },\mathcal{A})$, then by the triangular inequality we have $d _{H}(\mathbf{X} ^{\prime \prime },\mathbf{X}^{\dag})\le \mathcal{S}(\mathbf{X}^{\prime },\mathcal{A})+1$, and consequently $\mathcal{S}(\mathbf{X} ^{\prime \prime },\mathcal{A})\le \mathcal{S}(\mathbf{X} ^{\prime },\mathcal{A})+1$. In the same way we also have $\mathcal{S}(\mathbf{X} ^{\prime \prime },\mathcal{A})\ge \mathcal{S}(\mathbf{X} ^{\prime },\mathcal{A})-1$. The statement is proved.

By the statement above we know that the constructed algorithm $\widetilde{\mathcal{A}}$ is a $(0,\delta )$-differentially private algorithm since for any two dataset $\mathbf{X}$ and $\mathbf{X} ^{\prime }$ with $d _{H}(\mathbf{X}, \mathbf{X} ^{\prime })\le 1$, we have
\begin{equation*}
    \begin{aligned}
        |\mathbb{P}\left(\widetilde{\mathcal{A}}(\mathbf{X})=1\right)-\mathbb{P}\left(\widetilde{\mathcal{A}}(\mathbf{X} ^{\prime })=1\right)| & =|\max \left\{0,1-\delta \cdot \mathcal{S}(\mathbf{X},\mathcal{A})\right\}-\max \left\{0,1-\delta \cdot \mathcal{S}(\mathbf{X} ^{\prime },\mathcal{A})\right\}|\\
        & \le |\delta \cdot \mathcal{S}(\mathbf{X},\mathcal{A})-\delta \cdot \mathcal{S}(\mathbf{X} ^{\prime },\mathcal{A})|\\
        & \le \delta .
    \end{aligned}
\end{equation*}
Now assume $\tilde{\mathbf{X}}$ is from $H _{0}$. By the definition of robust testing, with probability at least $0.9$, we can arbitrarily modify $\tilde{\mathbf{X}}$ and obtain the same result as long as the number of modifications does not exceed $\epsilon N$. This equivalently says that $\mathcal{S}(\tilde{\mathbf{X}},\mathcal{A})\ge \epsilon N$, which, by the construction of $\widetilde{\mathcal{A}}$, implies that $\widetilde{\mathcal{A}}$ will output $1$ with probability $0$ in this case. On the other hand, if $\tilde{\mathbf{X}}$ is from $H _{1}$, then with probability at least $0.9$, we will have $\mathcal{A}(X)=1$, in which case $\mathcal{S}(X,\mathcal{A})=0$ and consequently $\widetilde{\mathcal{A}}$ will output $1$ with probability $1$. This shows that $\widetilde{\mathcal{A}}$ is a $(0,\delta )$-differential privacy with error probabilities at most $0.1$ under both hypotheses.
\end{proof}

We are now prepared to prove the \hyperref[lemma: main lemma in the proof of lower bound 2]{core lemma} and Theorem \ref{theorem: lower bound 2}. 

\begin{proof}[Proof of Lemma \ref{lemma: main lemma in the proof of lower bound 2}]
    Suppose we already have a robust test $\mathcal{A}$ for the mean testing problem \eqref{eq: testing problem} with the fraction of contamination at most $\epsilon $. Then by Lemma \ref{lemma: construction of a diff private algorithm}, we can contruct a $(0,\delta )$-differential privacy $\widetilde{\mathcal{A}}$ with errors at most $0.1$. Therefore, the constructed algorithm $\widetilde{\mathcal{A}}$ satifies the condition of Lemma \ref{lemma: requirement of robust diff private algorithm}. Take $D=\mathbb{E}d _{H}(\mathbf{X},\mathbf{X} ^{\prime })$ in Lemma \ref{lemma: requirement of robust diff private algorithm}, and by the conclusion of Lemma \ref{lemma: requirement of robust diff private algorithm} we should have $\mathbb{E}d _{H}(\mathbf{X},\mathbf{X} ^{\prime })\gtrsim \frac{1}{1/\epsilon N}=\epsilon N$. This necessary condition finishes the proof.
\end{proof}

\begin{proof}[Proof of Theorem \ref{theorem: lower bound 2}]

Take $\mathcal{U}$ to be the law of $N$ i.i.d. samples from $\mathcal{N}\left(0,\sigma ^{2}\mathbf{I}\right)$ and $\mathcal{V}$ to be the law of $N$ i.i.d. samples from $\mathcal{N}\left(v,\sigma ^{2}\mathbf{I}\right)$ with $v$ drawn from the uniform distribution on the set
\begin{equation*}
    V :=\left\{v \left\lvert\right. v _{i}=\gamma _{i}\theta _{i}, \gamma _{i}\in \left\{-1,1\right\}, i \ge 1,i \in I\right\},
\end{equation*}
where $\theta $ is the vector from Lemma \ref{lemma: existence of the optimal vector of a QCO set} setting $c=\frac{(k _{2}^{\star })^{1/4}\sqrt{\epsilon }}{N ^{1/4}}\sigma $, i.e., $\left\lVert \theta \right\rVert_{2}^{}= \frac{(k _{2}^{\star })^{1/4}\sqrt{\epsilon }}{N ^{1/4}}\sigma $ and $\left\lVert \theta \right\rVert_{\infty }^{}\le \frac{\sqrt{\epsilon }}{(k _{2}^{\star })^{1/4}N ^{1/4}}\sigma $. We are going to compute $d _{H}(\mathbf{X},\mathbf{X} ^{\prime })$ using Lemma \ref{lemma: total variation and coupling}. Specifically, we compare (in terms of total variation) a distribution of $k$ exchangeable variables generated from $\mathcal{V}$ to a joint distribution of $k-1$ exchangeable variables generated from $\mathcal{N}\left(v ,\sigma ^{2}\mathbf{I}\right)$ and the last one from $\mathcal{N}\left(0,\sigma ^{2}\mathbf{I}\right)$ for $1 \le k \le N$. We approach this by decomposing both $\mathbf{X}$ and $\mathbf{X} ^{\prime }$ by their sequence of conditional probabilities and constructing the optimal coupling between $X _{k}$ and $X ^{\prime }_{k}$ conditional on $X _{1},\dots,X _{k-1}$ using Lemma \ref{lemma: total variation and coupling} such that
\begin{equation*}
    \mathbb{P}\left(X _{k}\neq X ^{\prime }_{k}|X _{1},\dots,X _{k-1}\right)=\text{TV}(p(x _{k}|x _{1},\dots,x _{k-1}),q(x _{k})),
\end{equation*}
where $p (x _{k}|x _{1},\dots,x _{k-1})$ is the density of conditional distribution of the first $k$ coordinates of $X ^{\prime }$ and $q(x _{k})$ is the density of standard normal distribution. If we denote the optimal coupling above by $\mathcal{P}_{k}=\mathcal{P}_{k}(X _{k},X ^{\prime }_{k}|X _{1},...,X _{k-1})$, it is not hard to verify that
\begin{equation*}
    \tilde{\mathcal{P}}(\mathbf{X},\mathbf{X}^{\prime }):=\mathcal{P}_{1}(X _{1},X ^{\prime }_{1})\cdot \prod\limits_{i=2}^{N}\mathcal{P}_{i}(X _{k},X ^{\prime }_{k}|X _{1},\dots,X _{k-1})
\end{equation*}
is a valid coupling between $\mathcal{U}$ and $\mathcal{V}$. By the property of the optimal coupling and the definition of total variation distance, we have under $\tilde{\mathcal{P}}$ that
\begin{equation*}
    \begin{aligned}
        \mathbb{E}_{\mathbf{X},\mathbf{X} ^{\prime }\sim \tilde{\mathcal{P}}}\mathbf{1}_{\left\{X _{k}\neq X ^{\prime }_{k}\right\}} & =\mathbb{E}\text{TV}(p(x _{k}|x _{1},\dots,x _{k-1}),q(x _{k})) \\
                                                                                                                    & =\text{TV}(p(x _{1},\dots,x _{k}),q(x _{1},\dots,x _{k})),
    \end{aligned}
\end{equation*}
where $p(x _{1},\dots,x _{k})$ is the joint distribution of the first $k$ coordinates of $X ^{\prime }$ and $q(x _{1},\dots,x _{k}):=p(x _{1},\dots,x _{k-1})\cdot q(x _{k})$. Since 
\begin{equation*}
    \mathbb{E}_{(\mathbf{X},\mathbf{X} ^{\prime })\sim \tilde{\mathcal{P}}}d _{H}(\mathbf{X},\mathbf{X} ^{\prime })=\mathbb{E}_{\mathbf{X},\mathbf{X} ^{\prime }\sim \tilde{\mathcal{P}}}\sum\limits_{i=1}^{N}\mathbf{1}_{\left\{X _{i}\neq X ^{\prime }_{i}\right\}},
\end{equation*}
to evaluate $\mathbb{E}_{(\mathbf{X},\mathbb{X} ^{\prime })\sim \tilde{\mathcal{P}}}d _{H}(\mathbf{X},\mathbf{X} ^{\prime })$, it suffices to evaluate 
\begin{equation*}
    \sum\limits_{k=1}^{N}\text{TV}(p(x _{1},\dots,x _{k}), q(x _{1},\dots,x _{k})).
\end{equation*}
We bound each individual term by the \hyperref[lemma: Pinsker's inequality]{Pinsker's inequality}, where we have
\begin{equation*}
    \begin{aligned}
        2\text{TV}^{2}(p(x _{1},\dots,x _{k}), q(x _{1},\dots,x _{k})) & \le \text{KL}(p(x _{1},\dots,x _{k})|q(x _{1},\dots,x _{k}))\\
        &\overset{\text{(i)}}{=}\mathbb{E}\ln \left(\frac{\frac{1}{2 ^{|V|}}\sum\limits_{v \in V}^{}\exp\left\{\left(\sum\limits_{i=1}^{k}X _{i}\right)^{\top } v/\sigma ^{2}\right\}}{\frac{1}{2 ^{|V|}}\sum\limits_{v \in V}^{}\exp\left\{\left(\sum\limits_{i=1}^{k-1}X _{i}\right)^{\top } v/\sigma ^{2}\right\}}\right)-\frac{\left\lVert v \right\rVert_{2}^{2}}{2 \sigma ^{2}}\\
        &\overset{\text{(ii)}}{=}\mathbb{E}\ln \left(\frac{\prod\limits_{j=1}^{d}\cosh \left\{\left(\sum\limits_{i=1}^{k}X _{i}\right)_{j}v _{j}/\sigma ^{2}\right\}}{\prod\limits_{j=1}^{d}\cosh \left\{\left(\sum\limits_{i=1}^{k-1}X _{i}\right)_{j}v _{j}/\sigma ^{2}\right\}}\right)-\frac{\left\lVert v \right\rVert_{2}^{2}}{2 \sigma ^{2}},
    \end{aligned}
\end{equation*}
where $(\cdot )_{j}$ denotes the $j$-th coordinate, (i) is from that $\left\lVert v\right\rVert_{2}^{2}$ does not depend on the specific choice of $v \in V$ by the construction of $V$ and (ii) is from an important observation that
\begin{equation*}
    \frac{1}{2 ^{|V|}}\sum\limits_{v \in V}^{}\exp\left\{\left(\sum\limits_{i=1}^{k}X _{i}\right)^{\top } v\right\}=\prod\limits_{j=1}^{d}\cosh \left\{\left(\sum\limits_{i=1}^{k}X _{i}\right)_{j}v _{j}\right\}.
\end{equation*}
Therefore, we have
\begin{equation*}
    \begin{aligned}
        \sum\limits_{k=1}^{N}\text{KL}(p(x _{1},\dots,x _{k})|q(x _{1},\dots,x _{k}))&=\sum\limits_{k=1}^{N}\mathbb{E}\left[\ln \left(\frac{\prod\limits_{j=1}^{d}\cosh \left\{\left(\sum\limits_{i=1}^{k}X _{i}\right)_{j}v _{j}/\sigma ^{2}\right\}}{\prod\limits_{j=1}^{d}\cosh \left\{\left(\sum\limits_{i=1}^{k-1}X _{i}\right)_{j}v _{j}/\sigma ^{2}\right\}}\right)-\frac{\left\lVert v \right\rVert_{2}^{2}}{2 \sigma ^{2}}\right]\\  
        &\overset{\text{(iii)}}{=}\sum\limits_{j=1}^{d}\mathbb{E}\ln \left(\cosh \left\{\left(\sum\limits_{i=1}^{N}X _{i}\right)_{j}v _{j}/\sigma ^{2}\right\}\right)-N \cdot \frac{\left\lVert v\right\rVert_{2}^{2}}{2 \sigma ^{2}},
    \end{aligned}
\end{equation*}
where (iii) is from telescoping and the exchange of summation and expectation. Using the inequality that $\ln \left(\cosh x\right)\le \frac{x ^{2}}{2}, \forall x \in \mathbb{R}$, which can be easily proved by inspecting the Taylor's expansion of $\ln \left(\cosh x\right)$, and the fact that $\left(\sum\limits_{i=1}^{N}X _{i}\right)_{j}\sim \mathcal{N}\left(Nv _{j},N \sigma ^{2}\right)$ for any $j$, we can bound the expression above by 
\begin{equation*}
    \begin{aligned}
        \sum\limits_{k=1}^{N}\text{KL}(p(x _{1},\dots,x _{k})|q(x _{1},\dots,x _{k}))&\le \frac{1}{2}\sum\limits_{j=1}^{d}\mathbb{E}\left[\left(\sum\limits_{i=1}^{N}X _{i}\right)_{j}v _{j}/\sigma ^{2}\right]^{2}-N \cdot \frac{\left\lVert v\right\rVert_{2}^{2}}{2 \sigma ^{2}}\\ 
        &\le \frac{1}{2}\sum\limits_{j=1}^{d}\frac{v _{j}^{2}}{\sigma ^{4}}\left(N ^{2}v _{j}^{2}+N \sigma ^{2}\right)-N \cdot \frac{\left\lVert v\right\rVert_{2}^{2}}{2\sigma ^{2}}\\ 
        &=\frac{N ^{2}\left\lVert v\right\rVert_{4}^{4}}{2\sigma ^{4}}.
    \end{aligned}
\end{equation*}
Combining all of the expressions and inequalities above together, we arrive at the following 
\begin{equation*}
    \begin{aligned}
        \frac{1}{N}\sum\limits_{k=1}^{N}\text{TV}(p(x _{1},\dots,x _{k})|q(x _{1},\dots,x _{k}))&\lesssim \frac{1}{N}\sum\limits_{k=1}^{N}\sqrt{\text{KL}(p(x _{1},\dots,x _{k})|q(x _{1},\dots,x _{k}))}\\ 
        &\overset{\text{(i)}}{\le }\sqrt{\frac{1}{N}\cdot \sum\limits_{k=1}^{N}\text{KL}(p(x _{1},\dots,x _{k})|q(x _{1},\dots,x _{k}))}\\ 
        &\lesssim \frac{\sqrt{N}\left\lVert v\right\rVert_{4}^{2}}{\sigma ^{2}},
    \end{aligned}
\end{equation*}
where (i) follows from the Jensen's inequality. Recall that by Lemma \ref{lemma: existence of the optimal vector of a QCO set} we have $\left\lVert \theta \right\rVert_{2}^{}=\frac{(k _{2}^{\star }) ^{1/4}\sqrt{\epsilon }}{N ^{1/4}}\sigma $ and $\left\lVert \theta \right\rVert_{\infty }^{}\le \frac{\sqrt{\epsilon }}{(k _{2}^{\star }) ^{1/4}N ^{1/4}}\sigma $. Now that $\theta _{i}^{4}\le \theta _{i}^{2}\cdot \left\lVert \theta \right\rVert_{\infty }^{2}$ for all $i \in \left\{1,2,\dots,d\right\}$, we know $\left\lVert \theta \right\rVert_{4}^{2}\le \left\lVert \theta \right\rVert_{2}^{}\left\lVert \theta \right\rVert_{\infty }^{}\le \frac{\epsilon }{\sqrt{N}}\sigma ^{2}$. By the definition of $V$, we know that $\left\lVert v\right\rVert_{4}^{2}=\left\lVert \theta \right\rVert_{4}^{2}$. Therefore, we have $\frac{\sqrt{N}\left\lVert v\right\rVert_{4}^{2}}{\sigma ^{2}}\le \frac{\sqrt{N}}{\sigma ^{2}}\cdot \frac{\epsilon }{\sqrt{N}}\sigma ^{2}\le \epsilon $. By the \hyperref[lemma: main lemma in the proof of lower bound 2]{core lemma}, we have 
\begin{equation*}
    \inf\limits_{\phi \in A _{s}(0.1)}\sup\limits_{\mu \in K}\sup\limits_{\mathcal{C}}\mathbb{P}_{\mu }(\phi (\mathcal{C}(\tilde{\mathbf{X}}))=0)\ge 0.1.
\end{equation*}
The proof is completed.
\end{proof}

\subsection{Proof of Theorem \ref{theorem: lower bound 3}}\label{subsection: proof of the main lower bound 3}

In this section, we prove Theorem \ref{theorem: lower bound 3}. We reach the result by constructing two mixture distributions based on $H _{0}$ and $H _{1}$ which are, in fact, identical. We then evaluate the total variation distance between two truncated binomial distributions, which enables us to construct an adaptive adversary, leading high type-II error. Similar techniques also appear in \cite{prasadan2025informationtheoreticlimitsrobust} and \cite{chen2017robustcovariancescattermatrix}.

\begin{proof}[Proof of Theorem \ref{theorem: lower bound 3}]
We prove the theorem by several steps.\opar

\textbf{Step I: control the TV distance.} Set $\epsilon ^{\prime }=\epsilon -1/\sqrt{2N}$. Since $\epsilon \gtrsim 1/\sqrt{N}$, we have $\epsilon \asymp \epsilon ^{\prime }$. Now fix a $\mu _{0}\in K$ such that $\left\lVert \mu _{0}\right\rVert_{2}^{2}\le \frac{4 \sigma ^{2}\epsilon ^{\prime 2}}{(1-\epsilon ^{\prime })^{2}}$. Consider two Gaussian measures $\mathbb{P}_{0}\sim \mathcal{N}(0,\sigma ^{2}\mathbf{I}_{d})$ and $\mathbb{P}_{\mu _{0}}\sim \mathcal{N}(\mu _{0},\sigma ^{2}\mathbf{I}_{d})$. By the Pinsker's inequality we have $\text{TV}\left(\mathbb{P}_{0},\mathbb{P}_{\mu _{0}}\right)\le \sqrt{\frac{1}{2}D _{\text{KL}}\left(\mathbb{P}_{0}\left\lVert\right. \mathbb{P}_{\mu _{0}}\right)}=\frac{\epsilon ^{\prime }}{1-\epsilon ^{\prime }}$. Since $r(x):=\frac{x}{1-x}$ is an increasing function for $x \in (0,1)$, there exists some $\epsilon ^{\prime \prime }\in (0,\epsilon ^{\prime })$ such that $\text{TV}\left(\mathbb{P}_{0},\mathbb{P}_{\mu _{0}}\right)=\frac{\epsilon ^{\prime \prime }}{1-\epsilon ^{\prime \prime }}$.\opar

\textbf{Step II: construct a pair of identical mixtures.} Let $p _{0}$ and $p _{\mu _{0}}$ be the probability density functions of $\mathbb{P}_{0}$ and $\mathbb{P}_{\mu _{0}}$, construct two distributions $Q _{1},Q _{2}$ with density
\begin{equation}
q _{1}=\frac{(p _{\mu _{0}}-p _{0})\mathbf{1}_{\left\{p _{\mu _{0}}\ge p _{0}\right\}}}{\text{TV}\left(\mathbb{P}_{0},\mathbb{P}_{\mu _{0}}\right)},\ q _{2}=\frac{(p _{0}-p _{\mu _{0}})\mathbf{1}_{\left\{p _{0}>p _{\mu _{0}}\right\}}}{\text{TV}(\mathbb{P}_{0},\mathbb{P}_{\mu _{0}})}.
\end{equation}
Consider the mixture measures 
\begin{equation}
R_{1}:=(1-\epsilon ^{\prime \prime })\mathbb{P}_{0}+\epsilon ^{\prime \prime }Q_{1},\ R_{2}:=(1-\epsilon ^{\prime \prime })\mathbb{P}_{\mu _{0}}+\epsilon ^{\prime \prime }Q _{2}.
\end{equation}
It is not hard to verify that actually $R _{1}=R _{2}$, i.e., they represent the same distribution. Consequently $\text{TV}\left(R _{1}^{\otimes N},R _{2}^{\otimes N}\right)=0$, where $R _{j}^{\otimes N}, j \in \left\{1,2\right\}$ is the $N$-fold product of $R _{j}$.\opar

\textbf{Step III: expand the TV distance.} Define 
\begin{equation*}
\binom{[N]}{s}:=\left\{(i _{1},\dots,i _{s}):i _{1},\dots,i _{s}\in \left\{1,2,\dots,N\right\}, 1 \le i _{1}<i _{2}<\dots < i _{s}\le N\right\}.
\end{equation*}
Further define the quantities
\begin{equation*}
\begin{aligned}
    &S _{1}^{L}:=\sum\limits_{s=0}^{N \left(\epsilon ^{\prime \prime }+\frac{1}{\sqrt{2N}}\right)}\sum\limits_{I \in \binom{[N]}{s}}^{}(1-\epsilon ^{\prime \prime })^{N-s}\epsilon ^{\prime \prime s}\prod\limits_{i \in [N]\backslash I}^{}p _{0}(x _{i})\prod\limits_{i \in I}^{}q _{1}(x _{i}),\\ 
    &S _{2}^{L}:=\sum\limits_{s=0}^{N \left(\epsilon ^{\prime \prime }+\frac{1}{\sqrt{2N}}\right)}\sum\limits_{I \in \binom{[N]}{s}}^{}(1-\epsilon ^{\prime \prime })^{N-s}\epsilon ^{\prime \prime s}\prod\limits_{i \in [N]\backslash I}^{}p _{\mu _{0}}(x _{i})\prod\limits_{i \in I}^{}q _{2}(x _{i}),\\ 
    &S _{1}^{U}:=\sum\limits_{s=N \left(\epsilon ^{\prime \prime }+\frac{1}{\sqrt{2N}}\right)+1}^{N}\sum\limits_{I \in \binom{[N]}{s}}^{}(1-\epsilon ^{\prime \prime })^{N-s}\epsilon ^{\prime \prime s}\prod\limits_{i \in [N]\backslash I}^{}p _{0}(x _{i})\prod\limits_{i \in I}^{}q _{1}(x _{i}),\\ 
    &S _{2}^{U}:=\sum\limits_{s=N \left(\epsilon ^{\prime \prime }+\frac{1}{\sqrt{2N}}\right)+1}^{N}\sum\limits_{I \in \binom{[N]}{s}}^{}(1-\epsilon ^{\prime \prime })^{N-s}\epsilon ^{\prime \prime s}\prod\limits_{i \in [N]\backslash I}^{}p _{\mu _{0}}(x _{i})\prod\limits_{i \in I}^{}q _{2}(x _{i}),\\ 
\end{aligned}
\end{equation*}
where without loss of generality we assume that all bounds of summation are integers. For $S _{1}^{U}$ and $S _{2}^{U}$, we have 
\begin{equation*}
\frac{1}{2}\displaystyle\int_{\mathbb{R}^{d}}^{}\left|S _{1}^{U}-S _{2}^{U}\right|d \mathbf{x} \le \frac{1}{2}\displaystyle\int_{\mathbb{R}^{d}}^{}(S _{1}^{U}+S _{2}^{U})d \mathbf{x}\le \mathbb{P}\left(B >N \left(\epsilon ^{\prime \prime }+\frac{1}{\sqrt{2N}}\right)\right)\overset{\text{(i)}}{\le }e ^{-1}.
\end{equation*}
Here $B$ is a binomial random variable with parameters $N$ and $\epsilon ^{\prime \prime }$, and (i) follows from the Hoeffding's inequality.

By expanding the TV distance between the $N$-fold products of $R _{1}$ and $R _{2}$ and after some calculation, we have 
\begin{equation*}
0=\text{TV}\left(R _{1}^{\otimes N},R _{2}^{\otimes N}\right)\overset{\text{(ii)}}{\ge }\frac{1}{2}\displaystyle\int_{\mathbb{R}^{d}}^{}\left|S _{1}^{L}-S _{2}^{L}\right|d \mathbf{x}-\frac{1}{2}\displaystyle\int_{\mathbb{R}^{d}}^{}\left|S _{1}^{U}-S _{2}^{U}\right|d \mathbf{x},
\end{equation*}
where (ii) is from the definition of TV distance and the triangle inequality. Combining with previous argument, we obtain 
\begin{equation*}
\frac{1}{2}\displaystyle\int_{\mathbb{R}^{d}}^{}\left|S _{1}^{L}-S _{2}^{L}\right|d \mathbf{x}\le e ^{-1}.
\end{equation*}

\textbf{Step IV: conditional argument.} Consider a adversary $\mathcal{C}_{B,\epsilon ^{\prime \prime }}$ with following strategy. Upon observing the original data $\tilde{\mathbf{X}}$ and its underlying law (which is either $\mathbb{P}_{0}$ or $\mathbb{P}_{\mu _{0}}$ in our case), the adversary first samples a number $n _{c}$ from the conditional law of $B|B \le N(\epsilon ^{\prime \prime }+1/\sqrt{2N})$, then sample $n _{c}$ replacements from the mixture distribution corresponding to the true underlying law. Specifically, if $\tilde{\mathbf{X}}$ is from $\mathbb{P}_{0}$, then replacements will be from $Q _{1}$, otherwise the replacements will be from $Q _{2}$. In this way, if we denote $\tilde{\mathbb{P}}_{0}$ to be the law of the contaminated samples from $\mathbb{P}_{0}$ and $Q _{1}$, and $\tilde{\mathbb{P}}_{\mu _{0}}$ to be the law of the contaminated samples from $\mathbb{P}_{\mu _{0}}$ and $Q _{2}$, then with minor calculation we have 
\begin{equation*}
\text{TV}(\tilde{\mathbb{P}}_{0},\tilde{\mathbb{P}}_{\mu _{0}})=\frac{\frac{1}{2}\displaystyle\int_{\mathbb{R}^{d}}^{}\left|S _{1}^{L}-S _{2}^{L}\right|d \mathbf{x}}{\mathbb{P}\left(B \le N \left(\epsilon ^{\prime \prime }+\frac{1}{\sqrt{2N}}\right)\right)}\overset{\text{(iii)}}{\le }\frac{e ^{-1}}{1-e ^{-1}}<1,
\end{equation*}
where (iii) is again from the Hoeffding's inequality and the argument from the step III.\opar

\textbf{Step V: finish the proof.} Now for any $\tilde{\mu }$ that $\left\lVert \tilde{\mu }\right\rVert_{2}^{}\le \left\lVert \mu _{0}\right\rVert_{2}^{}=\frac{2 \sigma \epsilon ^{\prime }}{1-\epsilon ^{\prime }}$. Following a similar argument in the \hyperref[subsection: proof of the main lower bound 1]{proof} of Theorem \ref{theorem: lower bound 1}, for any $\phi \in A _{s}$, we have
\begin{equation*}
\begin{aligned}
    \sup\limits_{\mu \in K,\left\lVert \mu \right\rVert_{2}^{}\ge \rho }\sup\limits_{\mathcal{C}}\mathbb{P}_{\mu }\left(\phi (\mathcal{C}(\mathbf{X}^{\prime }))=0\right)&\ge \mathbb{P}_{\tilde{\mu }}\left(\phi (\mathcal{C}_{B,\epsilon ^{\prime \prime }}(\mathbf{X}^{\prime }))=0\right)\\ 
    &=1-\mathbb{P}_{\tilde{\mu }}\left(\phi (\mathcal{C}_{B,\epsilon ^{\prime \prime }}(\mathbf{X}^{\prime }))=1\right)+\mathbb{P}_{0}\left(\phi (\mathcal{C}_{B,\epsilon ^{\prime \prime }}(\mathbf{X}^{\prime }))=1\right)\\
    &-\mathbb{P}_{0}\left(\phi (\mathcal{C}_{B,\epsilon ^{\prime \prime }}(\mathbf{X}^{\prime }))=1\right)\\ 
    &\ge 1-\alpha -\left(\mathbb{P}_{\tilde{\mu }}\left(\phi (\mathcal{C}_{B,\epsilon ^{\prime \prime }}(\mathbf{X}^{\prime }))=1\right)-\mathbb{P}_{0}\left(\phi (\mathcal{C}_{B,\epsilon ^{\prime \prime }}(\mathbf{X}^{\prime }))=1\right)\right)\\ 
    &=1-\alpha -\left(\tilde{\mathbb{P}}_{\tilde{\mu }}\left(\phi =1\right)-\tilde{\mathbb{P}}_{0}(\phi =1)\right)\\ 
    &\ge 1-\alpha -\text{TV}\left(\tilde{\mathbb{P}}_{0},\tilde{\mathbb{P}}_{\tilde{\mu }}\right)\\ 
    &\ge 1-\alpha -\text{TV}\left(\tilde{\mathbb{P}}_{0},\tilde{\mathbb{P}}_{\mu _{0}}\right)=\varTheta (1).
\end{aligned}
\end{equation*}
Taking infimum over $\phi \in A _{s}$ on both sides completes the proof.
\end{proof}

\subsection{Proof of Theorem \ref{theorem: upper bound 1}}\label{subsection: proof of the main upper bound 1}

We start with a recall of the definition of the consistent subset, which plays a central role in the proof.

\begin{definition}[Consistent subset]\label{def_consistent_subset}
    A subset $S \subset \mathbf{X}=\left\{X _{1},\dots,X _{N}\right\}$ is called a consistent subset regarding the contamination fraction $\epsilon $ and a test $\phi :2 ^{S}\mapsto \left\{0,1\right\}$ if\\ 
    \hspace*{0.5em}(i), $\left|S\right|\ge (1-\epsilon )N$,\\ 
    \hspace*{0.5em}(ii), $\phi (S ^{\prime })=\phi (S)$ for any $S ^{\prime }\subset S$ with $\left|S ^{\prime }\right|\ge (1-2\epsilon )N$.
\end{definition}

We are now ready to prove Theorem \ref{theorem: upper bound 1}. To provide a high-level outline, we first establish the high-probability existence of a consistent subset with respect to $\phi _{e}$ implied by $E _{e}$ defined in \eqref{eq: definition of the indicator event}. This is done by expanding the $\ell _{2}$ norm of the projection vector in $E _{e}$ and bounding each component of the expansion. The rationale of $\phi _{e}$ then follows from the nature of the consistent subset. In the proof, $k ^{\star }\in [1,d]$\footnote{$1\le k ^{\star }<\infty $ if $K \subset \ell _{2}$.} is any fixed dimension that does not necessary to be $k _{1}^{\star }$ or $k _{2}^{\star }$, and $P ^{\star }$ is its corresponding projection operator determined by the Kolmogorov widths. The pseudo code of $\phi _{e}$ is attached in Appendix \ref{section: algorithms}.

\begin{proof}[Proof of Theorem \ref{theorem: upper bound 1}]

We will first prove the fact that for a properly fixed error tolerance $\alpha $, the consistent subset exists with probability at least $1-\alpha $. Specifically, consider the following $\chi ^{2}$-type event
\begin{equation}\label{eq: definition of the indicator event}
    E _{e}:=\left\{\left\lVert P ^{\star }\sum\limits_{i \in T}^{}X _{i}\right\rVert_{2}^{2}-k ^{\star }\left|T\right|\sigma ^{2}\ge c\left|T\right|^{2}E^{2}(\epsilon ,K,N,\sigma ^{2})\right\},
\end{equation}
where $T$ is any subset of $\mathbf{X}$, and $c \in \mathbb{R}$ is some universal, properly chosen constant. On the subset $T$, the test $\phi _{e}$ rejects $H _{0}$ when $E _{e}$ holds and fails to reject otherwise.

For the rest of the proof, let us first focus on $\tilde{\mathbf{X}}$ and any its subset $\tilde{\mathbf{X}}_{S ^{\prime }}$ with $\left|S ^{\prime }\right|\ge (1-2\epsilon )N$. We will prove that with $\tilde{\mathbf{X}}$ is consistent with all such $S ^{\prime }$ if \eqref{eq: required condition of the main upper bound 1} is satisfied. The key to the proof is carefully analyzing the $\ell _{2}$ norm of the sum of the projected original samples 
\begin{equation}\label{eq: decomposition of projection}
    \left\lVert P ^{\star }\sum\limits_{i \in S ^{\prime }}^{}\tilde{X} _{i}\right\rVert_{2}^{2}=\underbrace{\left\lVert P ^{\star }\sum\limits_{i \in S ^{\prime }}^{}(\tilde{X} _{i}-\mu )\right\rVert_{2}^{2}}_{:=\text{I}}+\underbrace{\left|S ^{\prime }\right|^{2}\left\lVert P ^{\star }\mu \right\rVert_{2}^{2}}_{:=\text{II}}+2 \underbrace{\left|S ^{\prime }\right|\sum\limits_{i \in S ^{\prime }}^{}\left[P ^{\star }(\tilde{X} _{i}-\mu )\right]^{\top } \mu }_{:=\text{III}}.
\end{equation}
We attempt to show that the parts I and III, which are $\ell _{2}$ norm of the centralized projected vectors and the crossing term, are both asymptotically less than the part II under the condition (\ref{eq: required condition of the main upper bound 1}), which vanishes under the null.

\textbf{Part I.} By representing $S ^{\prime }=[N]\backslash \left([N]\backslash S ^{\prime }\right)$, we have 
\begin{equation*}
    \text{I}=\underbrace{\left\lVert P ^{\star }\sum\limits_{i \in [N]}^{}(\tilde{X} _{i}-\mu )\right\rVert_{2}^{2}}_{:=\text{I}_{1}}+\underbrace{\left\lVert P ^{\star }\sum\limits_{i \in [N] \backslash S ^{\prime }}^{}(\tilde{X} _{i}-\mu )\right\rVert_{2}^{2}}_{:=\text{I}_{2,S ^{\prime }}}-2 \underbrace{\sum\limits_{i \in [N]}^{}(\tilde{X} _{i}-\mu )^{\top } P ^{\star }\sum\limits_{i \in [N]\backslash S ^{\prime }}^{}(\tilde{X} _{i}-\mu )}_{:=\text{I}_{3,S ^{\prime }}}.
\end{equation*}
For $\text{I}_{1}$, from the property of the Gaussian distribution, we know $\text{I}_{1}\sim N \sigma ^{2}\cdot \chi _{k ^{\star }}^{2}$. Therefore $\mathbb{E}\text{I} _{1}=k ^{\star }N \sigma ^{2}, \text{Var}\left(\text{I} _{1}\right)=2k ^{\star }N ^{2}\sigma ^{4}$. By the Chebyshev's inequality we know that 
\begin{equation*}
    \left|\text{I}_{1}-k ^{\star }N \sigma ^{2}\right|\le \sqrt{\frac{2k ^{\star }}{\alpha }} N \sigma ^{2}:=t _{\text{I}_{1}}.
\end{equation*}
with probability greater than $1-\alpha $.

For $\text{I}_{2,S ^{\prime }}$, since we allow $S ^{\prime }$ to be any subset of $[N]$, refined analysis is needed. Using the Hanson--Wright inequality \ref{lemma: hw_inequality} and observing that $\left\lVert P ^{\star }\right\rVert_{F}^{}=\sqrt{k ^{\star }},\left\lVert P ^{\star }\right\rVert_{2}^{}\le 1$, for each fixed $S ^{\prime }\subset [N]$, we have 
\begin{equation*}
    \mathbb{P}\left(\left|\text{I} _{2}-k ^{\star }\left|[N]\backslash S ^{\prime }\right|\sigma ^{2}\right| \ge t\right)\le 2\exp\left\{-c \min\limits \left\{\frac{t ^{2}}{k ^{\star }\left|[N]\backslash S ^{\prime }\right|^{2}\sigma ^{4}},\frac{t}{\left|[N]\backslash S ^{\prime }\right|\sigma ^{2}}\right\}\right\}.
\end{equation*}
Since $\left|S ^{\prime }\right|\ge (1-2\epsilon )N$, we have $\sum\limits_{i=0}^{2\epsilon N}\binom{N}{i}=\left(\frac{e}{2\epsilon }\right)^{2\epsilon N}$ possible $[N]\backslash S ^{\prime }$ by Lemma \ref{lemma: combination number}. Using the union bound, we have $\text{I} _{2,S ^{\prime }}-k ^{\star }\left|[N]\backslash S ^{\prime }\right|\sigma ^{2}\ge t$ holds for some $S ^{\prime }$ with probability less than
\begin{equation*}
    \sum\limits_{S ^{\prime }}^{}\exp\left\{-c \min\limits \left\{\frac{t ^{2}}{k ^{\star }\left|[N]\backslash S ^{\prime }\right|^{2}\sigma ^{4}},\frac{t}{\left|[N]\backslash S ^{\prime }\right|\sigma ^{2}}\right\}\right\}\le \left(\frac{e}{2\epsilon }\right)^{2\epsilon N}\cdot \exp\left\{-c \min\limits \left\{\frac{t ^{2}}{k ^{\star }\epsilon ^{2}N^{2}\sigma ^{4}},\frac{t}{\epsilon N\sigma ^{2}}\right\}\right\}.
\end{equation*}
If we select $t _{\text{I}_{2}}$ to be
\begin{equation*}
    t _{\text{I}_{2}}=\max\limits \left\{\sqrt{\frac{4\left[2\epsilon N \ln \left(\frac{e}{2\epsilon }\right)+\ln \left(\frac{2}{\alpha }\right)\right]k ^{\star }\epsilon ^{2}N ^{2}\sigma ^{4}}{c}},\frac{2\left[2\epsilon N \ln \left(\frac{e}{\epsilon }\right)+\ln \left(\frac{2}{\alpha }\right) \right]\epsilon N \sigma ^{2}}{c}\right\}
\end{equation*}
then it can be verified that $\left|\text{I} _{2,S ^{\prime }}-k ^{\star }\left|[N]\backslash S ^{\prime }\right|\sigma ^{2}\right|\le  t _{\text{I}_{2}}$ holds for all $S ^{\prime }$ with probability greater than $1-\alpha $.\npar

For $\text{I}_{3,S ^{\prime }}$, we use the following argument of conditional probability to effectively bound its tail. First, observe that the joint distribution of $\sum\limits_{i \in [N]\backslash S ^{\prime }}^{}(\tilde{X} _{i}-\mu )$ and $\sum\limits_{i \in [N]}^{}(\tilde{X} _{i}-\mu )$ is a multivariate Gaussian
\begin{equation*}
    \left(\sum\limits_{i \in [N]\backslash S ^{\prime }}^{}(\tilde{X} _{i}-\mu )^{\top } ,\sum\limits_{i \in [N]}^{}(\tilde{X} _{i}-\mu )^{\top } \right)\sim \mathcal{N}\left(0,\begin{pmatrix}
        \left|[N]\backslash S ^{\prime }\right| \mathbf{I}_{d} & \left|[N]\backslash S ^{\prime }\right| \mathbf{I}_{d}\\ 
        \left|[N]\backslash S ^{\prime }\right| \mathbf{I}_{d} & N \mathbf{I}_{d}
    \end{pmatrix}\sigma ^{2}\right).
\end{equation*}
By the property of multivariate Gaussian distribution, conditional on the event $\left\{\sum\limits_{i \in [N]}^{}(\tilde{X} _{i}-\mu )=x\right\}$, the conditional distribution of $\sum\limits_{i \in [N]\backslash S ^{\prime }}^{}(\tilde{X} _{i}-\mu )$ is
\begin{equation*}
    \left[\sum\limits_{i \in [N]\backslash S ^{\prime }}^{}(\tilde{X} _{i}-\mu )\left\lvert\right. \sum\limits_{i \in [N]}^{}(\tilde{X} _{i}-\mu )=x \right]\sim \mathcal{N}\left(\frac{\left|[N]\backslash S ^{\prime }\right|}{N}x,\sigma ^{2}\left(\left|[N]\backslash S ^{\prime }\right|-\frac{\left|[N]\backslash S ^{\prime }\right|^{2}}{N}\right)\mathbf{I}_{d}\right).
\end{equation*}
By the maximal inequality for (sub-)Gaussian random variable that $\mathbb{P}\left(\max\limits _{1 \le i \le N}\tilde{X} _{i}\ge \sqrt{2 \sigma ^{2}(\ln N+t)}\right)\le e ^{-t}$, and conditional on the event $\left\{\sum\limits_{i \in [N]}^{}(\tilde{X} _{i}-\mu )=x\right\}$, we have the probability of the following event 
\begin{equation*}
    \left\{\text{I}_{3,S ^{\prime }}\ge \frac{\left|[N]\backslash S ^{\prime }\right|}{N}x ^{\top } P ^{\star }x+\sqrt{2x ^{\top } P ^{\star }x \left(\left|[N]\backslash S ^{\prime }\right|-\frac{\left|[N]\backslash S ^{\prime }\right|^{2}}{N}\right)\sigma ^{2}\cdot \left(2\epsilon N\ln \left(\frac{e}{2\epsilon }\right)+t\right)}\text{ for some $S ^{\prime }\subset S$}\right\}
\end{equation*}
less than $e ^{-t}$. We also need the following basic inequality of conditional probability for any random variables $X$ and $Y$
\begin{equation*}
    \begin{aligned}
        \mathbb{P}\left(Y \ge a\right)&=\mathbb{E}\mathbf{1}_{\left\{Y \ge a\right\}}\\ 
        &=\mathbb{E}\left[\mathbf{1}_{\left\{Y \ge a\right\}}\cdot \left(\mathbf{1}_{\left\{f(X)\ge a\right\}}+\mathbf{1}_{\left\{f(X)<a\right\}}\right)\right]\\ 
        &=\mathbb{E}\left[\mathbb{E}\left[\mathbf{1}_{\left\{Y \ge a\right\}}\cdot \mathbf{1}_{\left\{f(X)\ge a\right\}}\left|\right.X\right]+\mathbb{E}\left[\mathbf{1}_{\left\{Y \ge a\right\}}\cdot \mathbf{1}_{\left\{f(X)<a\right\}}\left|\right.X\right]\right]\\ 
        &\le \mathbb{E}\left[\mathbb{E}\left[\mathbf{1}_{\left\{f(X)\ge a\right\}}\left|\right.X\right]+\mathbb{E}\left[\mathbf{1}_{\left\{Y \ge f(X)\right\}}\cdot \mathbf{1}_{\left\{f(X)<a\right\}}\left|\right.X\right]\right]\\ 
        &\le \mathbb{P}\left(f(X)\ge a\right)+\mathbb{E}\left[\mathbb{P}\left(Y \ge f(X)\left|\right.X\right)\right].
    \end{aligned}
\end{equation*}
Now observe that $\left|[N]\backslash S ^{\prime }\right|\le 2\epsilon N$ for all $S ^{\prime }$, and take $Y=\max\limits _{S ^{\prime }}\left\{\text{I}_{3,S ^{\prime }}\right\}$, $X=\sum\limits_{i \in [N]}^{}(\tilde{X} _{i}-\mu )$, $f(x)=g(x ^{\top } P ^{\star }x)$ where $g(x)=\frac{\left|[N]\backslash S ^{\prime }\right|}{N}x+\sqrt{4x \cdot \epsilon N \left(\epsilon N \ln \left(\frac{e}{\epsilon }\right)+t\right)}$ and $a=g(Nk ^{\star }\sigma ^{2}+\sqrt{\frac{4k ^{\star }}{\alpha }}N \sigma ^{2})$. By the argument above and the fact that $g(\cdot )$ is an increasing function on the positive real line, we know
\begin{equation*}
    \mathbb{P}\left(f(X)\ge a\right)=\mathbb{P}\left(x ^{\top } P ^{\star }x \ge Nk ^{\star }\sigma ^{2}+\sqrt{\frac{4k ^{\star }}{\alpha }}N \sigma ^{2}\right)\le \frac{\alpha }{2} .
\end{equation*}
From the maximal inequality above, we also have $\mathbb{E}\left[\mathbb{P}\left(Y \ge f(X)\left|\right.X\right)\right]\le \mathbb{E}e ^{-t}=e ^{-t}$. Taking $t=\ln \left(\frac{2}{\alpha }\right)$, and combining the results above, we conclude that with probability greater than $1-\alpha $ the following event holds for all $S ^{\prime }$.
\begin{equation*}
    \left\{\text{I}_{3,S ^{\prime }}\le \left|[N]\backslash S ^{\prime }\right|k ^{\star }\sigma ^{2}+\left|[N]\backslash S ^{\prime }\right|\sqrt{\frac{4k ^{\star }}{\alpha }}\sigma ^{2}+\sqrt{8 \left(Nk ^{\star }\sigma ^{2}+\sqrt{\frac{4k ^{\star }}{\alpha }}N \sigma ^{2}\right)\epsilon N \left(2\epsilon N \ln \left(\frac{e}{2\epsilon }\right)+\ln \left(\frac{2}{\alpha }\right)\right)\sigma ^{2}}\right\}.
\end{equation*}
By the symmetrical nature of $\text{I}_{3,S ^{\prime }}$, same technique and argument can also be applied to the probability of the lower tail of $\text{III}_{3,S ^{\prime }}$.

\textbf{Part III.} We rewrite the part III as
\begin{equation*}
    \text{III}=\underbrace{\left|S ^{\prime }\right|\sum\limits_{i \in [N]}^{}(\tilde{X} _{i}-\mu )^{\top }P ^{\star }\mu }_{\text{III}_{1,S ^{\prime }}}-\underbrace{\left|S ^{\prime }\right|\sum\limits_{i \in [N]\backslash S ^{\prime }}^{}(\tilde{X} _{i}-\mu )^{\top } P ^{\star }\mu }_{\text{III}_{2,S ^{\prime }}}.
\end{equation*}
Dealing with $\text{III}_{1,S ^{\prime }}$ is relatively easier. Observe that $\text{III}_{1,S ^{\prime }}$ actually follows a Gaussian distribution with $\mathbb{E}\text{III}_{1,S ^{\prime }}=0$ and $\text{Var}\left(\text{III}_{1,S ^{\prime }}\right)=\left|S ^{\prime }\right|^{2}N \sigma ^{2}\left\lVert P ^{\star }\mu \right\rVert_{2}^{2}\le N ^{3}\sigma ^{2}\left\lVert P ^{\star }\mu \right\rVert_{2}^{2}$. By the Chebyshev's inequality, we have
\begin{equation*}
    \mathbb{P}\left(\left|\text{III}_{1, S ^{\prime }}\right|\le \sqrt{\frac{N ^{3}\sigma ^{2}\left\lVert P ^{\star }\mu \right\rVert_{2}^{2}}{\alpha }}\right)\ge 1-\alpha .
\end{equation*}
For $\text{III}_{2,S ^{\prime }}$, since the summation involves $S ^{\prime }$, we again need to exploit the union bound (or equivalently the maximal inequality). $\text{III}_{2,S ^{\prime }}$ also follows a Gaussian distribution with $\mathbb{E}\text{III}_{2,S ^{\prime }}=0$ and $\text{Var}\left(\text{III}_{2,S ^{\prime }}\right)=\left|S ^{\prime }\right|^{2} \left|[N]\backslash S ^{\prime }\right|\left\lVert P ^{\star }\mu \right\rVert_{2}^{2}\sigma ^{2}\le 2\epsilon N ^{3}\left\lVert P ^{\star }\mu \right\rVert_{2}^{2}\sigma ^{2}$. By the Chebyshev's inequality and the union bound, we know that the following holds.
\begin{equation*}
    \mathbb{P}\left(\left|\text{III}_{2,S ^{\prime }}\right|\le \sqrt{2\epsilon N ^{3}\left\lVert P ^{\star }\mu \right\rVert_{2}^{2}\sigma ^{2}}\cdot \sqrt{2\epsilon N \ln \left(\frac{e}{2\epsilon }\right)+\ln \left(\frac{1}{\alpha }\right)},\forall S ^{\prime }\subset S\right)\ge 1-\alpha .
\end{equation*}
In summary, each tail behaviour of the term above consists of two parts: the mean of the term and the variation around the mean. To clearly display the results above, the following table can be helpful.
\begin{table}
    \begin{tabular}{|l|l|l|l|}
        \hline
        \textbf{Terms} & \textbf{Mean} & \textbf{Variation} & \textbf{Tail}\\ 
        \hline
        $\text{I}_{1}$ & $k ^{\star }N \sigma ^{2}$& $\sqrt{\frac{2k ^{\star }}{\alpha }}N \sigma ^{2}$ & $\alpha $\\ 
        \hline
        $\text{I}_{2,S ^{\prime }}$ & $k ^{\star }\left|[N]\backslash S ^{\prime }\right|\sigma ^{2}$ & $\max\limits \left\{\sqrt{\frac{4\left[2\epsilon N \ln \left(\frac{e}{2\epsilon }\right)+\ln \left(\frac{2}{\alpha }\right)\right]k ^{\star }\epsilon ^{2}N ^{2}\sigma ^{4}}{c}},\frac{2\left[2\epsilon N \ln \left(\frac{e}{\epsilon }\right)+\ln \left(\frac{2}{\alpha }\right) \right]\epsilon N \sigma ^{2}}{c}\right\}$ & $\alpha $\\ 
        \hline
        $\text{I}_{3,S ^{\prime }}$ & $k ^{\star }\left|[N]\backslash S ^{\prime }\right|\sigma ^{2}$ &$\left[\left|[N]\backslash S ^{\prime }\right|\sqrt{\frac{4k ^{\star }}{\alpha }}+\sqrt{8 \left(Nk ^{\star }+\sqrt{\frac{4k ^{\star }}{\alpha }}N \right)\epsilon N \left(2\epsilon N \ln \left(\frac{e}{2\epsilon }\right)+\ln \left(\frac{2}{\alpha }\right)\right)}\right]\sigma ^{2} $ & $2\alpha $\\ 
        \hline
        $\text{III}_{1,S ^{\prime }}$ & $0$ & $\sqrt{\frac{N ^{3}\sigma ^{2}\left\lVert P ^{\star }\mu \right\rVert_{2}^{2}}{\alpha }}$ & $\alpha $\\ 
        \hline
        $\text{III}_{2,S ^{\prime }}$ & $0$ & $\sqrt{2\epsilon N ^{3}\left\lVert P ^{\star }\mu \right\rVert_{2}^{2}\sigma ^{2}}\cdot \sqrt{2\epsilon N \ln \left(\frac{e}{2\epsilon }\right)}$ & $\alpha $\\
        \hline
    \end{tabular}
    \caption{corresponding means and variations of the terms from \eqref{eq: decomposition of projection}.}
\end{table}
To complete the proof, we finally need to compare the variation of each term above with the part II.\npar

\hspace*{1em}$\text{I}_{1}$. Require: $\left|S ^{\prime }\right|^{2}\left\lVert P ^{\star }\mu \right\rVert_{2}^{2}\gtrsim \sqrt{\frac{2k ^{\star }}{\alpha }}N \sigma ^{2}$. This is equivalent to $\left\lVert P ^{\star }\mu \right\rVert_{2}^{2}\gtrsim \frac{\sqrt{k ^{\star }}}{N}\sigma ^{2}$.\npar

\hspace*{1em}$\text{I}_{2,S ^{\prime }}$. Require: $\left|S ^{\prime }\right|^{2}\left\lVert P ^{\star }\mu \right\rVert_{2}^{2}\gtrsim \max\limits \left\{\sqrt{\frac{4\left[2\epsilon N \ln \left(\frac{e}{\epsilon }\right)+\ln \left(\frac{2}{\alpha }\right)\right]k ^{\star }\epsilon ^{2}N ^{2}\sigma ^{4}}{c}}, \frac{2\left[2\epsilon N \ln \left(\frac{e}{\epsilon }\right)+\ln \left(\frac{2}{\alpha }\right) \right]\epsilon N \sigma ^{2}}{c}\right\}$. This is equivalent to $\left\lVert P ^{\star }\mu \right\rVert_{2}^{2}\gtrsim \sqrt{\frac{\epsilon ^{3}\ln \left(\frac{1}{\epsilon }\right)k ^{\star }\sigma ^{4}}{N}}$ and $\left\lVert P ^{\star }\mu \right\rVert_{2}^{2}\gtrsim \epsilon ^{2}\ln \left(\frac{1}{\epsilon }\right)\sigma ^{2}$.\npar

\hspace*{1em}$\text{I}_{3,S ^{\prime }}$. Require: $\left|S ^{\prime }\right|^{2}\left\lVert P ^{\star }\mu \right\rVert_{2}^{2}\gtrsim \left|[N]\backslash S ^{\prime }\right|\sqrt{\frac{4k ^{\star }}{\alpha }}\sigma ^{2}\vee \sqrt{8 \left(Nk ^{\star }+\sqrt{\frac{4k ^{\star }}{\alpha }}N\right)\epsilon N \left(2\epsilon N \ln \left(\frac{1}{2\epsilon }\right)+\ln \left(\frac{2}{\alpha }\right)\right)}\sigma ^{2}$. This is equivalent to $\left\lVert P ^{\star }\mu \right\rVert_{2}^{2}\gtrsim \frac{\epsilon \sqrt{k ^{\star }}}{N}\sigma ^{2}$ and $\left\lVert P ^{\star }\mu \right\rVert_{2}^{2}\gtrsim \sqrt{\frac{\epsilon ^{2}\ln \left(\frac{1}{\epsilon }\right)k ^{\star }\sigma ^{4}}{N}}$.\npar

\hspace*{1em}$\text{III}_{1,S ^{\prime }}$. Require: $\left|S ^{\prime }\right|^{2}\left\lVert P ^{\star }\mu \right\rVert_{2}^{2}\gtrsim \sqrt{\frac{N ^{3}\sigma ^{2}\left\lVert P ^{\star }\mu \right\rVert_{2}^{2}}{\alpha }}$. This is equivalent to $\left\lVert P ^{\star }\mu \right\rVert_{2}^{2}\gtrsim \frac{\sigma ^{2}}{N}$.\npar

\hspace*{1em}$\text{III}_{2,S ^{\prime }}$. Require: $\left|S ^{\prime }\right|^{2} \left\lVert P ^{\star }\mu \right\rVert_{2}^{2}\gtrsim \sqrt{2\epsilon N ^{3}\left\lVert P ^{\star }\mu \right\rVert_{2}^{2}\sigma ^{2}}\cdot \sqrt{2\epsilon N \ln \left(\frac{e}{2\epsilon }\right)+\ln \left(\frac{1}{\alpha }\right)}$. This is equivalent to $\left\lVert P ^{\star }\mu \right\rVert_{2}^{2}\gtrsim \max\limits \left\{\epsilon ^{2}\ln \left(\frac{1}{\epsilon }\right)\sigma ^{2},\frac{\epsilon \sigma ^{2}}{N}\right\}$.\npar

In the argument above, we have used the assumption that $(1-2\epsilon )N \le \left|S ^{\prime }\right|\le N$ and therefore $0 \le \left|[N]\backslash S ^{\prime }\right|\le 2\epsilon N$. Summarizing the table and the argument above, it is required that 
\begin{equation*}
    \left\lVert P ^{\star }\mu \right\rVert_{2}^{2}\gtrsim \sigma ^{2}\max\limits \left\{\frac{\sqrt{k ^{\star }}}{N},\epsilon ^{2}\ln \left(\frac{1}{\epsilon }\right),\sqrt{\frac{\epsilon ^{2}\ln \left(\frac{1}{\epsilon }\right)k ^{\star }}{N}}\right\}.
\end{equation*}
Finally, by the Pythagorean theorem, we have $\left\lVert P ^{\star }\mu \right\rVert_{2}^{2}=\rho ^{2}-\left\lVert \mu -P^{\star }\mu \right\rVert_{2}^{2}\ge \rho ^{2}-D _{k ^{\star }}^{2}(K)$. Therefore, in order that $\left\lVert P _{k}^{\star }\mu \right\rVert_{2}^{2}$ satifies the condition above, it suffices that 
\begin{equation*}
    \rho ^{2}\gtrsim D _{k ^{\star }}^{2}(K)+\sigma ^{2}\max\limits \left\{\frac{\sqrt{k ^{\star }}}{N},\epsilon ^{2}\ln \left(\frac{1}{\epsilon }\right),\sqrt{\frac{\epsilon ^{2}\ln \left(\frac{1}{\epsilon }\right)k ^{\star }}{N}}\right\}=E _{\text{raw}}^{2}(\epsilon ,K,N,\sigma ^{2}).
\end{equation*}

In fact, we have proved that there exists consistent subsets within $\mathbf{X}$ (at least $\mathbf{X}_{[N]\backslash C}$) for $\phi _{e}$ with high probability. The remaining part goes naturally. For the observed dataset $\mathbf{X}$, we scan through all $S \subset [N]$ with $\left|S\right|\ge (1-\epsilon )N$ and perform $\phi _{e}$ on $\mathbf{X}_{S}$ and $\mathbf{X}_{S ^{\prime }}$ with $S ^{\prime }\subset S, \left|S\right|^{\prime }\ge (1-2\epsilon )N$ to check if the subset is consistent. The proof above guarantees that we can eventually find a consistent subset $S _{0}$ (with high probability), though $\mathbf{X}_{S _{0}}$ may contain corrupted samples. However, once found, we simply perform the test $\phi _{e}$ on $\mathbf{X}_{S _{0}}$ and accept or reject $H _{0}$ depending on the testing result. By the assumption on the cardinality, we know $\left|S _{0} \cap [N]\backslash C\right|\ge (1-2\epsilon )N$. The proof and the assumption above show that with probability greater than $1-\alpha $, we have $\phi _{e}(\mathbf{X}_{S _{0}})=\phi _{e}(\mathbf{X}_{S _{0} \cap [N]\backslash C})=\phi _{e}(\mathbf{X}_{[N]\backslash C})=\phi _{e}(\tilde{\mathbf{X}})$. Therefore, testing on $\mathbf{X}_{S _{0}}$ is as testing on $\tilde{\mathbf{X}}$. Since $\phi _{e}$ itself is a valid test for the uncontaminated testing problem \eqref{eq: original testing problem}, the proof is completed.
\end{proof}

\subsection{Proof of Theorem \ref{theorem: upper bound 2}}\label{subsection: proof of the main upper bound 2}

Motivated by \cite{10353143}, we demonstrate the rationale of Algorithm \ref{algorithm: polynomial algorithm} in this section. Same as the \hyperref[subsection: proof of the main upper bound 2]{proof} of Theorem \ref{theorem: upper bound 1}, $k ^{\star }$ is any fixed dimension and $P ^{\star }$ is its corresponding projection operator. For the conciseness of the notation, we define $\tilde{Y}_{i}:=P ^{\star }\tilde{X}_{i}$, $Y _{i}:=P ^{\star }X _{i}$, and $\nu :=P ^{\star }\mu $. We assume without loss of generality that $\tilde{Y}_{i}, Y _{i}\in \mathbb{R}^{k ^{\star }}$ and $\text{Var}\left(\tilde{Y}_{i}\right)=\mathbf{I}_{k ^{\star }}$ with the fact in mind that any projected multivariate Gaussian random variable can be regarded as embedded into a lower dimension subspace with covariance matrix $\text{diag}\{(1,...,1,0,\dots,0)\}$ by rotating the axes. We also assume $\sigma =1$ in this section (one can scale the results by $\sigma $ anytime if needed).

To clarify the structure of this section and highlight how the definitions, lemmas, and algorithms fit together, we provide the following roadmap. We depart from the definition of regularity, and arrive at the final complete testing procedure, Algorithm \ref{algorithm: polynomial algorithm}. The pseudo code is provided in Appendix \ref{section: algorithms}.

\begin{figure}[htbp]
    \centering
    \includegraphics[width=\textwidth]{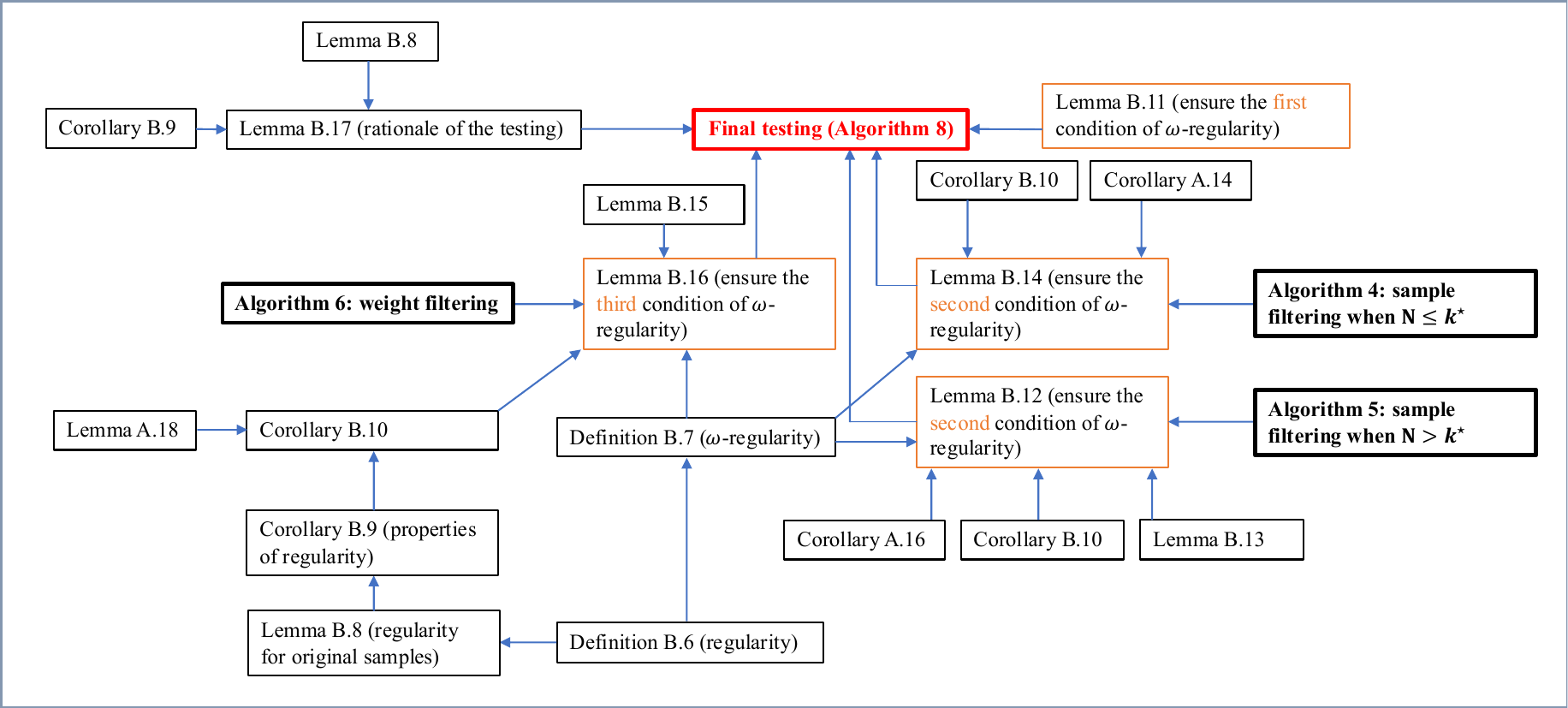}
    \caption{roadmap of this section.}\label{fig: tour map}
\end{figure}

We start with the definition of regularity and $\omega $-regularity, which controls the behaviour of the $\ell _{2}$ norm of the individual and sum, and the inner product.

\begin{definition}[regularity]\label{def_regularity}
    $\mathbf{Y}=\left\{Y _{1},\dots,Y _{N}\right\}$ is said to be $(\epsilon ,\beta _{1},\beta _{2})$-regular if for all subsets $S \subset [N]$ with $\left|S\right|\le \epsilon N$, we have the following properties:\\
    \begin{tabular}{rl}
        (i), & $\left|\sum\limits_{i \in S}^{}\left\lVert Y _{i}\right\rVert_{2}^{2}-\left|S\right|k ^{\star }\right|\le c \beta _{1}$,\\ 
        (ii), & $\left|\left\lVert \sum\limits_{i \in S}^{}Y _{i}\right\rVert_{2}^{2}-\left|S\right|k ^{\star }\right|\le c \beta _{2}$, and\\ 
        (iii), & $\left|\left\langle \sum\limits_{i \in S}^{}Y _{i}, \sum\limits_{j \in [N]}^{}Y _{j} \right\rangle-\left|S\right|k ^{\star }\right|\le c \sqrt{N}\beta _{1}.$
        \end{tabular}
\end{definition}

$\omega $-regularity is a natural generalization of the definition above, which involves the interaction of a weight vector $\omega $ with $\mathbf{Y}$.
\begin{definition}[$\omega $-regularity]\label{def_omega_regularity}
    Given a weight vector $\omega =(\omega _{1},\dots,\omega _{N})^{\top } $, $\mathbf{Y}$ is said to be $(\epsilon ,\beta _{1},\beta _{2})$-regular if for all subsets $S \subset [N]$ with $\left|S\right|\le \epsilon N$, we have the following properties:\\ 
    \begin{tabular}{rl}
        (i), & $\left|\sum\limits_{i \in S}^{}\left\lVert Y _{i}\right\rVert_{2}^{2}-\left|S\right|k ^{\star }\right|\le c \beta _{1}$,\\ 
        (ii), & $\left|\left\lVert \sum\limits_{i \in S}^{}\sqrt{\omega }_{i} Y _{i}\right\rVert_{2}^{2}-\left\lVert \omega _{S}\right\rVert_{1}^{}k ^{\star }\right|\le c \beta _{2}$, and\\ 
        (iii), & $\left|\left\langle \sum\limits_{i \in S}^{}\sqrt{\omega }_{i}Y _{i}, \sum\limits_{j \in [N]}^{}\sqrt{\omega }_{j}Y _{j} \right\rangle-\left\lVert \omega _{S}\right\rVert_{1}^{} k ^{\star }\right|\le c \sqrt{N}\beta _{1}.$
    \end{tabular}

    Recall that $\omega _{S} \in \mathbb{R}^{N}$ means the restriction of $\omega $ on the set $S$.
\end{definition}
\ref{def_omega_regularity} immediately recovers \ref{def_regularity} if $\omega =\mathbf{1}_{S}$.

The next lemma states that any uncontaminated $\mathbf{Y}$ from $\mathcal{N}\left(\nu ,\mathbf{I}_{k ^{\star }}\right)$ is in fact $(\epsilon ,\beta _{1},\beta _{2})$-regular with high probability and specified $\beta _{1}$ and $\beta _{2}$ below depending on $\epsilon ,\left\lVert \nu \right\rVert_{2}^{},N,k ^{\star }$.

\begin{lemma}[]\label{lemma: good_samples_regularity}
    Assume $Y _{1},\dots,Y _{N}$ are i.i.d. samples from $\mathcal{N}\left(\nu ,\mathbf{I}_{k ^{\star }}\right)$, then with probability greater than $1-\alpha $, $\mathbf{Y}=\left\{Y _{1},\dots,Y _{N}\right\}$ is $(\epsilon ,\beta _{1},\beta _{2})$-regular with 
    \begin{equation}\label{eq: specified beta_1 and beta_2}
        \begin{aligned}
            &\beta _{1}=\epsilon N \left(\sqrt{k ^{\star }}+\sqrt{N}\left\lVert \nu \right\rVert_{2}^{}\right)\sqrt{\ln \left(\frac{N}{\alpha }\right)}+\epsilon N \ln \left(\frac{N}{\alpha }\right)+\epsilon N \sqrt{N}\left\lVert \nu \right\rVert_{2}^{2},\\
            &\beta _{2}=\epsilon N \sqrt{\epsilon Nk ^{\star } \ln \left(\frac{1}{\epsilon }\right)}+(\epsilon N)^{2}\ln \left(\frac{1}{\epsilon }\right)+\left\lVert \nu \right\rVert_{2}^{}(\epsilon N)^{2}\sqrt{\ln \left(\frac{1}{\epsilon }\right)}+\left\lVert \nu \right\rVert_{2}^{2}(\epsilon N)^{2},
        \end{aligned}
    \end{equation}
    given the condition that $N \ge \frac{\ln (1/\alpha )}{\epsilon \ln (1/\epsilon )}$.
\end{lemma}

\begin{proof}
    The last condition is of no importance since it only depends on the predetermined constants $\alpha $, $\epsilon $, and the fact that $\epsilon \ge \frac{1}{N}$. We prove this lemma by direct check of the definition.

    \textbf{Check of (i).} By adding and subtracting $\nu $, for a single sample $Y$, we have
    \begin{equation*}
        \left\lVert Y\right\rVert_{2}^{2}-k ^{\star }=\left\lVert Y-\nu \right\rVert_{2}^{2}-k ^{\star }+\left\lVert \nu \right\rVert_{2}^{2}+2\nu ^{\top } (Y-\nu ).
    \end{equation*}
    By assumption of the distribution of $Y$, $\left\lVert Y-\mu\right\rVert_{2}^{2}$ follows a chi-square distribution with degree of freedom $k ^{\star }$. By standard concentration bound, we have:
    \begin{equation*}
        \mathbb{P}\left(\left|\left\lVert Y-\nu \right\rVert_{2}^{2}-k ^{\star }\right|\ge t\right)\le 2 \exp\left\{-C \min\limits \left\{\frac{t ^{2}}{k ^{\star }},t\right\}\right\}.
    \end{equation*}
    Now set the probability to be $\frac{\alpha }{6N}$, we have $t=c \max\limits \left\{\sqrt{k ^{\star }\ln \left(\frac{N}{\alpha }\right)},\ln \left(\frac{N}{\alpha }\right)\right\}\le c \left(\sqrt{k ^{\star } \ln \left(\frac{N}{\alpha }\right)}+\ln \left(\frac{N}{\alpha }\right)\right)$ for some appropriate $c$.\\ 
    $\nu ^{\top } (Y-\nu )$ follows a Gaussian distribution with mean zero and variance $\left\lVert \nu \right\rVert_{2}^{2}$. By the concentration bound for the Gaussian, we have 
    \begin{equation*}
        \mathbb{P}\left(\left|2\nu ^{\top } (Y-\nu )\right|\ge t\right)\le 2 \exp\left\{-\frac{t ^{2}}{8 \left\lVert \nu \right\rVert_{2}^{2}}\right\}.
    \end{equation*}
    Setting the probability to be $\frac{\alpha }{6N}$ as the same, we find $t=c \sqrt{\left\lVert \nu \right\rVert_{2}^{2}\ln \left(\frac{N}{\alpha }\right)}$.
    Combining the results above, we have
    \begin{equation*}
        \left|\left\lVert Y _{i}\right\rVert_{2}^{2}-k ^{\star }\right|\le c \left(\left(\sqrt{k ^{\star }}+\left\lVert \nu \right\rVert_{2}^{}\right)\sqrt{\ln \left(\frac{N}{\alpha }\right)}+\ln \left(\frac{N}{\alpha }\right)+\left\lVert \nu \right\rVert_{2}^{2}\right)
    \end{equation*}
    with probability greater than $1-\frac{\alpha }{3N}$. Consequently, this holds for all $Y _{i},1 \le i \le N$ with probability greater than $1-\frac{\alpha }{3}$ by the union bound. Therefore, we have 
    \begin{equation}
        \left|\sum\limits_{i \in S}^{}\left\lVert Y _{i}\right\rVert_{2}^{2}-\left|S\right|k ^{\star }\right|\le c \left[\epsilon N \sqrt{(k ^{\star }+\left\lVert \nu \right\rVert_{2}^{2})\ln \left(\frac{N}{\alpha }\right)}+\epsilon N \ln \left(\frac{N}{\alpha }\right)+\epsilon N \left\lVert \nu \right\rVert_{2}^{2}\right]
        \label{regularity_con_i}
    \end{equation}

    \textbf{Check of (ii).} By the same trick, we have $\left\lVert \sum\limits_{i \in S}^{}Y _{i}\right\rVert_{2}^{2}=\left\lVert \sum\limits_{i \in S}^{}(Y _{i}-\nu )\right\rVert_{2}^{2}+\left\lVert \nu \right\rVert_{2}^{2}\left|S\right|^{2}+2 \left|S\right|\nu ^{\top } \sum\limits_{i \in S}^{}(Y _{i}-\nu )$. Similarly, we have $\left\lVert \sum\limits_{i \in S}^{}(Y _{i}-\nu )\right\rVert_{2}^{2}=\left|S\right|\left\lVert \frac{1}{\sqrt{\left|S\right|}}\sum\limits_{i \in S}^{}(Y _{i}-\nu )\right\rVert_{2}^{2}$, and 
    \begin{equation*}
        \mathbb{P}\left(\left|\left\lVert \frac{1}{\sqrt{\left|S\right|}}\sum\limits_{i \in S}^{}(Y _{i}-\nu )\right\rVert_{2}^{2}-k ^{\star }\right|\ge t\right)\le 2 \exp\left\{-c \min\limits \left\{\frac{t ^{2}}{k ^{\star }},t\right\}\right\}.
    \end{equation*}
    By Lemma \ref{lemma: combination number}, we have at most $\left(\frac{e}{\epsilon }\right)^{\epsilon N}$ different possible $S$ in total. Setting the probability to be $\frac{\alpha }{6}\left(\frac{\epsilon }{e}\right)^{\epsilon N}$ and solving for $t$, we obtain 
    \begin{equation*}
        \begin{aligned}
            t&=c \max\limits \left\{\sqrt{k ^{\star } \left[\ln \left(\frac{1}{\alpha }\right)+\epsilon N \ln \left(\frac{1}{\epsilon }\right)\right]},\ln \left(\frac{1}{\alpha }\right)+\epsilon N \ln \left(\frac{1}{\epsilon }\right)\right\}\\ 
            &\le c \left[\sqrt{k ^{\star } \left[\ln \left(\frac{1}{\alpha }\right)+\epsilon N \ln \left(\frac{1}{\epsilon }\right)\right]}+\ln \left(\frac{1}{\alpha }\right)+\epsilon N \ln \left(\frac{1}{\epsilon }\right)\right]\\ 
            &\le c \left[\sqrt{k ^{\star }\epsilon N \ln \left(\frac{1}{\epsilon }\right)}+\epsilon N \ln \left(\frac{1}{\epsilon }\right)\right],
        \end{aligned}
    \end{equation*}
    where the last step use the condition that $N \ge \frac{\ln (1/\alpha )}{\epsilon \ln (1/\epsilon )}$.

    For $\nu ^{\top } \sum\limits_{i \in S}^{}(Y _{i}-\nu )$, by the same concentration bound and technique, we will obtain
    \begin{equation*}
        \mathbb{P}\left(\left|\nu ^{\top } \sum\limits_{i \in S}^{}(Y _{i}-\nu )\right|\ge c \left\lVert \nu \right\rVert_{2}^{}\sqrt{\epsilon N}\cdot \sqrt{\epsilon N \ln \left(\frac{1}{\epsilon }\right)}\right)\le \frac{\alpha }{6}\left(\frac{\epsilon }{e}\right)^{\epsilon N}.
    \end{equation*} 
    
    Combining the results together, and using the union bound over all possible $S$, we have
    \begin{equation*}
        \left|\left\lVert \sum\limits_{i \in S}^{}Y _{i}\right\rVert_{2}^{2}-\left|S\right|k ^{\star }\right|\le c \left[\epsilon N \sqrt{\epsilon Nk ^{\star } \ln \left(\frac{1}{\epsilon }\right)}+(\epsilon N)^{2}\ln \left(\frac{1}{\epsilon }\right)+\left\lVert \nu \right\rVert_{2}^{}(\epsilon N)^{2}\sqrt{\ln \left(\frac{1}{\epsilon }\right)}+\left\lVert \nu \right\rVert_{2}^{2}(\epsilon N)^{2}\right]
    \end{equation*}
    with probability greater than $1-\frac{\alpha }{3}$. The RHS is exactly $\beta _{2}$ in \eqref{eq: specified beta_1 and beta_2}.

    \textbf{Check of (iii).} Finally, we are going to evaluate the term $\left\langle \sum\limits_{i \in S}^{}Y _{i}, \sum\limits_{j \in [N]}^{}Y _{j}\right\rangle$. Since this is equal to $\sum\limits_{i \in S}^{}\left\lVert Y _{i}\right\rVert_{2}^{2}+\sum\limits_{i \in S}^{}\left\langle Y _{i}, \sum\limits_{j \in [N] \backslash i}^{}Y _{j} \right\rangle$. We have obtained the bound for $\left\lVert Y _{i}\right\rVert_{2}^{2}$ in (i). For the term of inner product, we have 
    \begin{equation*}
        \left\langle Y _{i}, \sum\limits_{j \in [N] \backslash i}^{}Y _{j} \right\rangle=\left\langle Y _{i}-\nu , \sum\limits_{j \in [N] \backslash i}^{}(Y _{j}-\nu ) \right\rangle+\left\langle Y _{i}-\nu , (N-1)\nu  \right\rangle+\left\langle \nu , \sum\limits_{j \in [N] \backslash i}^{}(Y _{j}-\nu ) \right\rangle+(N-1)\left\lVert \nu \right\rVert_{2}^{2}.
    \end{equation*}
    Applying the sub-exponential concentration bound for the first term and sub-Gaussian concentration bound for the second as the previous handling, this is bounded by 
    \begin{equation*}
        c \left[\sqrt{Nk ^{\star } \ln \left(\frac{N}{\alpha }\right)}+\ln \left(\frac{N}{\alpha }\right)+(N-1)\left\lVert \nu \right\rVert_{2}^{}\sqrt{\ln \left(\frac{N}{\alpha }\right)}+(N-1)\left\lVert \nu \right\rVert_{2}^{2}\right].
    \end{equation*}
    With the bound of $\left\lVert Y _{i}\right\rVert_{2}^{2}$ in (i), the whole thing is bounded by 
    \begin{equation*}
        \begin{aligned}
            &c _{1}\left[\epsilon N \left(\sqrt{k ^{\star }}+\left\lVert \nu \right\rVert_{2}^{}\right)\sqrt{\ln \left(\frac{N}{\alpha }\right)}+\epsilon N \ln \left(\frac{N}{\alpha }\right)+\epsilon N \left\lVert \nu \right\rVert_{2}^{2}\right]\\ 
            +&c _{2}\left[\epsilon N \sqrt{Nk ^{\star } \ln \left(\frac{N}{\alpha }\right)}+\epsilon N \ln \left(\frac{N}{\alpha }\right)+\epsilon N ^{2}\left\lVert \nu \right\rVert_{2}^{}\sqrt{\ln \left(\frac{N}{\alpha }\right)}+\epsilon N ^{2}\left\lVert \nu \right\rVert_{2}^{2}\right]\\ 
            \le &c \left[\epsilon N \left(\sqrt{Nk ^{\star }}+N \left\lVert \nu \right\rVert_{2}^{}\right)\sqrt{\ln \left(\frac{N}{\alpha }\right)}+\epsilon N \ln \left(\frac{N}{\alpha }\right)+\epsilon N ^{2}\left\lVert \nu \right\rVert_{2}^{2}\right]
        \end{aligned}
    \end{equation*}
    Comparing with \eqref{regularity_con_i} and Definition \ref{def_regularity}, we find that the required $\beta _{1}$ is
    \begin{equation*}
        \epsilon N \left(\sqrt{k ^{\star }}+\sqrt{N}\left\lVert \nu \right\rVert_{2}^{}\right)\sqrt{\ln \left(\frac{N}{\alpha }\right)}+\epsilon N \ln \left(\frac{N}{\alpha }\right)+\epsilon N \sqrt{N}\left\lVert \nu \right\rVert_{2}^{2},
    \end{equation*}
    which is exactly $\beta _{1}$ in \eqref{eq: specified beta_1 and beta_2}. The proof is completed.
\end{proof}

By Lemma \ref{lemma: good_samples_regularity}, we should be able to control the inner product between the sum of samples on two sets $S$ and $S ^{\prime }$ by the polarization identity. The next corollary exactly states this benefit.

\begin{corollary}\label{corollary: generalized result of regularity 1}
    Assume $\mathbf{Y}$ is $(2\epsilon ,\beta _{1},\beta _{2})$-regular. Then for any subsets $S,S ^{\prime }\subset [N]$ with $\left|S\right|\le \epsilon N, \left|S ^{\prime }\right|\le \epsilon N$, we have with probability greater than $1-\alpha $:
    \begin{equation*}
        \left|\left\langle \sum\limits_{i \in S}^{}Y _{i}, \sum\limits_{j \in S ^{\prime }}^{}Y _{j} \right\rangle-\left|S \cap S ^{\prime }\right|k ^{\star }\right|\le c \beta _{2}.
    \end{equation*}
\end{corollary}

\begin{proof}
    By the polarization identity, we know 
    \begin{equation*}
        \left\langle \sum\limits_{i \in S}^{}Y _{i}, \sum\limits_{j \in S ^{\prime }}^{}Y _{j} \right\rangle=\frac{1}{2}\left[\left\lVert \sum\limits_{i \in S \cup S ^{\prime }}^{}Y _{i}\right\rVert_{2}^{2}+\left\lVert \sum\limits_{i \in S \cap S ^{\prime }}^{}Y _{i}\right\rVert_{2}^{2}-\left\lVert \sum\limits_{i \in S \backslash S ^{\prime }}^{}Y _{i}\right\rVert_{2}^{2}-\left\lVert \sum\limits_{i \in S ^{\prime }\backslash S}^{}Y _{i}\right\rVert_{2}^{2}\right].
    \end{equation*}
    Since each of $S \cup S ^{\prime }, S \cap S ^{\prime },S \backslash S ^{\prime },S ^{\prime }\backslash S$ has cardinality no more than $2\epsilon N$ by assumption, by Definition \ref{def_regularity} and Lemma \ref{lemma: good_samples_regularity}, we know 
    \begin{equation*}
        \left|\left\langle \sum\limits_{i \in S}^{}Y _{i}, \sum\limits_{j \in S ^{\prime }}^{}Y _{j} \right\rangle-\frac{1}{2} \left(\left|S \cup S ^{\prime }\right|+\left|S \cap S ^{\prime }\right|-\left|S \backslash S ^{\prime }\right|-\left|S ^{\prime }\backslash S\right|\right)k ^{\star }\right|\le c \beta _{2}.
    \end{equation*}
    Finally, since $\left|S \cup S ^{\prime }\right|-\left|S \backslash S ^{\prime }\right|-\left|S ^{\prime }\backslash S\right|=\left|S \cap S ^{\prime }\right|$, the proof is completed.
\end{proof}

With Lemma \ref{lemma: convex combination}, it is not hard to extend the result in Corollary \ref{corollary: generalized result of regularity 1} to the case dealing with the convex combination of the indicator vectors.

\begin{corollary}[]\label{corollary: generalized result of regularity 2}
    Assume $\mathbf{Y}$ is $(2\epsilon ,\beta _{1},\beta _{2})$-regular, then for any vectors $a,b \in [0,1]^{N}$ with $\left\lVert a\right\rVert_{1}^{}\le \epsilon N, \left\lVert b\right\rVert_{1}^{}\le \epsilon N$, we have with probability greater than $1-\alpha $:\\ 
    \begin{tabular}{rl}
        (i), & $\left|\left\langle \sum\limits_{i=1}^{N}a _{i}Y _{i}, \sum\limits_{j=1}^{N}b _{j}Y _{j}\right\rangle-\left\langle a, b \right\rangle k ^{\star }\right|\le c _{1}\beta _{2}$;\\ 
        (ii), & $\left|\left\langle \sum\limits_{i=1}^{N}a _{i}Y _{i}, \sum\limits_{j=1}^{N}Y _{j} \right\rangle-\left\lVert a\right\rVert_{1}^{}k ^{\star }\right|\le c _{2}\sqrt{N}\beta _{1}$
    \end{tabular}
\end{corollary}
In Corollary \ref{corollary: generalized result of regularity 1}, there is no corresponding (ii) since that is exactly the third definition of regularity.

\begin{proof}
    We proof by direct calculation. By Lemma \ref{lemma: convex combination}, the vectors $a$ and $b$ can be represent as the convex combination of some indicator vectors, respectively. Let $a=\sum\limits_{k=1}^{M _{a}}\xi _{k}\mathbf{1}_{S _{k}}, b=\sum\limits_{l=1}^{M _{b}}\eta _{l}\mathbf{1}_{S _{l}}$, where $S _{k},S _{l}\subset [N],1 \le k \le N,1 \le l \le N$. We have:
    \begin{equation*}
        \begin{aligned}
            \left\langle \sum\limits_{i=1}^{N}a _{i}Y _{i}, \sum\limits_{j=1}^{N}b _{j}Y _{j} \right\rangle&=\left\langle \sum\limits_{i=1}^{N}\sum\limits_{k=1}^{M _{a}}\xi _{k}\mathbf{1}_{S _{k},i}Y _{i}, \sum\limits_{j=1}^{N}\sum\limits_{l=1}^{M _{b}}\eta _{l}\mathbf{1}_{S _{l},j} Y _{j}\right\rangle\\  
            &=\sum\limits_{k=1}^{M _{a}}\sum\limits_{l=1}^{M _{b}}\xi _{k}\eta _{l}\left\langle \sum\limits_{i \in S _{k}}^{}Y _{i}, \sum\limits_{j \in S _{l}}^{}Y _{j} \right\rangle\\ 
            &\overset{(1)}{=}\sum\limits_{k=1}^{M _{a}}\sum\limits_{l=1}^{M _{b}}\xi _{k}\eta _{l}\left(\left|S _{k}\cap S _{l}\right|k ^{\star }+R _{k,l}\right)\\ 
            &=k ^{\star } \cdot \sum\limits_{k=1}^{M _{a}}\sum\limits_{l=1}^{M _{b}}\xi _{k}\eta _{l}\sum\limits_{i=1}^{N}\mathbf{1}_{S _{k},i}\cdot \mathbf{1}_{S _{l},i}+\sum\limits_{k=1}^{M _{a}}\sum\limits_{l=1}^{M _{b}}\xi _{k}\eta _{l}R _{k,l}\\ 
            &=k ^{\star } \cdot \sum\limits_{i=1}^{N}\left(\sum\limits_{k=1}^{M _{a}}\xi _{k}\mathbf{1}_{S _{k},i}\right)\left(\sum\limits_{l=1}^{M _{b}}\eta _{l}\mathbf{1}_{S _{l},i}\right)+\sum\limits_{k=1}^{M _{a}}\sum\limits_{l=1}^{M _{b}}\xi _{k}\eta _{l}R _{k,l}\\ 
            &=k ^{\star } \cdot \sum\limits_{i=1}^{N}a _{i}b _{i}+\sum\limits_{k=1}^{M _{a}}\sum\limits_{l=1}^{M _{b}}\xi _{k}\eta _{l}R _{k,l}.
        \end{aligned}
    \end{equation*}
    Here in the equation, (1) is from Corollary \ref{corollary: generalized result of regularity 1}, and $R _{k,l}\le c \beta _{2}$ for some constant $c$. Now that by Lemma \ref{lemma: convex combination}, $0 \le \sum\limits_{k=1}^{M _{a}}\xi _{k}\le 1,0 \le \sum\limits_{l=1}^{M _{b}}\eta _{l}\le 1$, we know $\sum\limits_{k=1}^{M _{a}}\sum\limits_{l=1}^{M _{b}}\xi _{k}\eta _{l}R _{k,l}\le c \beta _{2}$. Therefore, (i) is checked. With very similar argument, we can also check (ii). The proof is completed.
\end{proof}

At this stage, we have all tools needed for the proof. The goal of the following Algorithms \ref{algorithm: prefiltering}, \ref{algorithm n > k}, \ref{algorithm n <= k} is to filter the weight $\omega $ and the samples $\mathbf{Y}$ such that the resulting $\omega $ and $\mathbf{Y}$ meet the requirements of \hyperref[def_omega_regularity]{$\omega $-regularity}. We start with the first requirement.

\begin{algorithm*}[htbp]
    \caption{Prefiltering.}\label{algorithm: prefiltering}
    Set $\gamma _{1}=c\left[\sqrt{k ^{\star }\ln \left(\frac{N}{\alpha }\right)}+\ln \left(\frac{N}{\alpha }\right)\right]$, $\texttt{count}=0,\texttt{i}=0$

    \While {$\texttt{i} < N$}{
        \If {$\left|\left\lVert Y _{i}\right\rVert_{2}^{2}-k ^{\star }\right|> \gamma _{1}$}{

            $\texttt{count}=\texttt{count}+1$

            \If {$\texttt{count}>\epsilon N$}{
                \Return \texttt{None}
            } 
            Delete $Y _{i}$ from $\mathbf{Y}$      
        }
        $\texttt{i}=\texttt{i}+1$
    }

    \Return $\mathbf{Y}$
\end{algorithm*}

\begin{lemma}[]\label{lemma: guarantee of omega regularity 1}
    Upon finishing Algorithm \ref{algorithm: prefiltering}, all of samples satisfy
    \begin{equation*}
        \left|\left\lVert Y _{i}\right\rVert_{2}^{2}-k ^{\star }\right| \le \gamma _{1}.
    \end{equation*}
    Moreover, under the null hypothesis $H _{0}$, the algorithm will terminate without rejecting with probability at least $1-\alpha $.
\end{lemma}

\begin{proof}
    The satisfaction of the inequality comes directly from the definition of Algorithm \ref{algorithm: prefiltering}. According to Lemma \ref{lemma: good_samples_regularity} where setting $\left\lVert \nu \right\rVert_{2}^{}=0$, all uncontaminated samples will satisfy the inequality with probability at least $1-\alpha $ and therefore will not be deleted. The remaining follows.
\end{proof}

\begin{remark*}
    Rigorously speaking, after the prefiltering, the number of samples left should be adjusted to $N-\texttt{count}$, where $\texttt{count}$ is possible to be positive. However, if the prefiltering finishes without rejection, \texttt{count} is at most $\epsilon N$, upon which, and for the conciseness of the notations, we continue to use $N$ in the proofs later. One can verify that this does not affect the correctness of the following proofs.
\end{remark*}

The next two algorithms both deal with the second requirement of \hyperref[def_omega_regularity]{$\omega $-regularity}. However, different techniques are required for the classic scenario $N>k ^{\star }$ and the high-dimensional scenario $N \le k ^{\star }$. Assume first $N>k ^{\star }$.

Let
\begin{equation}\label{eq: definition of gamma_2}
    \gamma _{2}:=c \left[\sqrt{Nk ^{\star }}+\sqrt{N\ln \left(\frac{1}{\alpha }\right)}+\ln \left(\frac{1}{\alpha }\right)+\epsilon N \ln \left(\frac{1}{\epsilon }\right)\right].
\end{equation}

Denote $D(\omega )=k ^{\star } \cdot \text{diag}\{\omega \}$, and $\lambda (\omega )=\left\lVert Y ^{\top } D(\omega )Y-N \mathbf{I}_{k ^{\star }}\right\rVert_{2}^{}$. The major part of the algorithm is a ``while'' iteration, whose control condition is $\lambda (\omega ^{(t)})\ge \gamma _{2}$, where $\omega (t)$ represents the value of $\omega $ in the $t$-th iteration. Inside each iteration, let $v ^{(t)}$ denote the singular vector associated with $\lambda (\omega ^{(t)})$. We first compute $\tau _{i}^{(t)}=\left\langle v, Y _{i} \right\rangle ^{2}\mathbf{1}_{\left\{\omega _{i}^{(t)}>0\right\}}$, then sort the samples according to the decreasing order of $\tau _{i}^{(t)}$. (Without any loss of generality, we assume $\tau _{i}^{(t)}$ are already sorted.) Next, let $I$ be the smallest index such that $\sum\limits_{i=1}^{I}\omega _{i}^{(t)}\ge 2\epsilon N$ and update $\omega _{i}^{(t+1)}$ as following:
\begin{equation}\label{eq: update of omega}
    \omega _{i}^{(t+1)}=\left\{
    \begin{tabular}{ll}
        $\left(1-\frac{\tau _{i}^{(t)}}{\tau _{1}^{(t)}}\right)\omega _{i}^{(t)}$ & if $i \le I$,\\
        $\omega _{i}^{(t)}$& if $i>I$.
    \end{tabular}
    \right.
\end{equation}
Finally, we check the condition whether $\left\lVert \omega \right\rVert_{1}^{}\le N(1-2\epsilon )$ and reject the null hypothesis directly if this holds and pass to the following step if not.

\begin{algorithm}[htbp]
    \caption{Sample filtering when $N>k ^{\star }$.}\label{algorithm n > k}
    Set $\gamma _{2}$ as \eqref{eq: definition of gamma_2}, and $\lambda =\left\lVert Y ^{\top } D(\omega )Y-N \mathbf{I}_{k ^{\star }}\right\rVert_{2}^{}$. ($\omega $ is initialized as $\mathbf{1}$.)

    \While {$\lambda \ge \gamma _{2}$}{
        Set $v$ to be the unit singular vector associated with $\lambda $

        Compute $\tau _{i}=\left\langle v, Y _{i} \right\rangle ^{2}\mathbf{1}_{\left\{\omega _{i}>0\right\}}$ for $1 \le i \le N$

        Set $I$ be the smallest index such that $\sum\limits_{i=1}^{I}\omega _{i}\ge 2\epsilon N$

        Update $w$ according to \eqref{eq: update of omega}

        \If {$\left\lVert \omega \right\rVert_{1}^{}<N(1-2\epsilon )$} {

            \Return \texttt{None}
        }

        Set $\lambda =\left\lVert Y ^{\top } D(\omega )Y-N \mathbf{I}_{k ^{\star }}\right\rVert_{2}^{}$
    }

    \Return $\omega $
\end{algorithm}

We are going to prove that if the algorithm finishes without rejecting the null, the second condition of \hyperref[def_omega_regularity]{$\omega $-regularity} is guaranteed with high probability. This is formalized as the following lemma.

\begin{lemma}[]\label{lemma: guarantee of omega regularity 2 1}
    Given the algorithm \ref{algorithm n > k}, we have the following facts:\\
    \hspace*{0.5em}(i), the algorithm will terminate within finite rounds,\\ 
    \hspace*{0.5em}(ii), under $H _{0}$, the algorithm will not output rejection with high probability;\\
    \hspace*{0.5em}(iii), once finished without rejection, for any $S \subset [N]$ with $\left|S\right|\le \epsilon N$, we have $\left\lVert \sum\limits_{i \in S}^{}\sqrt{\omega _{i}}Y _{i}\right\rVert_{2}^{2}\le 2 \gamma _{2}\epsilon N$ with probability higher than $1-\alpha $.
\end{lemma}

\begin{proof}
Proof of (i) is relatively easy. By the update \eqref{eq: update of omega}, the $\omega _{i}$ corresponding to the maximal $\tau _{i}$ will be set to zero in each iteration. Consequently, $\left\lVert \omega \right\rVert_{1}^{}$ will be less than $N(1-2\epsilon )$ in at most $2\epsilon N$ rounds, followed by the rejection.

For (ii), define $\Lambda _{1}=\left\{\omega :\omega \in [0,1]^{N},\left\lVert \mathbf{1}_{[N] \backslash C}-\omega _{[N] \backslash C}\right\rVert_{1}^{}\le \left\lVert \mathbf{1}_{C}-\omega _{C}\right\rVert_{1}^{}\right\}$, and denote $\omega ^{(t)}, \tau ^{(t)}$ as the value of $\omega ,\tau $ in the $t$-th iteration. We attempt to prove that under $H _{0}$, if $\omega ^{(t)}\in \Lambda _{1}$ and we are still in the iteration, then so does $\omega ^{(t+1)}$. Suppose this holds, we have $\left\lVert \mathbf{1}-\omega ^{(t+1)}\right\rVert_{1}^{}\le 2 \left\lVert \mathbf{1}_{C}-\omega _{C} ^{(t+1)}\right\rVert_{1}^{}\le 2\epsilon N$, therefore the condition for rejection will never be triggered.

Now it remains to prove the statement. Suppose $\omega ^{(t)}\in \Lambda _{1}$ and we are still in the iteration. By the definition of $\Lambda _{1}$ and the update \eqref{eq: update of omega}, it suffices to prove $\sum\limits_{i \in [I]\cap  C}^{}\tau _{i}^{(t)}\omega _{i}^{(t)}\ge \sum\limits_{i \in [I]\cap ([N] \backslash C)}^{}\tau _{i}^{(t)}\omega ^{(t)}$, i.e.,
\begin{equation*}
    v ^{\top } \left(\sum\limits_{i \in [I]\cap C}^{}\omega _{i}^{(t)}Y _{i}Y _{i}^{\top }\right)v \ge v ^{\top } \left(\sum\limits_{i \in [I]\cap ([N] \backslash C)}^{}\omega _{i}^{(t)}Y _{i}Y _{i}^{\top }\right)v.
\end{equation*}
By the control condition, we have
\begin{equation*}
    v ^{\top } \left(\sum\limits_{i=1}^{N}\omega _{i}^{(t)}Y _{i}Y _{i}^{\top } -N \mathbf{I}_{k ^{\star }}\right)v \ge \gamma _{2}.
\end{equation*}
By Corollary \ref{corollary: bounds on the l2-operator norm of X^TX and XX^T} (set $\left\lVert \nu \right\rVert_{2}^{}=0$) and the definition of $\gamma _{2}$ \eqref{eq: definition of gamma_2}, with high probability, we have 
\begin{equation*}
    \begin{aligned}
        v ^{\top } \left(\sum\limits_{i \in C}^{}\omega _{i}^{(t)}Y _{i}Y _{i}^{\top } \right)v&=v ^{\top } \left(\sum\limits_{i=1}^{N}\omega _{i}^{(t)}Y _{i}Y _{i}^{\top } -N \mathbf{I}_{k ^{\star }}\right)v-v ^{\top } \left(\sum\limits_{i \in [N] \backslash C}^{}\omega _{i}^{(t)}Y _{i}Y _{i}^{\top } -N \mathbf{I}_{k ^{\star }}\right)v\\ 
        &\ge v ^{\top } \left(\sum\limits_{i=1}^{N}\omega _{i}^{(t)}Y _{i}Y _{i}^{\top } -N \mathbf{I}_{k ^{\star }}\right)v-v ^{\top } \left(\sum\limits_{i=1}^{N}\tilde{Y}_{i}\tilde{Y}_{i}^{\top } -N \mathbf{I}_{k ^{\star }}\right)v\\
        &\ge (1-c)\gamma _{2}.
    \end{aligned}
\end{equation*}
By the definition of $I$ we know $2\epsilon N \le \sum\limits_{i \in [I]}^{}\omega _{i}<2\epsilon N+1$. Because $\omega ^{(t)}\in \Lambda _{1}$, we know $\left\lVert \mathbf{1}-\omega ^{(t)}\right\rVert_{1}^{}\le 2 \left\lVert \mathbf{1}_{C}-\omega _{C}^{(t)}\right\rVert_{1}^{}\le 2\epsilon N$. Hence $I-(2\epsilon N+1)\le \sum\limits_{i \in [I]}^{}(1-\omega _{i}^{(t)})\le 2\epsilon N$, which leads to $I \le 4\epsilon N+1$. Therefore, according to Lemma \ref{lemma: bound on the max eigenvalue of the covariance matrix}, with high probability, we have 
\begin{equation}\label{eq: proof of the lemma 6 1}
    v ^{\top } \left(\sum\limits_{i \in [I]\cap ([N] \backslash C)}^{}\omega _{i}^{(t)}Y _{i}Y _{i}^{\top }\right)v \le c \gamma _{2}.
\end{equation}
Now that we sorted $\tau _{i}$ by decreasing order, $\sum\limits_{i \in [I]}^{}\omega _{i}\ge 2\epsilon N$, and $\left|C\right|\le \epsilon N$, we have
\begin{equation*}
    \sum\limits_{i \in C}^{}\omega _{i}^{(t)}\tau _{i}^{(t)}-\sum\limits_{i \in [I]}^{}\omega _{i}^{(t)}\tau _{i}^{(t)}\le \sum\limits_{i \in C \backslash [I]}^{}\omega _{i}^{(t)}\tau _{i}^{(t)}\le \frac{1}{2}\sum\limits_{i \in [I]}^{}\omega _{i}^{(t)}\tau _{i}^{(t)}.
\end{equation*}
Therefore $\sum\limits_{i \in [I]}^{}\omega _{i}^{(t)}\tau _{i}^{(t)}\ge \frac{2}{3}\sum\limits_{i \in C}^{}\omega _{i}^{(t)}\tau _{i}^{(t)}\ge \frac{2(1-c)}{3}\gamma _{2}$. Note that by selecting the constant in \eqref{eq: definition of gamma_2} sufficiently large, we can let $c$ sufficiently small. Hence, we can assume $c \le \frac{1}{4}$, by which we have $\sum\limits_{i \in [I]}^{}\omega _{i}^{(t)}\tau _{i}^{(t)}\ge \frac{1}{2}\gamma _{2}$. Combine with \eqref{eq: proof of the lemma 6 1}, we find 
\begin{equation*}
    \sum\limits_{i \in [I]\cap C}^{}\omega _{i}^{(t)}\tau _{i}^{(t)}\ge \frac{1}{4}\gamma _{2}\ge \sum\limits_{i \in [I]\cap ([N] \backslash C)}^{}\omega _{i}^{(t)}\tau _{i}^{(t)}.
\end{equation*}

Hence, (ii) is proved.

To prove (iii), we need the following lemma, which exploits the property of the norm of the mean from the algorithm.

\begin{lemma}[]\label{lemma: property on the mean}
    If $\left\lVert \nu \right\rVert_{2}^{}\ge c$ for some sufficient large but universal constant $c$, then Algorithm \ref{algorithm n > k} will reject $H _{0}$ with high probability.
\end{lemma}

\begin{proof}

    First, we fix the vector $v=\frac{\nu }{\left\lVert \nu \right\rVert_{2}^{}}$. If the algorithm exits without rejection, we have 
    \begin{equation*}
        \left\lVert \sum\limits_{i=1}^{N}\omega _{i}Y _{i}Y _{i}^{\top } -N \mathbf{I}_{k ^{\star }}\right\rVert_{2}^{}<\gamma _{2}.
    \end{equation*}
    This suggests that $\sum\limits_{i=1}^{N}\omega _{i}(v ^{\top } Y _{i})^{2}-N < \gamma _{2}$, based on which we have 
    \begin{equation*}
        \frac{\sum\limits_{i \in [N] \backslash C}^{}\omega _{i}(v ^{\top } Y _{i})^{2}}{N}-1 \le \frac{\sum\limits_{i=1}^{N}\omega _{i}(v ^{\top } Y _{i})^{2}}{N}-1 < \frac{\gamma _{2}}{N}\le c,
    \end{equation*}
    where the last step is from the definition of $\gamma _{2}$ \eqref{eq: definition of gamma_2}. If we expand $v ^{\top } Y _{i}$, this merely says 
    \begin{equation}\label{eq: proof of the lemma 6 2}
        \sum\limits_{i \in [N] \backslash C}^{}\omega _{i}\left[v ^{\top } (Y _{i}-\nu )\right]^{2}+\sum\limits_{i \in [N] \backslash C}^{}\omega _{i}\left\lVert \nu \right\rVert_{2}^{2}+2 \sum\limits_{i \in [N] \backslash C}^{}\omega _{i}\left[v ^{\top } (Y _{i}-\nu )\right]\left\lVert \nu \right\rVert_{2}^{}\le cN.
    \end{equation}
    By the fact that $ab \ge -\frac{1}{4}a ^{2}-b ^{2}$ for $a,b \in \mathbb{R}$, we know 
    \begin{equation*}
        2 \sum\limits_{i \in [N] \backslash C}^{}\omega _{i}\left[v ^{\top } (Y _{i}-\nu )\right]\left\lVert \nu \right\rVert_{2}^{}\ge -\frac{1}{2}\sum\limits_{i \in [N] \backslash C}^{}\omega _{i}\left\lVert \nu \right\rVert_{2}^{2}-2 \sum\limits_{i \in [N] \backslash C}^{}\omega _{i}\left[v ^{\top } (Y _{i}-\nu )\right]^{2}.
    \end{equation*}
    Adding both sides with $\sum\limits_{i \in [N] \backslash C}^{}\omega _{i}\left\lVert \nu \right\rVert_{2}^{2}$ and plugging \eqref{eq: proof of the lemma 6 2} in, we have 
    \begin{equation*}
        \frac{1}{2}\sum\limits_{i \in [N] \backslash C}^{}\omega _{i}\left\lVert \nu \right\rVert_{2}^{2}\le c N+\sum\limits_{i \in [N] \backslash C}^{}\omega _{i}\left[v ^{\top } (Y _{i}-\nu )\right]^{2}\le c N+\sum\limits_{i \in [N] \backslash C}^{}\left[v ^{\top } (Y _{i}-\nu )\right]^{2}.
    \end{equation*}
    Now since $\left\lVert v\right\rVert_{2}^{}=1$, we know that $\left[v ^{\top } (Y _{i}-\nu )\right]^{2}$ is subject to a chi-square distribution with degree of freedom 1. By the standard argument for the sub-exponential random variables (see \ref{def_sub_exp_var}), we know that $\sum\limits_{i \in [N] \backslash C}^{}\left[v ^{\top } (Y _{i}-\nu )\right]^{2}\le N+ c \sqrt{N}$ with (high) constant probability. On the other hand, by the definition of the algorithm we have $\left\lVert \omega _{[N] \backslash C}\right\rVert_{1}^{}\ge \left\lVert \omega \right\rVert_{1}^{}-\epsilon N \ge (1-3\epsilon )N$. Therefore, we have with (high) constant probability,
    \begin{equation*}
        \frac{1}{2}(1-3\epsilon )N \left\lVert \nu \right\rVert_{2}^{2}\le c _{1} \cdot N+c _{2} \cdot \sqrt{N}.
    \end{equation*}
    Now deviding both sides by $\frac{1}{2}(1-3\epsilon )N$, we see that Lemma \ref{lemma: property on the mean} holds if we take $c=\frac{2(c _{1}+c _{2})}{1-3\epsilon }$
\end{proof}

We are going to prove (iii) by contradiction. Suppose that at the end of the algorithm without rejection, we have $\left\lVert \sum\limits_{i \in S}^{}\sqrt{\omega _{i}}Y _{i}\right\rVert_{2}^{2}> 2 \gamma _{2}\epsilon N$ for some $S \subset [N]$ with $\left|S\right|\le \epsilon N$. This indicates that there exists a vector $v \in \mathbb{R}^{k ^{\star }}$ with $\left\lVert v\right\rVert_{2}^{}=1$ such that $\sum\limits_{i \in S}^{}\sqrt{\omega _{i}}v ^{\top } Y _{i}>\sqrt{2 \gamma _{2}\epsilon N}$. Therefore, by the Cauchy--Schwarz inequality, we have $\sum\limits_{i \in S \cup C}^{}\omega _{i}(v ^{\top } Y _{i})^{2}\ge 2 \gamma _{2}$. This implies that 
\begin{equation}\label{eq: proof of the lemma 5 1}
    \left\lVert \sum\limits_{i \in S \cup C}^{}\omega _{i}Y _{i}Y _{i}^{\top } +\left|[N] \backslash (S \cup C)\right|\nu \nu ^{\top } \right\rVert_{2}^{}\ge \left\lVert \sum\limits_{i \in S \cup C}^{}\omega _{i}Y _{i}Y _{i}^{\top } \right\rVert_{2}^{}\ge 2 \gamma _{2}.
\end{equation}
Now consider the quantity $\left\lVert \sum\limits_{i \in [N] \backslash (S \cup C)}^{}\omega _{i}Y _{i}Y _{i}^{\top } -\left|[N] \backslash (S \cup C)\right|\nu \nu ^{\top } - N \mathbf{I}_{k ^{\star }}\right\rVert_{}^{}$, by the triangular inequality, we have
\begin{equation*}
    \begin{aligned}
        \left\lVert \sum\limits_{i \in [N] \backslash (S \cup C)}^{}\omega _{i}Y _{i}Y _{i}^{\top } -\left|[N] \backslash (S \cup C)\right|\nu \nu ^{\top } - N \mathbf{I}_{k ^{\star }}\right\rVert_{}^{}&\le \left\lVert \sum\limits_{i=1}^{N}\omega _{i}\tilde{Y}_{i}\tilde{Y}_{i}^{\top } -N\nu \nu ^{\top }-N \mathbf{I}_{k ^{\star }}\right\rVert_{2}^{}\\ 
        &+\left\lVert \sum\limits_{i \in S \cup C}^{}\omega _{i}\tilde{Y}_{i}\tilde{Y}_{i}^{\top }-\left|S \cup C\right|\nu \nu ^{\top } \right\rVert_{2}^{}\\
        &\le \text{I}+\text{II}+\text{III}+ \left|S \cup C\right|\left\lVert \nu \right\rVert_{2}^{2}.
    \end{aligned}
\end{equation*}
where I, II and III are defined as 
\begin{equation*}
    \text{I}=\left\lVert \sum\limits_{i=1}^{N}\tilde{Y}_{i}\tilde{Y}_{i}^{\top } -N\nu \nu ^{\top } - N \mathbf{I}_{k ^{\star }}\right\rVert_{2}^{},\text{II}=\left\lVert \sum\limits_{i=1}^{N}(1-\omega _{i})\tilde{Y}_{i}\tilde{Y}_{i}^{\top } \right\rVert_{2}^{},\text{III}=\left\lVert \sum\limits_{i \in S \cup C}^{}\omega _{i}\tilde{Y _{i}}\tilde{Y _{i}}^{\top }\right\rVert_{}^{}.
\end{equation*}
By Lemma \ref{lemma: property on the mean}, without loss of generality we can assume $\left\lVert \nu \right\rVert_{2}^{}\le c _{0}$ for some constant $c _{0}$. For the last term, since $\left|S \cup C\right|\le 2\epsilon N$, obviously it can be bounded by $c\gamma _{2}$. We are going to deal with the other three parts respectively. 

\textbf{Bound on I}. We have: 
\begin{equation*}
    \text{I}\le \left\lVert \sum\limits_{i=1}^{N}(\tilde{Y}_{i}-\nu )(\tilde{Y}_{i}-\nu )^{\top } -N \mathbf{I}_{k ^{\star }}\right\rVert_{2}^{}+2 \left\lVert \sum\limits_{i=1}^{N}(\tilde{Y}_{i}-\nu )^{\top } \nu \right\rVert_{2}^{}.
\end{equation*}
Set $\left\lVert \nu \right\rVert_{2}^{}=0$ in Corollary \ref{corollary: bounds on the l2-operator norm of X^TX and XX^T}, we have 
\begin{equation}\label{eq: proof of the lemma 5 2}
    \left\lVert \sum\limits_{i=1}^{N}(\tilde{Y}_{i}-\nu )(\tilde{Y}_{i}-\nu )^{\top } -N \mathbf{I}_{k ^{\star }}\right\rVert_{2}^{}\le C _{1}\left(\sqrt{Nk ^{\star }}+\sqrt{N \ln \left(\frac{1}{\alpha }\right)}+\ln \left(\frac{1}{\alpha }\right)\right).
\end{equation}
By normalization and standard concentration tools (see \ref{def_sub_Gaussian_var}), we have 
\begin{equation}\label{eq: proof of the lemma 5 3}
    \left\lVert \sum\limits_{i=1}^{N}(\tilde{Y}_{i}-\nu )^{\top } \nu \right\rVert_{2}^{}=\left\lVert \nu \right\rVert_{2}^{}\cdot \left\lVert \sum\limits_{i=1}^{N}(\tilde{Y}_{i}-\nu )\right\rVert_{2}^{}=C _{0}\cdot \sqrt{Nk ^{\star }}\cdot \left\lVert \frac{1}{\sqrt{N}}\sum\limits_{i=1}^{N}\frac{\tilde{Y}_{i}-\nu }{\sqrt{k ^{\star }}}\right\rVert_{2}^{}\le c _{2}\sqrt{Nk ^{\star }}
\end{equation}
with (high) constant probability. Now combine \eqref{eq: proof of the lemma 5 2} and \eqref{eq: proof of the lemma 5 3} to get 
\begin{equation*}
    \text{I}\le (c _{1}+2c _{2})\left(\sqrt{Nk ^{\star }}+\sqrt{N \ln \left(\frac{1}{\alpha }\right)}+\ln \left(\frac{1}{\alpha }\right)\right)\le c \gamma _{2}.
\end{equation*}

\textbf{Bound on II}. Since upon the finish the algorithm, we have $\left\lVert \mathbf{1}-\omega \right\rVert_{1}^{}\le 2\epsilon N$, by Corollary \ref{corollary: bound on the max eigenvalue of the covariance matrix in general}, we know 
\begin{equation*}
    \text{II}\le c _{3}\left[\epsilon N \ln \left(\frac{1}{\epsilon }\right)+k ^{\star }+\ln \left(\frac{1}{\alpha }\right)\right]\le c \gamma _{2}.
\end{equation*}

\textbf{Bound on III}. Since $\left|S \cup C\right|\le 2\epsilon N$, by Corollary \ref{corollary: bound on the max eigenvalue of the covariance matrix in general} and an argument similar to the one in II, we also have 
\begin{equation*}
    \text{III}\le c \gamma _{2}.
\end{equation*}
Therefore, we conclude that 
\begin{equation}\label{eq: proof of the lemma 5 4}
    \left\lVert \sum\limits_{i \in [N] \backslash (S \cup C)}^{}\omega _{i}Y _{i}^{\top } Y _{i}-\left|[N] \backslash (S \cup C)\right|\nu \nu ^{\top } -N \mathbf{I}_{k ^{\star }}\right\rVert_{2}^{}\le 4c \gamma _{2}.
\end{equation}
Now combine \eqref{eq: proof of the lemma 5 1} and \eqref{eq: proof of the lemma 5 4} and use the triangular inequality, we find that 
\begin{equation*}
    \left\lVert \sum\limits_{i=1}^{N}\omega _{i}Y _{i}^{\top } Y _{i}-N \mathbf{I}_{k ^{\star }}\right\rVert_{2}^{}\ge (2-4c)\gamma _{2}.
\end{equation*}
Again, by choosing the constant in the definition of $\gamma _{2}$ \eqref{eq: definition of gamma_2} sufficiently large, we can make $c$ above sufficiently small. Therefore, we can assume $c<\frac{1}{4}$, which leads to $\left\lVert \sum\limits_{i=1}^{N}\omega _{i}Y _{i}^{\top } Y _{i}-N \mathbf{I}_{k ^{\star }}\right\rVert_{2}^{}>\gamma _{2}$, a contradiction of the control condition of Algorithm \ref{algorithm n > k}. The proof is completed.
\end{proof}

We now turn to the high-dimensional scenario where $N \le k ^{\star }$. In this case, we let
\begin{equation}\label{eq: definition of gamma_3}
    \gamma _{3}=c \left(\sqrt{Nk ^{\star }}+\sqrt{k ^{\star } \ln \left(\frac{1}{\alpha }\right)}+\ln \left(\frac{1}{\alpha }\right)+\epsilon N \ln \left(\frac{1}{\epsilon }\right)\right)
\end{equation}
for some universal constant $c$. The control condition of the iteration, in this case, is
\begin{equation*}
    \left\lVert \sqrt{D(\omega)}\mathbf{Y}\mathbf{Y} ^{\top } \sqrt{D(\omega)}-D(\omega )\right\rVert_{2}^{}\ge \gamma _{3},
\end{equation*}
where $D(\omega )=k ^{\star } \cdot \text{diag}\{\omega \}$. Inside each iteration, similar as before, we find $v$, which is the eigenvector associated with $\lambda $, and compute $\tau _{i}=\frac{v _{i}^{2}}{\omega _{i}}\mathbf{1}_{\left\{\omega _{i}>0\right\}}, 1 \le i \le N$. Next, we update $\omega $ as following
\begin{equation}\label{eq: update of omega 2}
    \omega _{i}=\left(1-\frac{\tau _{i}}{\max\limits _{i}\tau _{i}}\right)\omega _{i}, 1 \le i \le N.
\end{equation}
Finally, we check the condition whether $\left\lVert \omega \right\rVert_{1}^{}\le N(1-6\epsilon )$ and reject the null directly if this holds and pass to the following step if not.

\begin{algorithm}[htbp]
    \caption{Sample filtering when $N \le k ^{\star }$.}\label{algorithm n <= k}
    Set $\gamma _{3}$ as \eqref{eq: definition of gamma_3}, and $\lambda =\left\lVert \sqrt{D(\omega)}\mathbf{Y}\mathbf{Y} ^{\top } \sqrt{D(\omega)}-D(\omega )\right\rVert_{2}^{}$. ($\omega $ is initialized as $\mathbf{1}$.)

    \While {$\lambda \ge \gamma _{3}$}{
        Set $v$ to be the unit singular vector associated with $\lambda $

        Compute $\tau _{i}=\frac{v _{i}^{2}}{\omega _{i}}\mathbf{1}_{\left\{w _{i}>0\right\}}$

        Update $w$ according to \eqref{eq: update of omega 2}

        \If {$\left\lVert \omega \right\rVert_{1}^{}<N(1-6\epsilon )$} {

            \Return \texttt{None}
        }

        Set $\lambda =\left\lVert \sqrt{D(\omega)}\mathbf{Y}\mathbf{Y} ^{\top } \sqrt{D(\omega)}-D(\omega )\right\rVert_{2}^{}$
    }

    \Return $\omega $
\end{algorithm}

Similar to the algorithm \ref{algorithm n > k}, we will have the following results to guarantee the satisfaction of the second condition of \hyperref[def_omega_regularity]{$\omega $-regularity}.

\begin{lemma}[]\label{lemma: guarantee of omega regularity 2 2}
    Given the algorithm \ref{algorithm n <= k}, the following statements hold:\\
    \hspace*{0.5em}(i), the algorithm will terminate within finite rounds;\\ 
    \hspace*{0.5em}(ii), under $H _{0}$, the algorithm will not output rejection with high probability;\\
    \hspace*{0.5em}(iii), upon the finish of the algorithm, for any $S \subset [N]$ with $\left|S\right|\le \epsilon N$, we have 
    \begin{equation*}
        \left|\left\lVert \sum\limits_{i \in S}^{}\sqrt{\omega _{i}}Y _{i}\right\rVert_{2}^{2}-\left\lVert \omega _{S}\right\rVert_{1}^{}k ^{\star } \right|\le c \epsilon N \gamma _{3}
    \end{equation*}
    with probability higher than $1-\alpha $.
\end{lemma}

\begin{proof}
    Proof of (i) is same as the one in the proof of Lemma \ref{lemma: guarantee of omega regularity 2 1}.\opar

    To prove (ii), for the fixed samples $S$, define $\Lambda _{2}=\left\{\omega :\omega \in [0,1]^{N},\left\lVert \mathbf{1}_{[N] \backslash C}-\omega _{[N] \backslash C}\right\rVert_{1}^{}\le 5 \left\lVert \mathbf{1}_{C}-\omega _{C}\right\rVert_{1}^{}\right\}$. We will prove that $\omega ^{(t)}\in \Lambda _{2}$ before exiting the iteration. Suppose first this holds, we have $\left\lVert \mathbf{1}-\omega \right\rVert_{1}^{}\le 6 \left\lVert \mathbf{1}_{C}-\omega _{C}\right\rVert_{1}^{}\le 6\epsilon N$. This implies that we will always exit the iteration before the algorithm has chance to reject.

    Now assume $\omega ^{(t)} \in \Lambda _{2}$ and we are still in the iteration, we will prove $\omega ^{(t+1)}\in \Lambda _{2}$. By the definition, we only need to prove $5\sum\limits_{i \in C}^{}v _{i} ^{2}\ge \sum\limits_{i \in [N] \backslash C}^{}v _{i}^{2}$. Suppose this does not hold, then $\left\lVert v _{C}\right\rVert_{2}^{2}<\frac{1}{6}$. By Corollary \ref{corollary: bounds on the l2-operator norm of X^TX and XX^T} (set $\left\lVert \nu \right\rVert_{2}^{}=0$), we have $\left\lVert \sqrt{D (\omega _{[N] \backslash C}^{(t)})}Y _{[N] \backslash C}Y _{[N] \backslash C}^{^{\top } }\sqrt{D (\omega _{[N] \backslash C}^{(t)})} -D(\omega _{[N] \backslash C}^{(t)})\right\rVert_{2}^{}\le c \gamma _{3}\le c \lambda ^{(t)}$ with high probability. However, since $v$ is the eigenvector corresponding to the maximal eigenvalue, we have 
    \begin{equation*}
        \begin{aligned}
            \lambda ^{(t)}&=\left|v ^{\top }\left(\sqrt{D(\omega ^{(t)})}Y ^{\top } Y \sqrt{D(\omega ^{(t)})}-D(\omega ^{(t)})\right)v\right|\\ 
            &\le c \lambda ^{(t)}\left\lVert v _{[N] \backslash C}\right\rVert_{2}^{2}+2 \lambda ^{(t)}\left\lVert v _{[N] \backslash C}\right\rVert_{2}^{}\left\lVert v _{C}\right\rVert_{2}^{}+\lambda ^{(t)}\left\lVert v _{C}\right\rVert_{2}^{2}.
        \end{aligned}
    \end{equation*}
    Noting that $\left\lVert v _{C}\right\rVert_{2}^{2}+\left\lVert v _{[N]\backslash C}\right\rVert_{2}^{2}=1$, define $f _{c}(x)=c+(1-c)x+2 \sqrt{x(1-x)}, 0 \le x<\frac{1}{6}$. Note that by selecting the constant in \eqref{eq: definition of gamma_3} sufficiently large, we can let $c \rightarrow 0$. It can be proved that in this case $0 \le f _{0}(x)<1$, which is a contradiction. (ii) is proved.\opar

    Upon the termination of the algorithm without rejection, we have $\left\lVert \sqrt{D(\omega)}\mathbf{Y}\mathbf{Y} ^{\top } \sqrt{D(\omega)}-D(\omega )\right\rVert_{2}^{}<\gamma _{3}$ and $\left\lVert \omega \right\rVert_{1}^{}\ge (1-6\epsilon )N$. Therefore, we have 
    \begin{equation*}
        \begin{aligned}
            \left|\left\lVert \sum\limits_{i \in S}^{}\sqrt{\omega }_{i}Y _{i}\right\rVert_{2}^{2}-k ^{\star } \left\lVert \omega _{S}\right\rVert_{1}^{}\right|&=\left|\sum\limits_{i \in S}^{}\sum\limits_{j \in S}^{}\sqrt{\omega _{i}\omega _{j}}Y _{i}^{\top } Y _{j}-k ^{\star } \left\lVert \omega _{S}\right\rVert_{1}^{}\right|\\ 
            &=\left|\mathbf{1}_{S}^{\top } \left[\sqrt{D(\omega )}Y Y ^{\top } \sqrt{D(\omega )}-D(\omega )\right]\mathbf{1}_{S}\right|\overset{\text{(i)}}{\le }c \epsilon N \gamma _{3}.
        \end{aligned}
    \end{equation*}
    Here (i) uses the condition that $\left\lVert \sqrt{D(\omega)}\mathbf{Y}\mathbf{Y} ^{\top } \sqrt{D(\omega)}-D(\omega )\right\rVert_{2}^{}\le \gamma _{3}$. The proof is completed.
\end{proof}

After Algorithms \ref{algorithm n > k} or \ref{algorithm n <= k}, we have filtered the samples to ensure the second condition of \hyperref[def_omega_regularity]{$\omega $-regularity}. Now we turn to the check of the final condition.

Actually, the weight $\omega $ output by Algorithms \ref{algorithm n > k} or \ref{algorithm n <= k} already provides us with some useful insight.

\begin{lemma}[]\label{lemma: weight filtering 1}
    Set $\beta _{1}$, $\beta _{2}$ according to \eqref{eq: specified beta_1 and beta_2} and $\gamma _{i}, i \in \left\{2,3\right\}$ according to \eqref{eq: definition of gamma_2} or \eqref{eq: definition of gamma_3} depending on whether $N \ge k ^{\star }$ or not. For the weight $\omega $ output by the Algorithms \ref{algorithm n > k} or \ref{algorithm n <= k} and any $S \subset [N] \backslash C$ with $\left|S\right|\le \epsilon N$, we have 
    \begin{equation*}
        \left|\sum\limits_{i \in S}^{}\sum\limits_{j=1}^{N}\sqrt{\omega _{i}\omega _{j}}\left\langle Y _{i}, Y _{j} \right\rangle-k ^{\star } \left\lVert \omega _{S}\right\rVert_{1}^{}\right|\le c (\sqrt{N}\beta _{1}+\beta _{2}+\epsilon N \gamma )
    \end{equation*}
    with high probability.
\end{lemma}

\begin{proof}
    By Lemma \ref{lemma: good_samples_regularity}, with high probability the original samples $\tilde{\mathbf{Y}}$ are $(2\epsilon ,\beta _{1},\beta _{2})$-regular. By spliting $[N]=([N] \backslash C)\cup C$, we have
    \begin{equation*}
        \begin{aligned}
            \left|\sum\limits_{i \in S}^{}\sum\limits_{j=1}^{N}\sqrt{\omega _{i}\omega _{j}}\left\langle Y _{i}, Y _{j} \right\rangle-k ^{\star } \left\lVert \omega _{S}\right\rVert_{1}^{}\right|&=\underbrace{\left|\sum\limits_{i \in S}^{}\sum\limits_{j=1}^{N}\sqrt{\omega _{i}\omega _{j}}\left\langle Y _{i}, \tilde{Y}_{j} \right\rangle-k ^{\star } \left\lVert \omega _{S}\right\rVert_{1}^{}\right|}_{\text{I}}\\ 
            &+\underbrace{\left|\sum\limits_{i \in S}^{}\sum\limits_{j \in C}^{}\sqrt{\omega _{i}\omega _{j}}\left\langle Y _{i}, \tilde{Y}_{j} \right\rangle\right|}_{\text{II}}+\underbrace{\left|\sum\limits_{i \in S}^{}\sum\limits_{j \in C}^{}\sqrt{\omega _{i}\omega _{j}}\left\langle Y _{i}, Y _{j} \right\rangle\right|}_{\text{III}}
        \end{aligned}
    \end{equation*}
For the term I, we further have 
\begin{equation*}
    \text{I}\le \underbrace{\left|\sum\limits_{i \in S}^{}\sum\limits_{j=1}^{N}\sqrt{\omega _{i}}\left\langle Y _{i},  \tilde{Y}_{j}\right\rangle-k ^{\star } \left\lVert \sqrt{\omega _{S}}\right\rVert_{1}^{}\right|}_{\text{I}_{1}}+\underbrace{\left|\sum\limits_{i \in S}^{}\sum\limits_{j=1}^{N}\sqrt{\omega _{i}}(1-\sqrt{\omega _{j}})\left\langle Y _{i}, \tilde{Y}_{j} \right\rangle-k ^{\star }\left\lVert \sqrt{\omega _{S}}(1-\sqrt{\omega _{S}})\right\rVert_{1}^{}\right|}_{\text{I}_{2}}
\end{equation*}
We now bound the terms $\text{I}_{1},\text{I}_{2},\text{II},\text{III}$ respectively.

\textbf{Bound on $\text{I}_{1}$.} Because $\left|S\right|\le \epsilon N$, from \hyperref[def_regularity]{$\omega $-regularity} and Corollary \ref{corollary: generalized result of regularity 2}, we know $\text{I}_{1}\le c _{1}\sqrt{N}\beta _{1}$.

\textbf{Bound on $\text{I}_{2}$.}  Because $\omega $ is output by the Algorithms \ref{algorithm n > k} or \ref{algorithm n <= k}, we have $\left\lVert 1- \sqrt{\omega }\right\rVert_{1}^{}\le \left\lVert 1-\omega \right\rVert_{1}^{}\le \epsilon N$. By Corollary \ref{corollary: generalized result of regularity 2} again, we have $\text{I}_{2}\le c _{2}\beta _{2}$.

\textbf{Bound on II.} From Corollary \ref{corollary: generalized result of regularity 2} and the fact $S \cap C=\varnothing $, we know $\text{II}\le c _{3}\beta _{2}$.

\textbf{Bound on III.} From Lemma \ref{lemma: guarantee of omega regularity 2 1}, \ref{lemma: guarantee of omega regularity 2 2}, and the polarization equality, we have $\text{III}\le c _{4}\epsilon N \gamma $.

Combining them together, we come to the conclusion of the lemma.
\end{proof}

We are now prepared for the final algorithm, which is designed to ensure the third condition of \hyperref[def_omega_regularity]{$\omega $-regularity}. Given samples $S$ and weight $\omega $ by the algorithms \ref{algorithm n > k} or \ref{algorithm n <= k}, we compute $\tau _{i}=\left|\left\langle \sqrt{\omega _{i}}Y _{i}, \sum\limits_{j=1}^{N}\sqrt{\omega _{j}}Y _{j} \right\rangle-\omega _{i}k ^{\star }\right|\mathbf{1}_{\left\{\omega _{i}>0\right\}}$. Then, we set $\omega _{i}$ to zero for the indices $i$ corresponding to the $\epsilon N$ maximal $\tau _{i}$.

\begin{algorithm}[htbp]
    \caption{Weight filtering.}\label{algorithm: weight filtering}
    Compute $\tau _{i}=\left|\left\langle \sqrt{\omega _{i}}Y _{i}, \sum\limits_{j=1}^{N}\sqrt{\omega _{j}}Y _{j} \right\rangle-\omega _{i}k ^{\star }\right|\mathbf{1}_{\left\{\omega _{i}>0\right\}}, 1 \le i \le N$

    Sort $\tau _{i}$ by desreasing order and find the indices $\left\{i _{1},i _{2},\dots,i _{\epsilon N}\right\}$ corresponding to the first $\epsilon N$ maximal $\tau _{i}$

    Set $\omega _{i _{1}},\omega _{i _{2}},\dots,\omega _{i _{\epsilon N}}$ to zero

    \Return $\omega $
\end{algorithm}

\begin{lemma}[]\label{lemma: guarantee of omega regularity 3}
    For the weight output by Algorithm \ref{algorithm: weight filtering} and any $S \subset [N]$ with $\left|S\right|\le \epsilon N$, we have 
    \begin{equation*}
        \left|\sum\limits_{i \in S}^{}\sum\limits_{j=1}^{N}\sqrt{\omega _{i}\omega _{j}}\left\langle Y _{i}, Y _{j} \right\rangle-k ^{\star } \left\lVert \omega _{S}\right\rVert_{1}^{}\right|\le c(\sqrt{N}\beta _{1}+\beta _{2}+\epsilon N \gamma ).
    \end{equation*}
\end{lemma}

This lemma is very similar to Lemma \ref{lemma: weight filtering 1} but stronger in the sense that we allow any $S \subset [N]$ rather than $S \subset [N]\backslash C$. This advantage is exactly granted by the sorting and pruning processes in Algorithm \ref{algorithm: weight filtering}.

\begin{proof}
    The proof of this lemma is not hard but depends on the following fact. Because we have proved in Lemma \ref{lemma: weight filtering 1} that for any $S \subset [N] \backslash C$ and $\left|S\right|\le \epsilon N$, we have $\left|\sum\limits_{i \in S}^{}\sum\limits_{j=1}^{N}\sqrt{\omega _{i}\omega _{j}}\left\langle Y _{i}, Y _{j} \right\rangle-k ^{\star } \left\lVert \omega _{S}\right\rVert_{1}^{}\right|\le c (\sqrt{N}\beta _{1}+\beta _{2}+\epsilon N \gamma )$, we know that this bound also holds for any $S \subset [N] \backslash \left\{1,2,\dots,\epsilon N\right\}$. With the proning on $\omega $ to remove the first $\epsilon N$ coordinates, we can assume without loss of generality that $S \subset [N] \backslash \left\{1,2,\dots\epsilon N\right\}$. We have 
    \begin{equation*}
        \left|\sum\limits_{i \in S}^{}\sum\limits_{j=1}^{N}\sqrt{\omega _{i}\omega _{j}}\left\langle Y _{i}, Y _{j} \right\rangle-k ^{\star } \left\lVert \omega _{T}\right\rVert_{1}^{}\right|=\left|\sum\limits_{i \in S}^{}\sum\limits_{j=\epsilon N+1}^{N}\sqrt{\omega _{i}\omega _{j}}\left\langle Y _{i}, Y _{j} \right\rangle-k ^{\star } \left\lVert \omega _{T}\right\rVert_{1}^{}\right|\overset{\text{(1)}}{\le } c (\sqrt{N}\beta _{1}+\beta _{2}+\epsilon N\gamma ).
    \end{equation*}
    Here (1) can be obtained by considering the subtraction of summation on the original weight with indices between $1$ and $\epsilon N$ from the summation on the whole original weight, and using the fact that $S \cap \left\{1,2,\dots,\epsilon N\right\}=\varnothing $. 
\end{proof}

We have finished the proofs of Lemma \ref{lemma: guarantee of omega regularity 1}, \ref{lemma: guarantee of omega regularity 2 1}, \ref{lemma: guarantee of omega regularity 2 2}, and \ref{lemma: guarantee of omega regularity 3}, which means we have verified the correctness of \hyperref[def_omega_regularity]{$\omega $-regularity} for our samples with high probability. The final step is to exploit the property of \hyperref[def_omega_regularity]{$\omega $-regularity} to distinguish between $H _{0}$ and $H _{1}$.

\begin{lemma}\label{lemma: property of the filtered weights}
    Let $\omega $ be the weights after running through the algorithm \ref{algorithm: prefiltering}, \ref{algorithm n > k} or \ref{algorithm n <= k}, and \ref{algorithm: weight filtering}, and $\beta _{1},\beta _{2}$ seleced as in Lemma \ref{lemma: good_samples_regularity}. Then, with high probability we have 
    \begin{equation}\label{eq: property of the filtered weights}
        \left|\left\lVert \sum\limits_{i=1}^{N}\sqrt{\omega _{i}}Y _{i}\right\rVert_{2}^{2}-k ^{\star } \left\lVert \omega \right\rVert_{1}^{}-N ^{2}\left\lVert \nu \right\rVert_{2}^{2}\right|\le c \left(\sqrt{N}\beta _{1}+\beta _{2}+\epsilon N \gamma +\underbrace{\left(N ^{\frac{3}{2}}\left\lVert \nu \right\rVert_{2}^{}+\sqrt{Nk ^{\star }}\right)\sqrt{\ln \left(\frac{1}{\alpha }\right)}+\ln \left(\frac{1}{\alpha }\right)}_{:=\text{R}}\right).
    \end{equation}
\end{lemma}

\begin{proof}
    The proof of this lemma is from the definition of \hyperref[def_omega_regularity]{$\omega $-regularity} and straightforward calculations. To evaluate $\left\lVert \sum\limits_{i=1}^{N}\sqrt{\omega _{i}}Y _{i}\right\rVert_{2}^{2}$, we write it as:
    \begin{equation*}
        \begin{aligned}
            \left\lVert \sum\limits_{i=1}^{N}\sqrt{\omega _{i}}Y _{i}\right\rVert_{2}^{2}&=\left\lVert \sum\limits_{i \in [N] \backslash C}^{}\sqrt{\omega _{i}}Y _{i}+\sum\limits_{i \in C}^{}\sqrt{\omega _{i}}Y _{i}\right\rVert_{2}^{2}\\ 
            &=\underbrace{\left\lVert \sum\limits_{i \in [N] \backslash C}^{}\sqrt{\omega _{i}}Y _{i}\right\rVert_{2}^{2}}_{\text{I}}-\underbrace{\left\lVert \sum\limits_{i \in C}^{}\sqrt{\omega _{i}}Y _{i}\right\rVert_{2}^{2}}_{\text{II}}+\underbrace{2 \sum\limits_{i=1}^{N}\sum\limits_{j \in C}^{}\sqrt{\omega _{i}\omega _{j}}Y _{i}^{\top } Y _{j}}_{\text{III}}
        \end{aligned}
    \end{equation*}
We are going to bound the terms individually.

\textbf{Bound on I.} Define the vector $r \in \mathbb{R}^{N}$ where $r _{i}=1$ for $i \in C$ and $1-\sqrt{\omega _{i}}$ for $i \in [N] \backslash C$. By Algorithm \ref{algorithm n > k} or \ref{algorithm n <= k}, we have $\left\lVert r\right\rVert_{1}^{}\le \epsilon N+6\epsilon N=7\epsilon N$. For the term I, we have
\begin{equation*}
    \left\lVert \sum\limits_{i \in [N] \backslash C}^{}\sqrt{\omega _{i}}Y _{i}\right\rVert_{2}^{2}=\underbrace{\left\lVert \sum\limits_{i=1}^{N}\tilde{Y}_{i}\right\rVert_{2}^{2}}_{\text{I}_{1}}+\underbrace{\left\lVert \sum\limits_{i=1}^{N}r _{i}\tilde{Y}_{i}\right\rVert_{2}^{2}}_{\text{I}_{2}}-\underbrace{2 \sum\limits_{i=1}^{N}\sum\limits_{j=1}^{N}r _{j}\tilde{Y}_{i}^{\top }\tilde{Y}_{j}}_{\text{I}_{3}} 
\end{equation*}
For $\text{I}_{1}$, we have 
\begin{equation*}
    \text{I}_{1}=\sum\limits_{i=1}^{N}\left\lVert \tilde{Y}_{i}-\nu \right\rVert_{2}^{2}+N ^{2}\left\lVert \nu \right\rVert_{2}^{2}+2N\nu ^{\top } \sum\limits_{i=1}^{N}(\tilde{Y}_{i}-\nu ).
\end{equation*}
By the standard concentration bounds (see \ref{def_sub_Gaussian_var} and \ref{def_sub_exp_var}), with high probability we have 
\begin{equation*}
    \left|\text{I}_{1}-Nk ^{\star }-N ^{2}\left\lVert \nu \right\rVert_{2}^{2}\right|\le c _{1}\left(N ^{\frac{3}{2}}\left\lVert \nu \right\rVert_{2}^{}\sqrt{\ln \left(\frac{1}{\alpha }\right)}+\sqrt{Nk ^{\star } \ln \left(\frac{1}{\alpha }\right)}+\ln \left(\frac{1}{\alpha }\right).\right)
\end{equation*}
For $\text{I}_{2}$, by the argument above, Lemma \ref{lemma: good_samples_regularity}, and Corollary \ref{corollary: generalized result of regularity 2} with appropriate choice of $\epsilon $, with high probability we have 
\begin{equation*}
    \left|\text{I}_{2}-k ^{\star } \left\lVert r\right\rVert_{2}^{2}\right|\le c _{2} \beta _{2}.
\end{equation*}
Again by Lemma \ref{lemma: good_samples_regularity} and Corollary \ref{corollary: generalized result of regularity 2}, with high probability we have 
\begin{equation*}
    \left|\text{I}_{3}-2k ^{\star } \left\lVert r\right\rVert_{1}^{}\right|\le c _{3} \sqrt{N}\beta _{1}.
\end{equation*}
Combining the results above and using the fact that $Nk ^{\star }+k ^{\star } \left\lVert r\right\rVert_{2}^{2}+2k ^{\star } \left\lVert r\right\rVert_{1}^{}=k ^{\star } \left\lVert \mathbf{1}-r\right\rVert_{2}^{2}=k ^{\star }\left\lVert \omega _{[N] \backslash C}\right\rVert_{1}^{}$, we have
\begin{equation*}
    \left|\text{I}-k ^{\star }\left\lVert \omega _{[N] \backslash C}\right\rVert_{1}^{}-N ^{2}\left\lVert \nu \right\rVert_{2}^{2}\right|\le c \left[\sqrt{N}\beta _{1}+\beta _{2}+\left(N ^{\frac{3}{2}}\left\lVert \nu \right\rVert_{2}^{}+\sqrt{Nk ^{\star }}\right)\sqrt{\ln \left(\frac{1}{\alpha }\right)}+\ln \left(\frac{1}{\alpha }\right)\right].
\end{equation*}
\textbf{Bound on II.} By Lemma \ref{lemma: guarantee of omega regularity 2 1} (recall the definition of $\gamma _{2}$ \eqref{eq: definition of gamma_2} and the fact that $N>k ^{\star }$ in this case) or \ref{lemma: guarantee of omega regularity 2 2}, with high probability we have
\begin{equation*}
    \left|\text{II}-k ^{\star }\left\lVert \omega _{C}\right\rVert_{1}^{}\right|\le c\epsilon N \gamma .
\end{equation*}
\textbf{Bound on III.} By Lemma \ref{lemma: guarantee of omega regularity 3}, with high probability we have 
\begin{equation*}
    \left|\text{III}-2 k ^{\star }\left\lVert \omega _{C}\right\rVert_{1}^{}\right|\le c \left(\sqrt{N}\beta _{1}+\beta _{2}+\epsilon N \gamma \right).
\end{equation*}
Now combining the bounds of I, II, and III together, we have 
\begin{equation*}
    \left|\left\lVert \sum\limits_{i=1}^{N}\sqrt{\omega _{i}}Y _{i}\right\rVert_{2}^{2}-k ^{\star } \left\lVert \omega \right\rVert_{1}^{}-N ^{2}\left\lVert \nu \right\rVert_{2}^{2}\right|\le c \left(\sqrt{N}\beta _{1}+\beta _{2}+\epsilon N \gamma +\left(N ^{\frac{3}{2}}\left\lVert \nu \right\rVert_{2}^{}+\sqrt{Nk ^{\star }}\right)\sqrt{\ln \left(\frac{1}{\alpha }\right)}+\ln \left(\frac{1}{\alpha }\right)\right),
\end{equation*}
which is the desired results.
\end{proof}

\begin{proof}[Proof of Theorem \ref{theorem: upper bound 2}]
    Substituting $\nu $ with $P ^{\star }\mu $, the key is to require each term on the right-hand-side in \eqref{eq: property of the filtered weights} to be asymptotically smaller than the term $N ^{2}\left\lVert P ^{\star }\mu \right\rVert_{2}^{2}$. For convenience, we rewrite the definition of $\beta _{1}$ and $\beta _{2}$ here.
    \begin{equation*}
        \begin{aligned}
            &\beta _{1}=\epsilon N \left(\sqrt{k ^{\star }}+\sqrt{N}\left\lVert P ^{\star }\mu \right\rVert_{2}^{} \right)\sqrt{\ln \left(\frac{N}{\alpha }\right)}+\epsilon N \ln \left(\frac{N}{\alpha }\right)+\epsilon N \sqrt{N}\left\lVert P ^{\star }\mu \right\rVert_{2}^{2},\\
            &\beta _{2}=\epsilon N \sqrt{\epsilon N k ^{\star } \ln \left(\frac{1}{\epsilon }\right)}+(\epsilon N)^{2}\ln \left(\frac{1}{\epsilon }\right)+\left\lVert P ^{\star }\mu \right\rVert_{2}^{} (\epsilon N)^{2}\sqrt{\ln \left(\frac{1}{\epsilon }\right)}+\left\lVert P ^{\star }\mu \right\rVert_{2}^{2}(\epsilon N)^{2}.
        \end{aligned}
    \end{equation*}
    The table below presents the terms and corresponding required conditions for $\left\lVert P ^{\star }\mu \right\rVert_{2}^{2}$. Although there are many terms, only a few of them are useful indeed.
    \begin{table}[htbp]
    \centering
    \begin{tabular}{|c|l|l|l|}
        \hline
        \textbf{Sources} & \textbf{No.} & \textbf{Terms} \hspace*{8em} & \textbf{Required Conditions for $\left\lVert P ^{\star }\mu \right\rVert_{2}^{}$}\\
        \hline
        \multirow{4}{*}{$\sqrt{N}\beta _{1}$} & 1 & $\epsilon N \sqrt{Nk ^{\star }}\sqrt{\ln \left(\frac{N}{\alpha }\right)}$ & $\left\lVert P ^{\star }\mu \right\rVert_{2}^{2}\gtrsim \frac{\epsilon \sqrt{k ^{\star } \ln \left(\frac{N}{\alpha }\right)}}{\sqrt{N}}$ \\ 
        \cline{2-4}
        & 2 & $\epsilon N ^{2}\left\lVert P ^{\star }\mu \right\rVert_{2}^{}\sqrt{\ln \left(\frac{N}{\alpha }\right)}$ & $\left\lVert P ^{\star }\mu \right\rVert_{2}^{}\gtrsim \epsilon \sqrt{\ln \left(\frac{N}{\alpha }\right)}$\\ 
        \cline{2-4}
        & 3 & $\epsilon N \sqrt{N}\ln \left(\frac{N}{\alpha }\right)$ & $\left\lVert P ^{\star }\mu \right\rVert_{2}^{2}\gtrsim \frac{\epsilon  \ln \left(\frac{N}{\alpha }\right)}{\sqrt{N}}$\\
        \cline{2-4}
        & 4 & $\epsilon N ^{2}\left\lVert P ^{\star }\mu \right\rVert_{2}^{2}$ & $\epsilon \lesssim 1$\\ 
        \hline
        \multirow{4}{*}{$\beta _{2}$} & 5 & $\epsilon N \sqrt{\epsilon N k ^{\star } \ln \left(\frac{1}{\epsilon }\right)}$ & $\left\lVert P ^{\star }\mu \right\rVert_{2}^{2}\gtrsim \frac{\epsilon \sqrt{\epsilon k ^{\star } \ln \left(\frac{1}{\epsilon }\right)}}{\sqrt{N}}$\\ 
        \cline{2-4}
        & 6 & $\epsilon ^{2}N ^{2}\ln \left(\frac{1}{\epsilon }\right)$ & $\left\lVert P ^{\star }\mu \right\rVert_{2}^{2}\gtrsim \epsilon ^{2}\ln \left(\frac{1}{\epsilon }\right)$\\ 
        \cline{2-4}
        & 7 & $\left\lVert P ^{\star }\mu \right\rVert_{2}^{}\epsilon ^{2}N ^{2}\sqrt{\ln \left(\frac{1}{\epsilon }\right)}$ & $\left\lVert P ^{\star }\mu \right\rVert_{2}^{2}\gtrsim \epsilon ^{2}\sqrt{\ln \left(\frac{1}{\epsilon }\right)}$\\ 
        \cline{2-4}
        & 8 & $\left\lVert P ^{\star }\mu \right\rVert_{2}^{2}\epsilon ^{2}N ^{2}$ & $\epsilon \lesssim 1$\\ 
        \hline
        \multirow{4}{*}{$\epsilon N \gamma $} & 9 & $\epsilon N \sqrt{Nk ^{\star }}$ & $\left\lVert P ^{\star }\mu \right\rVert_{2}^{2}\gtrsim \frac{\epsilon \sqrt{k ^{\star }}}{\sqrt{N}}$\\ 
        \cline{2-4}
        & 10 & $\epsilon N \sqrt{N \ln \left(\frac{1}{\alpha }\right)}$ & $\left\lVert P ^{\star }\mu \right\rVert_{2}^{2}\gtrsim \frac{\epsilon \sqrt{\ln \left(\frac{1}{\alpha }\right)}}{\sqrt{N}}$\\ 
        \cline{2-4}
        & 11 & $\epsilon N \ln \left(\frac{1}{\alpha }\right)$ & Dominated by the term 3\\ 
        \cline{2-4}
        & 12 & $\epsilon ^{2}N ^{2}\ln \left(\frac{1}{\epsilon }\right)$ & $\left\lVert P ^{\star }\mu \right\rVert_{2}^{2}\gtrsim \epsilon ^{2}\ln \left(\frac{1}{\epsilon }\right)$\\
        \hline
        \multirow{3}{*}{R} & 13 & $N ^{3/2}\left\lVert P ^{\star }\mu \right\rVert_{2}^{}\sqrt{\ln \frac{1}{\alpha }}$ & $\left\lVert P ^{\star }\mu \right\rVert_{2}^{}\gtrsim \frac{\sqrt{\ln \left(\frac{1}{\alpha }\right)}}{\sqrt{N}}$\\ 
        \cline{2-4}
        & 14 & $N \sqrt{k ^{\star }}\ln \left(\frac{1}{\alpha }\right)$ & $\left\lVert P ^{\star }\mu \right\rVert_{2}^{2}\gtrsim \frac{\sqrt{k ^{\star }}\ln \left(\frac{1}{\alpha }\right)}{N}$\\ 
        \cline{2-4}
        & 15 & $\ln \left(\frac{1}{\alpha }\right)$ & Dominated by the term 3\\ 
        \hline
    \end{tabular}
    \caption{terms from \eqref{eq: property of the filtered weights} and corresponding required conditions on $\left\lVert P ^{\star }\mu \right\rVert_{2}^{}$.}
\end{table}

Inspecting the column ``\textbf{Required Conditions}'' and computing the maximal, we conclude the final requirement for $\left\lVert P ^{\star }\mu \right\rVert_{2}^{}$:
\begin{equation}\label{eq: raw required condition for the upper bound 2}
   \left\lVert P ^{\star }\mu \right\rVert_{2}^{2}\gtrsim \max\limits \left\{\frac{\epsilon \ln \left(\frac{N}{\alpha }\right)}{\sqrt{N}},\epsilon ^{2}\ln \left(\frac{N}{\alpha }\right),\sqrt{\frac{\epsilon ^{2}k ^{\star } \ln \left(\frac{N}{\alpha }\right)}{N}},\frac{\sqrt{k ^{\star }}\ln \left(\frac{1}{\alpha }\right)}{N}\right\}.
\end{equation}
Once this requirement is satisfied, by Lemma \ref{lemma: property of the filtered weights} and the inequality \eqref{eq: property of the filtered weights}, we know that $\left|\left\lVert \sum\limits_{i=1}^{N}\sqrt{\omega _{i}}P ^{\star }X _{i}\right\rVert_{2}^{2}-k ^{\star } \left\lVert \omega \right\rVert_{1}^{}\right|$ will be less than $c _{1} N ^{2}P(\epsilon ,k ^{\star },N,1)$ under $H _{0}$ and greater than $c _{2}N ^{2}P(\epsilon ,k ^{\star },N,1)$ under $H _{1}$ for some $0<c _{1}<c _{2}<1$, and both with high probability. $c _{1}$ and $c _{2}$ are universal and their values can be determined following proper choices of $\epsilon ,\alpha ,N,k ^{\star }$. From the argument at the end of the \hyperref[subsection: proof of the main upper bound 1]{proof} of Theorem \ref{theorem: upper bound 1}, we know \eqref{eq: raw required condition for the upper bound 2} leads to \eqref{eq: required condition of the main upper bound 2}. The proof is completed.
\end{proof}

\subsection{Proof of Corollary \ref{corollary: matching between bounds}}\label{subsection: proof of the main corollary}

\begin{proof}
By Definition \ref{def_first_optimal_dimension} and \ref{def_second_optimal_dimension}, we can rewrite $k _{1}^{\star }$ and $k _{2}^{\star }$ as 
\begin{equation*}
    \begin{aligned}
        & k _{1}^{\star }:=\max\limits _{k}\left\{k \left\lvert\right. k \ge 0, \frac{D _{k-1}(K)}{k ^{1/4}}\ge \frac{\sigma }{\sqrt{N}}\right\},\\ 
        & k _{2}^{\star }:=\max\limits _{k}\left\{k \left\lvert\right. k \ge 0, \frac{D _{k-1}(K)}{k ^{1/4}}\ge \frac{\sqrt{\epsilon }}{N ^{\frac{1}{4}}}\sigma \right\}.
    \end{aligned}
\end{equation*}
For a fixed $K$, define function $g(k)=\frac{D _{k-1}(K)}{k ^{1/4}}$. By the property of Kolmogorov widths, we know $g$ is a decreasing function. Therefore, in the alternative definition above, if $\frac{\sigma }{\sqrt{N}}\ge \frac{\sqrt{\epsilon }}{N ^{1/4}}\sigma $, i.e., $\epsilon \le \frac{1}{\sqrt{N}}$, we have $k _{1}^{\star }\le k _{2}^{\star }$; similarly, if $\epsilon \ge \frac{1}{\sqrt{N}}$, we have $k _{1}^{\star }\ge k _{2}^{\star }$. We discuss the two possibilities in the following.

\textbf{(i): $\epsilon \le \frac{1}{\sqrt{N}}$.} In this case, the dominant term in $E ^{2}(\epsilon ,K,N,\sigma ^{2})$ is $\frac{\sqrt{k _{1}^{\star }}}{N}\sigma ^{2}$. If $k _{1}^{\star }=k _{2}^{\star }$. Obviously, we know $\frac{\sqrt{k _{1}^{\star }}}{N}\sigma ^{2}$ in \eqref{eq: main theoretical lower bound} will also be the dominant. Otherwise, we have $k _{1}^{\star }+1\le  k _{2}^{\star }$. By the definition and the monotone property of Kolmogorov widths, we have 
\begin{equation*}
    \frac{\left(k _{2}^{\star }\right)^{\frac{1}{4}}\sqrt{\epsilon }}{N ^{\frac{1}{4}}}\sigma \le D _{k _{2}^{\star }-1}\le  D _{k _{1}^{\star }}(K)\le \frac{\left(k _{1}^{\star }+1\right)^{\frac{1}{4}}}{\sqrt{N}}\sigma .
\end{equation*}
This leads to 
\begin{equation*}
    \epsilon \le \frac{1}{\sqrt{N}}\sqrt{\frac{k _{1}^{\star }+1}{k _{2}^{\star }}}\lesssim \frac{1}{\sqrt{N}}\sqrt{\frac{k _{1}^{\star }}{k _{2}^{\star }}}.
\end{equation*}
With a little calculation, one can verify that this is equivalent to saying $\frac{\sqrt{k _{1}^{\star }}}{N}\sigma ^{2}$ is the dominant in \eqref{eq: main theoretical lower bound} again. Hence, \eqref{eq: main theoretical lower bound} matches $E (\epsilon ,K,N,\sigma ^{2})$ in this possibility.

\textbf{(ii): $\epsilon \ge \frac{1}{\sqrt{N}}$.} In this case, the dominant term in $E ^{2}(\epsilon ,K,N,\sigma ^{2})$ is either $\epsilon \sqrt{\ln \left(\frac{1}{\epsilon }\right)\frac{k _{2}^{\star }}{N}}\sigma ^{2}$ or $\epsilon ^{2}\ln \left(\frac{1}{\epsilon }\right)\sigma ^{2}$. If $k _{1}^{\star }=k _{2}^{\star }$, then $\epsilon \sqrt{\frac{k _{2}^{\star }}{N}}\sigma ^{2}$ will dominate over $\frac{\sqrt{k _{1}^{\star }}}{N}\sigma ^{2}$ and consequently the dominant term in \eqref{eq: main theoretical lower bound} is either $\epsilon \sqrt{\frac{k _{2}^{\star }}{N}}\sigma ^{2}$ or $\epsilon ^{2}\sigma ^{2}$, which nearly matches $E ^{2}(\epsilon ,K,N,\sigma ^{2})$ with only an exception of a $\ln \left(\frac{1}{\epsilon }\right)$ factor. Otherwise, we have $k _{1}^{\star }\ge k _{2}^{\star }+1$. Again, by the definition and the monotone property of Kolmogorov widths, we have 
\begin{equation*}
    \frac{\left(k _{1}^{\star }\right)^{\frac{1}{4}}}{\sqrt{N}}\sigma \le D _{k _{1}^{\star }-1}(K)\le D _{k _{2}^{\star }}(K)\le \frac{(k _{2}^{\star }+1)^{\frac{1}{4}}\sqrt{\epsilon }}{N ^{\frac{1}{4}}}\sigma ,
\end{equation*}
which leads to 
\begin{equation*}
    \epsilon \ge \frac{1}{\sqrt{N}}\sqrt{\frac{k _{1}^{\star }}{k _{2}^{\star }+1}}.
\end{equation*}
If $k _{2}^{\star }=0$, then the dominant terms in \eqref{eq: main theoretical lower bound} and $E ^{2}(\epsilon ,K,N,\sigma ^{2})$ are $\epsilon ^{2}\sigma ^{2}$ and $\epsilon ^{2}\ln \left(\frac{1}{\epsilon }\right)\sigma ^{2}$, respectively, which match except for a $\ln \left(\frac{1}{\epsilon }\right)$ term. Otherwise, we will have 
\begin{equation*}
    \epsilon \gtrsim \frac{1}{\sqrt{N}}\sqrt{\frac{k _{1}^{\star }}{k _{2}^{\star }}},
\end{equation*}
which again can be used to verify that $\epsilon \sqrt{\frac{k _{2}^{\star }}{N}}\sigma ^{2}$ dominates over $\frac{\sqrt{k _{1}^{\star }}}{N}\sigma ^{2}$ in \eqref{eq: main theoretical lower bound}. In both cases, \eqref{eq: main theoretical lower bound} and $E ^{2}(\epsilon ,K,N,\sigma ^{2})$ match with only an exception of a $\ln \left(\frac{1}{\epsilon }\right)$ factor. The other part of the corollary is directly from the expression of $E ^{2}(\epsilon ,K,N,\sigma ^{2})$ and $P ^{2}(\epsilon ,K,N,\sigma ^{2})$. The proof is completed.
\end{proof}

\subsection{Proof of Theorem \ref{theorem: lp main lower bound}}\label{subsection: proof of the lp main lower bound}

\begin{proof}
    The techniques are parallel to the proof of the lower bounds above. By Lemma \ref{lemma: general existence of the optimal vector in a QCO set}, for $k=k _{p,1}^{\star }, c=\frac{(k _{p,1}^{\star })^{1/p-1/4}}{\sqrt{N}}\sigma $, there exists a $\theta \in K$ such that $\left\lVert \theta \right\rVert_{p}^{}=\frac{(k _{p,1}^{\star })^{1/p-1/4}}{\sqrt{N}}\sigma $ and $\left\lVert \theta \right\rVert_{\infty }^{}\le \frac{1}{(k _{p,1}^{\star })^{1/4}\sqrt{N}}\sigma $. Again, we impose a uniform prior on the set $\left\{\mu \left\lvert\right. \mu \in K,\mu _{i}=\kappa \gamma _{i}\theta _{i},\gamma _{i}\in \left\{-1,1\right\},i \in I\right\}$ for some positive constant $\kappa $. By the same arguments of Appendix \ref{subsection: proof of the main lower bound 1}, if $\mathbb{E}_{\mu ,\mu ^{\prime }}\exp\left\{\frac{N\mu ^{\top } \mu ^{\prime }}{\sigma ^{2}}\right\}-1$, or equivalently $\frac{N ^{2}\left\lVert \theta \right\rVert_{4}^{4}\kappa ^{4}}{2 \sigma ^{4}}$ is sufficiently small, then there is no test with uniformly small type-I and type-II errors. We have 
    \begin{equation*}
        \frac{N ^{2}\left\lVert \theta \right\rVert_{4}^{4}\kappa ^{4}}{2\sigma ^{4}}\le \frac{N ^{2}\sum\limits_{i=1}^{d}\left\lVert \theta \right\rVert_{\infty }^{4-p}\left|\theta _{i}\right|^{p}\kappa ^{4}}{2\sigma ^{4}}=\frac{N ^{2}\left\lVert \theta \right\rVert_{\infty }^{4-p}\left\lVert \theta \right\rVert_{p}^{p}\kappa ^{4}}{2\sigma ^{4}}\le \kappa ^{4}.
    \end{equation*}
    This implies that it is necessary to ensure $\kappa \gtrsim 1$, from which we establish the first lower bound that $\rho _{\text{critical}}\gtrsim \frac{(k _{p,1}^{\star })^{1/p-1/4}}{\sqrt{N}}\sigma $.

    Similarly, it suffices to bound $\frac{\sqrt{N}\left\lVert v\right\rVert_{4}^{2}}{\sigma ^{2}}$ for the second term, where $v$ is uniform from the set
    \begin{equation*}
        V :=\left\{v \left\lvert\right. v _{i}=\gamma _{i}\theta _{i},\gamma _{i}\in \left\{-1,1\right\},i \in I\right\},
    \end{equation*}
    and $\theta \in K$ from Lemma \ref{lemma: general existence of the optimal vector in a QCO set} satisfies $\left\lVert \theta \right\rVert_{p}^{}=\frac{\sqrt{\epsilon }(k _{p,2}^{\star })^{1/p-1/4}}{N ^{1/4}}\sigma $ and $\left\lVert \theta \right\rVert_{\infty }^{}\le \frac{\sqrt{\epsilon }}{(Nk _{p,2}^{\star })^{1/4}}\sigma $. Consequently, it is not difficult to verify that $\frac{\sqrt{N}\left\lVert v\right\rVert_{4}^{2}}{\sigma ^{2}}\le \epsilon $. By Lemma \ref{lemma: main lemma in the proof of lower bound 2}, we establish the second lower bound that $\rho _{\text{critical}}\gtrsim \frac{\sqrt{\epsilon }(k _{p,2}^{\star })^{1/p-1/4}}{N ^{1/4}}\sigma $.

    Suppose now that $\theta $ is the vector from Lemma \ref{lemma: general existence of the optimal vector in a QCO set} when $k=k _{p,3}^{\star }$. Subsequently, we know $\left\lVert \theta \right\rVert_{2}^{2}\le \left\lVert \theta \right\rVert_{p}^{p}\cdot \left\lVert \theta \right\rVert_{\infty }^{2-p}\le \epsilon ^{2}\sigma ^{2}$. Following the same reasoning in Appendix \ref{subsection: proof of the main lower bound 3}, we know there always exists a strategy for $\mathcal{C}$ to discourage the testing. Therefore, we establish the third lower bound. The proof is completed.
\end{proof}

\subsection{Proof of Theorem \ref{theorem: lp main upper bound}}\label{subsection: proof of lp main upper bound}

Before we prove Theorem \ref{theorem: lp main upper bound}, we need the following observation of the p-convex orthosymmetric set.

\begin{lemma}[Truncation of a $p$-convex orthosymmetric set]\label{lemma: property of truncation}
    Assume $K \subset \mathbb{R}^{d}$ is a $p$-convex orthosymmetric set. For any $1 \le j \le d$ and $J:=\left\{i _{1},\dots,i _{j}\right\}\subset \left\{1,\dots,d\right\}$, define the Truncation of $K$ on $J$ as
    \begin{equation*}
        K _{J}:=\left\{\theta _{J}\left\lvert\right. \theta _{J}=\left(\theta _{i _{1}},\dots,\theta _{i _{j}}\right)\in \mathbb{R}^{j},\exists \beta \in K, \beta _{J}=\theta _{J}\right\},
    \end{equation*}
    where $\beta _{J}=(\beta _{i _{1}},\dots,\beta _{i _{j}})$ is the truncated $\beta $ on the index set $J$. Then $K _{J}\subset \mathbb{R}^{j}$ is also $p$-convex orthosymmetric. 
\end{lemma}

\begin{proof}
    Lemma \ref{lemma: property of truncation} is straightforward from the definition. For any $\theta _{J,1},\theta _{J,2}\in K _{J}$, suppose $\theta _{J,1}=\beta _{1,J},\theta _{J,2}=\beta _{2,J},\beta _{1},\beta _{2} \in K$, then we know for any $\lambda \in [0,1]$, we have $\lambda \theta _{J,1}^{p}+(1-\lambda )\theta _{J,2}^{p}=\lambda \beta _{1,J}^{p}+(1-\lambda )\beta _{2,J}^{p}=(\lambda \beta _{1}^{p}+(1-\lambda )\beta _{2}^{p})_{J}\overset{\text{(i)}}{=}\tilde{\beta }_{J}^{p}\in K _{J}^{p}$, where $\tilde{\beta }\in K ^{p}$ and (i) is from that $K$ is p-convex. Consequently, $K _{J}$ is p-convex.

    Suppose $\bm{\xi }=(\xi _{1},\dots,\xi _{d}), \xi _{i}\in \left\{-1,1\right\}$, then for any $\theta _{J}\in K _{J}$, we have $\bm{\xi } _{J}\odot \theta _{J}=\bm{\xi } _{J}\odot \beta _{J}=(\bm{\xi }\odot \beta )_{J}\overset{(\text{ii})}{=}\bar{\beta }_{J}\in K _{J}$, where $\bar{\beta }\in K$, $\odot $ is element-wise multiplication and (ii) is from that $K$ is orthosymmetric. Consequently, $K _{J}$ is also orthosymmetric.
\end{proof}

We are now prepared to prove the original theorem.

\begin{proof}
    Suppose first $\epsilon \le \frac{1}{\sqrt{N}}$. In this case, we have $\frac{j ^{1/p-1/4}}{\sqrt{N}}\sigma \ge \frac{\sqrt{\epsilon }j ^{1/p-1/4}}{N ^{1/4}}\sigma $ for any $0 \le j \le d$. By the monotonicity of Kolmogorov widths and the definitions of $k _{p,1}^{\star }$ and $k _{p,2}^{\star }$, we know $k _{p,1}^{\star }\le k _{p,2}^{\star }$. Later calculations show that in this case we also have $k _{p,2}^{\star }\le k _{p,3}^{\star }$. This implies $k _{p,1}^{\star }\le k _{p,2}^{\star }\le k _{p,3}^{\star }$. For $k _{p,1}^{\star }$ and its corresponding $P _{p,1}^{\star }$, without loss of generality we assume $P _{p,1}^{\star }\theta =(\theta _{1},\dots,\theta _{k _{p,1}^{\star }},0,\dots,0)$. Regard $P _{p,1}^{\star }\theta $ as a vector in $\mathbb{R}^{k _{p,1}^{\star }}$, by the H\"older's inequality \ref{lemma: holder's inequality}, we know $\left\lVert P _{p,1}^{\star }\theta \right\rVert_{2}^{}\gtrsim \frac{(k _{p,1}^{\star })^{1/4}}{\sqrt{N}}\sigma $ for any $\theta \in K$. Further define $k _{p,1,1}^{\star }:=\max\limits \left\{j \left\lvert\right. 0 \le j \le k _{p,1}^{\star }, D _{j-1}(K)>\frac{j ^{1/4}}{\sqrt{N}}\sigma \right\}$. Since $\epsilon \le \frac{1}{\sqrt{N}}$, by the \hyperref[subsection: proof of the main corollary]{proof} of Corollary \ref{corollary: matching between bounds} and Lemma \ref{lemma: property of truncation}, we know that if $\left\lVert P _{p,1}^{\star }\theta \right\rVert_{2}^{}\gtrsim \frac{(k _{p,1,1}^{\star })^{1/4}}{\sqrt{N}}\sigma $, then there exists a desired test. However, this is already guaranteed since $\left\lVert P _{p,1}^{\star }\theta \right\rVert_{2}^{}\gtrsim \frac{(k _{p,1}^{\star })^{1/4}}{\sqrt{N}}\sigma $ and $k _{p,1}^{\star }\ge k _{p,1,1}^{\star }$.

    Suppose now $\epsilon \ge \frac{1}{\sqrt{N}}$. Analogous to the discussion above, we know that $k _{p,2}^{\star }\le k _{p,1}^{\star }$. If $k _{p,2}^{\star }=k _{p,3}^{\star }$, then we can assume $P _{p,2}^{\star }=P _{p,3}^{\star }$, and we have 
    \begin{equation*}
        \begin{aligned}
            & \left\lVert P _{p,2}^{\star }\theta \right\rVert_{2}^{}\gtrsim \frac{\sqrt{\epsilon }(k _{p,2}^{\star })^{\frac{1}{4}}}{N ^{\frac{1}{4}}}\sigma ,\\ 
            & \left\lVert P _{p,2}^{\star }\theta \right\rVert_{2}^{}\gtrsim \epsilon \sigma .
        \end{aligned}
    \end{equation*}
    Since $\epsilon \ge \frac{1}{\sqrt{N}}$, we also have 
    \begin{equation*}
        \left\lVert P _{p,2}^{\star }\theta \right\rVert_{2}^{}\gtrsim \frac{(k _{p,2}^{\star })^{1/4}}{\sqrt{N}}\sigma .
    \end{equation*}
    Similarly, we can regard $P _{p,2}^{\star }\theta $ as a vector in $\mathbb{R}^{k _{p,2}^{\star }}$ and define $k _{p,2,1}^{\star }$ and $k _{p,2,2}^{\star }$ accordingly as we did for $k _{1}^{\star }$ and $k _{2}^{\star }$. Since $k _{p,2}^{\star }\ge \max\limits \left\{k _{p,2,1}^{\star },k _{p,2,2}^{\star }\right\}$, we are again guaranteed to have a desired test by the inequalities above and the results in the $\ell _{2}$ norm case.

    Assume now $k _{p,2}^{\star }\neq k _{p,3}^{\star }$. Suppose we have $k _{p,2}^{\star }\ge k _{p,3}^{\star }+1$. By the monotonicity of Kolmogorov widths and the definitions of $k _{p,2}^{\star }$ and $k _{p,3}^{\star }$, we have 
    \begin{equation*}
        \frac{\sqrt{\epsilon }(k _{p,2}^{\star })^{\frac{1}{p}-\frac{1}{4}}}{N ^{\frac{1}{4}}}\sigma \le \tilde{D} _{p,k _{p,2}^{\star }-1}(K)\le \tilde{D} _{p,k _{p,3}^{\star }}(K)<\epsilon (k _{p,3}^{\star })^{\frac{1}{p}-\frac{1}{2}}\sigma ,
    \end{equation*}
    from which we obtain 
    \begin{equation*}
        \epsilon \gtrsim \frac{(k _{p,2}^{\star })^{\frac{2}{p}-\frac{1}{2}}}{(k _{p,3}^{\star })^{\frac{2}{p}-1}}\cdot \frac{1}{\sqrt{N}}\gtrsim \sqrt{\frac{k _{p,3}^{\star }}{N}}.
    \end{equation*}
    For $k _{p,3}^{\star }$, we already have $\left\lVert P _{p,3}^{\star }\theta \right\rVert_{2}^{}\gtrsim \epsilon \sigma $. By the inequalities above, it is not hard to verify that $\left\lVert P _{p,3}^{\star }\theta \right\rVert_{2}^{}\gtrsim \frac{\sqrt{\epsilon }(k _{p,3}^{\star })^{1/4}}{N ^{1/4}}\sigma $ and $\left\lVert P _{p,3}^{\star }\theta \right\rVert_{2}^{}\gtrsim \frac{(k _{p,3}^{\star })^{1/4}}{\sqrt{N}}\sigma $. Therefore, again by the discussion above and the results in the $\ell _{2}$ norm case, we are guaranteed to have a desired test in $\mathbb{R}^{k _{p,3}^{\star }}$.

    The last case is $k _{p,2}^{\star }\le k _{p,3}^{\star }-1$. The analysis is essentially the same as the case above. We have
    \begin{equation*}
        \epsilon (k _{p,3}^{\star })^{\frac{1}{p}-\frac{1}{2}}\sigma \le \tilde{D} _{p,k _{p,3}^{\star }-1}(K)\le \tilde{D} _{p,k _{p,2}^{\star }}(K)\le \frac{\sqrt{\epsilon }(k _{p,2}^{\star }+1)^{\frac{1}{p}-\frac{1}{4}}}{N ^{\frac{1}{4}}}\sigma ,
    \end{equation*}
    which leads to
    \begin{equation*}
        \epsilon \lesssim \frac{(k _{p,2}^{\star })^{\frac{2}{p}-\frac{1}{2}}}{\sqrt{N}(k _{p,3}^{\star })^{\frac{2}{p}}-1}\lesssim \sqrt{\frac{k _{p,2}^{\star }}{N}}.
    \end{equation*}
    Since we have $\left\lVert P _{p,2}^{\star }\theta \right\rVert_{2}^{}\gtrsim \frac{\sqrt{\epsilon }(k _{p,2}^{\star })^{\frac{1}{4}}}{N ^{\frac{1}{4}}}\sigma $, we also have $\left\lVert P _{p,2}^{\star }\theta \right\rVert_{2}^{}\gtrsim \frac{(k _{p,2}^{\star })^{1/4}}{\sqrt{N}}\sigma $ and $\left\lVert P _{p,2}^{\star }\theta \right\rVert_{2}^{}\gtrsim \epsilon \sigma $, which are sufficient for us to have a desired test in the space $\mathbb{R}^{k _{p,2}^{\star }}$.
\end{proof}

\subsection{Verification of the Minimality in Section \ref{subsection: setting of constraints}}\label{subsection: verification of the property of minimal}

\begin{proof}
    Without any loss of generality we assume $a _{i}$'s are sorted in decreasing order. If $k ^{\star }=d$, the assertion is obvious since $\left\lVert P ^{\star }\mu \right\rVert_{2}^{} = \left\lVert \mu \right\rVert_{2}^{}=c$. Suppose $k ^{\star }<d$ in the following. If $c ^{2}\le a _{k ^{\star }+1}$, then the vector
    \begin{equation*}
        \mu = (\underbrace{0,\dots,0}_{k ^{\star } \text{ zeros}},c,0,\dots,0)
    \end{equation*}
    satisfies $P ^{\star }\mu =\mathbf{0}$. Therefore, the assertion holds again. Let us now assume that $c ^{2}>a _{k ^{\star }+1}$. Since $\left\lVert \mu \right\rVert_{2}^{}=c$, we have $\sum\limits_{i=k ^{\star }+1}^{d}\mu _{i}^{2}=c ^{2}-\sum\limits_{i=1}^{k ^{\star }}\mu _{i}^{2}$. On the other hand, since $\sum\limits_{i=1}^{d}\frac{\mu _{i}^{2}}{a _{i}}\le 1$ and $a _{i}$'s are decreasing, we have
    \begin{equation*}
        \frac{c ^{2}- \sum\limits_{i=1}^{k ^{\star }+1} \mu _{i}^{2}}{a _{k ^{\star }+1}}\le \sum\limits_{i=k ^{\star }+1}^{d} \frac{\mu _{i}^{2}}{a _{i}} \le 1 - \sum\limits_{i=1}^{k ^{\star }} \frac{\mu _{i}^{2}}{a _{i}}.
    \end{equation*}
    Rearranging gives
    \begin{equation*}
        \frac{c ^{2}}{a _{k ^{\star }+1}}-1 \le \sum\limits_{i=1}^{k ^{\star }}\mu _{i}^{2}\left(\frac{1}{a _{k ^{\star }+1}}-\frac{1}{a _{i}}\right)\le \sum\limits_{i=1}^{k ^{\star }}\mu _{i}^{2}\left(\frac{1}{a _{k ^{\star }+1}}-\frac{1}{a _{1}}\right).
    \end{equation*}
    Hence, we conclude that $\left\lVert P ^{\star }\mu \right\rVert_{2}^{2}=\sum\limits_{i=1}^{k ^{\star }}\mu _{i}^{2}\ge \frac{c ^{2}/a _{k ^{\star }+1}-1}{1/a _{k ^{\star }+1}-1/a _{1}}$. It is clear that the vector given in Section \ref{subsection: setting of constraints} achieves this bound. The only remaining thing is to see that the vector indeed has $\ell _{2}$ norm $c$ and it belongs to $K$, which can be verified by simple algebra. The proof is completed.
\end{proof}

\subsection{Proof of Lemma \ref{lemma: testing estimation relationship}}\label{subsection: proof of the testing estimation relationship}
\begin{proof}
By the definition of $k _{1}^{\star }$ and $k _{e}^{\star }$ and the decreasing property of Kolmogorov widths, since $\frac{k ^{1/4}}{\sqrt{N}}\le \sqrt{\frac{k}{N}}$ for any $k, N$, we have $k _{1}^{\star }\ge k _{e}^{\star }$. On the other side, for any $j+1 \ge (k _{e}^{\star }+1)^{2}$, we have $D _{j}(K)\le D _{k _{e}^{\star }}(K)\le \sqrt{k _{e}^{\star }+1}\sigma \le (j+1)^{\frac{1}{4}}\sigma $. Therefore, we know $k _{1}^{\star } \lesssim (k _{e}^{\star })^{2}$. To finish the proof, we discuss different cases of $\epsilon $.

\textbf{(i): $0 \le \epsilon \le \frac{1}{\sqrt{N}}$.} In this case, the dominant terms in \eqref{eq: main theoretical lower bound} and \eqref{eq: lower bound of estimation} are $\frac{\sqrt{k _{1}^{\star }}}{N}\sigma ^{2}$ and $\frac{k _{e}^{\star }}{N}\sigma ^{2}$, respectively. Since $k _{1}^{\star }\lesssim (k _{e}^{\star })^{2}$, the claim follows.

\textbf{(ii): $\epsilon \ge \frac{1}{\sqrt{N}}$.} In this case, we need to further discuss the relationship between $k _{2}^{\star }$ and $k _{e}^{\star }$.

\textbf{(1): $k _{2}^{\star }=k _{e}^{\star }$.} When $\frac{1}{\sqrt{N}}\le \epsilon \le \sqrt{\frac{k _{2}^{\star }}{N}}$, the dominant terms in \eqref{eq: main theoretical lower bound} and \eqref{eq: lower bound of estimation} are $\epsilon \sqrt{\frac{k _{2}^{\star }}{N}}\sigma ^{2}$ and $\frac{k _{e}^{\star }}{N}\sigma ^{2}$, respectively. However, since $\epsilon \le \frac{k _{2}^{\star }}{N}$, we know $\epsilon \sqrt{\frac{k _{2}^{\star }}{N}}\sigma ^{2}\le \frac{k _{e}^{\star }}{N}\sigma ^{2}$. Otherwise, the dominant terms will be both $\epsilon ^{2}\sigma ^{2}$. Therefore, the claim holds for this case.

\textbf{(2): $k _{2}^{\star }+1 \le k _{e}^{\star }$.} We have
\begin{equation*}
    \sqrt{\frac{k _{e}^{\star }}{N}}\sigma \le D _{k _{e}^{\star }-1}\le D _{k _{2}^{\star }}\le \left(\frac{k _{2}^{\star }+1}{N}\right)^{\frac{1}{4}}\sqrt{\epsilon }\sigma ,
\end{equation*}
which leads to 
\begin{equation*}
    \epsilon \ge \frac{1}{\sqrt{N}}\cdot \frac{k _{e}^{\star }}{\sqrt{k _{2}^{\star }+1}}\ge \sqrt{\frac{k _{e}^{\star }}{N}}.
\end{equation*}
Once this holds, the dominant terms will be both $\epsilon ^{2}\sigma ^{2}$. Therefore, the claim holds for this case.

\textbf{(3): $k _{e}^{\star }+1 \le k _{2}^{\star }$.} We have 
\begin{equation*}
    \left(\frac{k _{2}^{\star }}{N}\right)^{\frac{1}{4}}\sqrt{\epsilon }\sigma \le D _{k _{2}^{\star }-1}\le D _{k _{e}^{\star }}\le \sqrt{\frac{k _{e}^{\star }+1}{N}}\sigma ,
\end{equation*}
which leads to 
\begin{equation*}
    \frac{1}{\sqrt{N}}\le \epsilon \le \frac{1}{\sqrt{N}}\cdot \frac{k _{e}^{\star }+1}{\sqrt{k _{2}^{\star }}}\lesssim \sqrt{\frac{k _{e}^{\star }}{N}}.
\end{equation*}
Similarly, this implies that the dominant terms in \eqref{eq: main theoretical lower bound} and \eqref{eq: lower bound of estimation} are $\epsilon \sqrt{\frac{k _{2}^{\star }}{N}}\sigma ^{2}$ and $\frac{k _{e}^{\star }}{N}\sigma ^{2}$, respectively, and $\epsilon ^{2}\sigma ^{2}$ will never dominate. Since $\epsilon \lesssim \frac{1}{\sqrt{N}}\cdot \frac{k _{e}^{\star }}{\sqrt{k _{2}^{\star }}}$ by the inequalities above, we know $\epsilon \sqrt{\frac{k _{2}^{\star }}{N}}\sigma ^{2}\lesssim \frac{k _{e}^{\star }}{N}\sigma ^{2}$.

The discussion above shows that the minimax rate of estimation is greater than the minimax rate of testing in all possible cases. The proof is completed.
\end{proof}

\end{document}